%% file: Xue-Mei-Li.Thesis-arxiv.tex
	\documentclass[12pt]{report}
\usepackage{amsfonts,amsmath,amssymb, amsthm, times}
	\makeindex
	
\title{Stochastic flows on non-compact manifolds}
\author{Xue-Mei Li\\
{\small A thesis submitted for the degree of Doctor of Philosophy} \\
{\small  Mathematics Institute, The University of Warwick}}
\usepackage{imakeidx}
\makeindex
\date{1992}

\newcommand{\A}{{\bf \cal A}}
\newcommand{\B}{{ \bf \cal B }}
\newcommand{\C}{{\cal C}}
\newcommand{\E}{{\mathbb E}}
\newcommand{\F}{{\cal F}}

\newcommand{\h}{{\cal H}}

\newcommand{\half}{{  {1\over 2}  }}
\newcommand{\heatsemif}{{ {\rm e}^{ \half t\triangle^{h,1}}   }}
\newcommand{\heatsemi}{{ {\rm e}^{\half t \triangle^{h}}   }}

\newtheorem{theorem}{Theorem}[section]
\newtheorem{proposition}[theorem]{Proposition}
\newtheorem{lemma}[theorem]{Lemma}
\newtheorem{corollary}[theorem]{Corollary}
\newtheorem{definition}{Definition}[section]

\def\limsup{\mathop{\overline{\rm lim}}}
\def\liminf{\mathop{\underline{\rm lim}}}
\def\exp{{\rm e}}

\begin{document}
\maketitle
	\begin{titlepage}
\begin{center}
{\bf \large  Introduction to the  Arxiv version}
\end{center}

I have finally allowed myself persuaded  to make my PhD thesis available on the Arxiv, as some results are still relevant.
In particular, Chapters 2, 3, and  6 were not submitted for publication anywhere.  The thesis was completed at the University of
Warwick in 1992 under the supervision of Professor KD Elworthy, to whom I own my gratitude.  My study in the U.K. would not have been possible without the scholarship from the Sino-British Friendship Scholarship Scheme which fully funded my study.

One of the contributions of the thesis is to solve the open problem on the
so called strong completeness (the existence of a globally defined continuous solution
flow) for  stochastic differential equations (SDEs) on non-compact manifolds. Such results are well known for compact manifolds and for
the Euclidean space with Lipschitz continuous coefficients. At the same time we introduce
a new method: the derivative flow method which allows us to improve the results in 
$R^n$, not simply by weakening the Lipschitz conditions, but more profoundly. The idea is to take the holistic approach
and take the pair of stochastic processes: the solution flow and its derivative
with respect to the initial data (differentiated in the sense of probability only). This allows 
significant growth of the derivatives  (e.g. in fraction powers and even linearly) of the coefficients of the SDE
provided  the growth of the coefficients themselves are controlled. In light of the developments after the PhD thesis,
we mention also the martingale method
for obtaining integral representation formula for the weighted heat semigroup and the intertwining of heat semigroups.

An immediate and very well known demonstration of the connection between diffusion processes and 
geometry  is  given by the non-explosion of a Brownian
 motion. A $h$-Brownian motion on $M$  is a strong Markov process with generator $\Delta^h=\Delta+2 \nabla h$
 which exists naturally on a  Riemannian manifold. If $h=0$ this is the Brownian motion. 
 A Brownian motion has no-explosion from anywhere  if and only if 
the fundamental solution of the heat equations  is unique.  This is  related also to   the
 behaviour at infinity of the associated probability semigroup, c.f. Chapter 4. 
 We shall make more connections  between  properties of diffusions  and  that of the
manifolds.   For example,  a slightly modification of Bakry's method \cite{BA86} shows that
 if  the heat semigroup for 1-forms  is   
``continuous" on $L^\infty$, then the Brownian motion has no 
explosion. 
Indeed one can choose an SDE whose solutions $F_x(x)$  is an $h$-BM with initial value $x$.
We can take the special gradient SDE obtained with the isometric embedding of the manifold into $R^m$,
in most cases we could choose to use other equations as well.
The symmetric semi-group  on function is defined as usual: $P_tf(x)=\E(f(F_t(x)\chi_{\{t<\xi(x))\}})$, where $\xi(x)$ is the life time of the solution $F_x(x)$.
A differential 1-form is a section of the cotangent bundle $ \Gamma T^*M$. The semi-group on differential 1-forms is: for $v\in T_xM$,
$(\delta P_t)(\phi)(v)=\E\left(\phi(T_xF_t(v))\chi_{t<\xi(x)}\right)$.  If  applied to $\phi=df$ where $f:M\to R$ is a function ,
the inter-twinnining property 
$$(\delta P_t) (df)=d(P_tf)$$
holds  and if 
 there is a point $x_0$ and  an open neighbourhood $U$ of  $x_0$ such that  $|\heatsemif df|_x\le C_{t}(x)|df|_\infty$ 
 for any  $f\in C_K^\infty$ and for any $t\le t_0$  and any  $x\in U$,
we expect non-explosion. See Lemma \ref{le: homotopy 1}.  Let
$W_t(v)=\E\{T_xF_t(v) |\F_t^x\}$, where $\F_t^x$ is the $\sigma$-filtration generated by $F_t(x)$.
Since
 $\heatsemif df=\E [df(W_t)\chi_{t<\xi(x)}]$ and $W_t$ satisfies
 $$\frac d {dt}|W_t|^2 \le -\mathrm {Ric}_{F_t(x_0)}(W_t, W_t)+\nabla^2h(W_t, W_t),$$
 this makes connection with the Bakry-Emery condition, i.e. the deformed Ricci curvature is bounded from below. 
 We shall see that the Bakry-Emery condition,  is not  optimal for either the explosion problem or for the existence of an invariant probability measure. For the existence of an invariant probability measure see  the article `On extension of Myers' Theorem' (1995) for further clarification.

We observe that integrability conditions on the derivative flows
are also related to properties of the manifolds as well as on the diffusion 
themselves. Assuming this type of condition, we  have results on non-explosion problem, 
the existence of global smooth solution flow (strong p-completeness), as well as homotopy and
cohomology vanishing results, and integral representations of semigroups.

\medskip

Before going into detail, we would also like to point out that the existence
of a continuous flow is of basic importance in the study of ergodic properties of dynamical systems.
Our result on the strong completeness implies  the existence of perfect co-cycles.

\medskip

The first three chapters of the thesis are preliminary. Of these chapter 1 is to
establish notation, chapter 2 and chapter 3 deal 
with results which are essentially  known (especially in special cases) but not
 found anywhere in suitable form. A new result is the essential self-adjointness of 
 the Bismut-Witten Laplacian  $\triangle^h$, and $d+\delta^h$, on non-compact manifolds.

\medskip
Inspired by the idea of uniform covers, we introduce the weak uniform cover
method, in chapter 4, which allows us to look at many different problems from the
same point of view including the explosion problem, the behaviour at infinity 
of diffusion processes and  probability semigroups, along with the geometry
of the underlying manifold at infinity. As an example we conclude  that linear growth
 gives the $C_0$ property. We have a result on the non-explosion of SDEs on an
 open set of $R^n$. The results  on non-explosion are used in chapter 5 to obtain 
strong completeness. 

\medskip

In chapter 5,   we solve a question in SDE open for over twenty years:
Does a SDE has a smooth global solution flow on a non-compact manifold if the SDE does not explode?
The only  results, known at this points, are for linear spaces and for smooth SDEs on compact manifolds, 
the first assume globally Lipshcitz type  condition and latter  is obtained
by lifting the SDE to the diffeomorphism group. The first counter example of an SDE for which solutions from any point exist for all time
but for which the continuous dependence on the initial point fails, was given by Elworthy in 1978, prior to which there was the popular belief that  if the there is no explosion from any points then the solution
depends continuously on the initial data (provided coefficients are sufficiently smooth.
This incorrect belief was led by theories of ODE where there is no noise and theory of one dimensional SDEs.
In a recent article `Lack of strong completeness for stochastic flows'  with Schetzow we constructed an
SDE in $R^n$ with bounded smooth coefficients, the SDE is smooth but the continuous dependence on 
the initial data failed.

 Our results  also improve  and offer insights to the flow problem on the Euclidean space, beaking the
barrier of global Lipschitz continuity. We may allow polynomial growth of the derivatives of the coefficients of the SDE
 in some cases. 
Our observation is that  the existence of a global smooth solution flow is 
linked in a fundamental manner with  the derivative flows
 and with the rate of growth of the data in the linearised SDE along its solutions.

The concept of strong p-completeness is introduced due to the 
complexity of strong completeness.   A stochastic differential equation is said 
to be strongly p-complete if it has a version which is smooth on any given smooth 
singular p-simplex.  It is called strongly complete if there is a smooth flow
 of the solution on the whole manifold.   We start with a theorem on strong
 1-completeness which is applied to obtain the intertwing property $dP_tf=\delta P_t(df)$ in chapter 6.
%Here $\delta P_t(df)=Edf(TF_t\chi_{t<\xi})$, and $TF_t$ is the derivative 
%flow (c.f. section 1.2), a solution to the linearised SDE. 
 The intertwining property is an underlying condition for
our martingale method to differentiation formula of Markovian semi-groups.

Strong 1-completeness is needed in chapter 7 to obtain homotopy vanishing
results and is also used in chapter 8. Later we give a theorem on 
strong p-completeness, in particular the existence of a continuous flow
and flows of diffeomorphisms.  
 Strong p-completeness leads naturally  to   cohomology vanishing  given
 strong moment stability (see page ~\pageref{th: p-complete cohomology}). 
 These theorems are originally given in terms of integrability  conditions on 
the derivative flows $TF_t$, but those conditions can also be checked in terms of the  bounds on 
the coefficients of the stochastic differential equations.  See 
section ~\ref{applications}. In particular, as a simple example,   there is a 
smooth  flow of  Brownian motions on a submanifold of $R^n$  whose second 
fundamental form  is bounded.

\medskip

In chapter 6, we look at the probability semigroup for 1-forms, the heat
 semigroup for  1-forms, and the differentiation of the probability semigroups
 for functions and ask  a basic and yet important question: does $dP_tf$ equals
 $\delta P_t(df)$? This question is answered from different approaches, but all
  require conditions on $TF_t$.  We use this knowledge of the semigroups 
extensively in the next two chapters. Following from the discussions, we get
a result on the $L^p$ boundedness and contractivity of heat semigroups for forms. 
As a corollary, 
we give a cohomology vanishing result.

\medskip

Proceeding to chapter 7, the interplay  between diffusion processes and 
geometric and  topology properties of the underlying manifold becomes clearer.
The main theorems here are  on the vanishing of the first  homotopy group
$\pi_1(M)$   given conditions on the regularity of the diffusions and strong
 moment stability, i.e. the supremum over a compact set of the $p^{th}$ norm of
the derivative flow decays exponentially fast.  In particular we conclude 
that a certain class of manifold cannot have a strongly  moment stable s.d.s.
if its fundamental group is not trivial. We study, again briefly,  the vanishing of harmonic 1-forms and   cohomology vanishing.
The moment estimates was extended, after the completion of the thesis, 
to tensor field of $T_xF_t$, see the article of Elworthy and Rosenberg `Homotopy and Homology Vanishing Theorems and the Stability of Stochastic Flows'.

\medskip

Finally in chapter 8, we obtain a formula for $d(P_tf)$, for elliptic Markov generators, 
 in terms of $f$ and a martingale, but does not depend on $df$. For this we tweaked the basic idea, that  showing $P_tf$ is a
solution of the Kolmogorov equation,   with a simple martingale technique. 
 I happened on this method ` while  trying to extend a result of Elworthy
  \cite{ELflow}  for closed differential 1-forms to closed differential $q$-forms with $q>1$ and at the same time trying to understand Malliavin Calculus.
 At the line 4 of the proof,  one would normally take the expectation of the both sides of the identity and conclude Kolmogorov's equation. 
Instead, we multiplied both sides by a well chosen martingale, and applied the It\^o isometry. We then explored the semi-group on differential forms, especially its inter-twining property and apply the Markov property of the system of SDE on the tangent bundle. 
An easy application of this  formula is an implicit formula for the gradient of the logarithm of the probability kernels. For a compact manifold, the gradient logarithmic heat kernel  formula was proved by Bismut  using a perturbation of the driving  Brownian motion on the horizontal SDE on the orthonormal frame bundle.
The method introduced in this thesis has since been used in proving the Strong Feller property for stochastic partial differential equations (SPDE) and 
for a while it was the only method for proving uniqueness of invariant probability measures for SPDEs, the first  article showing strong Feller property by a differentiation formula for the semigroup generated by a SPDE was given by
Da Prato, Elworthy, and Zabczyk. More recently, it was shown by Hairer and Mattingly that an asymptotic Feller property can be used to replace the strong Feller property. The former can hold under less restrictions on the equation can can be obtained by an approximate or modified differentiation formula.
 
  We also obtain an integral representation  formula for  closed differential forms of degree $q$  where $q>1$,
  extending the result of Elworthy \cite{ELflow} where such a formula was obtained for closed differential 1-forms. In fact this was the motivation for the earlier work studying $d(P_tf)$, as explained earlier the intertwining property will imply that  $d(P_tf)=e^{\frac 12 \Delta^{h,a}}(df)$, the semigroup acting on $df$.
  For $q>1$, we took the geometric advantage  afforded by gradient SDEs. Gradient SDEs is not needed for differential 1-forms.
In addition, we develop also the martingale method for proving a formula for $P_t\phi$ where $\phi$ is a $q$-form, $q\ge 1$. This formula turned out to be a special version of the for-mentioend integral formula in the case  the differential form is exact.
Our use of  the gradient stochastic differential equation for the treatment of semi-groups on differential
forms of degree higher than $1$ motivated, with hindsight, the study of intuitive Riemmann connections associated to, elliptic and possibly degenerate,  stochastic differential equations of non-gradient form. See `Concerning the Geometry of Diffusion Operators' by Elworthy, LeJan, and nyself.
 Another theme  is to  relate the geometrical and
topological properties of  manifolds to that of  the diffusion processes
on it. 

Due to the problem of latex compatibility, some figures are missing in this version, 
 a comment or two are added.

\pagebreak

\begin{center}
Acknowledgement
\end{center}

I would like to thank  my supervisor professor D. Elworthy for leading me into
 the subject. I owe my interest in mathematics to his teaching and influence. 
I am grateful to professor Zhankan Nie of Xi'an Jiaotong university, and would like acknowledge
with thanks the full financial support  from  the Sino-British Friendship Scholarship Scheme during my entire 
period of study for the degree of Master of Science and for the degree of doctor of Philosophy. I would like to thank 
 professor J. Zabczyk for encouragement and many friends and colleagues for
 useful discussions and general help.  I also benefited from the EC programme 
SCI-0062-C(EDB).
 
	\end{titlepage}

\begin{abstract}

Here we look at the existence of solution flows of stochastic differential equations
on noncompact manifolds and the properties of the solutions in terms of the 
geometry and topology of the underlying manifold itself. We obtain some 
results  on  ``strong p-completeness" given conditions on the derivative
 flow, and thus given suitable conditions on the coefficients of the stochastic
 differential equations. In particular a smooth flow of  Brownian motion 
exists on submanifolds of $R^n$ whose second fundamental forms are bounded.
 Another class
of  results we obtain is on homotopy  vanishing given strong moment 
stability. We also have results on obstructions to moment stability by 
cohomology. Also we obtain formulae for $d(\heatsemi \phi)$ for differential
 form $\phi$ in terms of a martingale and the form itself, not just its
 derivative, extending Bismut's formula.
\end{abstract}

		\tableofcontents
	\addtocontents{chap1.toc}{chap.add.toc}

\part{Background}

\chapter{Introduction  and preliminaries}
\section{Introduction}

 Let $\{\Omega,{\cal F},{\cal F}_t, P\}$ be  a filtered probability  space satisfying the usual conditions including right continuity. Let $M$ be a $n$ dimensional smooth manifold. Consider the following stochastic differential 
equation(s.d.e.)\index{s.d.e.} on $M$: 

	\begin{equation}
dx_t=X(x_t)\circ dB_t+A(x_t)dt.   \label{eq: basic}
\end{equation}

\noindent
Here  $B_t$ is a $R^m$ valued Brownian motion($\F_t$-adapted),   $X$ is $C^2$
 from $ R^m\times M$ to
 the tangent bundle $TM$  with $X(x)$: $R^m$ $\to T_xM$ a linear map for each $x$ in 
$M$,  and $A$ is  a $C^1$ vector field on $M$.

By  $\circ$ we mean  the Stratonovich integral. 

\bigskip

  Let  $u$  be a random variable independent of $\F_0$.  Denote by  $F_t(u)$ the
 {\it solution}  starting from $u$, with  $\xi(u)$ the {\it explosion time}. 
By a solution we  mean  a maximal solution which is sample continuous unless 
otherwise stated.  Under our assumptions on coefficients the solution to equation 
$~\ref{eq: basic}$ is unique in the sense of that if $(x_t, \xi)$ and $(y_t, \eta)$
 are two solutions with same initial point, then they are versions of one another, 
i.e.   $\xi=\eta$ almostly surely, and   $x_t=y_t$ almostly surely for 
all $t$ on  $t<\xi$.  

Furthermore   the solutions to the stochastic differential equation 
$(~\ref{eq: basic})$ are 
{\it  diffusion}\index{diffusion} processes, i.e.  $\{F_t(x), t\ge 0\}$ is a path
 continuous  strong Markov process for each $x$. In fact most of the interesting 
diffusions can be given in this way. 
The importance of a diffusion process is largely related to its semigroup  
 $P_t$ (given by  
 $P_tf(x)= Ef(F_t(x))\chi_{t<\xi(x)}$), and the  corresponding infinitesimal 
 generator $\A$ (which we discuss in detail in chapter 4),  or to the consideration
of $(~\ref{eq: basic})$ as a random perturbation of the dynamical system given by
 the vector field $A$.

The pair $((X,A), (B_t,t))$, often simplified as $(X,A)$, is called a 
{\it stochastic dynamical system}\index{stochastic dynamical system}\index{s.d.s.}
(s.d.s.).

\bigskip

The thesis contains two themes. The first considers fundamental problems
of stochastic differential equations: completeness and strong
 completeness.  The main theme of the second part is to  relate geometrical and
topological properties of  manifolds to the diffusion processes
on it. The first  furnishes a start for the second. And both  make 
good use of the fact that integrability conditions on the derivative flows
give very strong conditions on the manifolds as well as on the diffusion 
themselves. We have nonexplosion, strong p-completeness results, homotopy and
cohomology vanishing results assuming this type of condition. 

As we shall see  properties of diffusions are directly related to  those of the
manifolds.  
An immediate demonstration of the connection between diffusion processes and 
geometry  is  given by the well known fact: the nonexplosion of a Brownian
 motion (which exists naturally on a  Riemannian manifold) is equivalent to 
the uniqueness  of solutions of the heat equations, which  relates to the
 behaviour at infinity of the associated probability semigroup as shown in 
chapter 4.

\bigskip

Before going into detail, we would also like to point out that the existence
of a continuous flow is of basic importance in the study of ergodic properties.

\bigskip

The first three chapters of the thesis are preliminary. Of these chapter 1 is to
establish notation and contains quoted results, chapter 2 and chapter 3 deal 
with results which are essentially  known (especially in special cases) but not
 found anywhere in suitable form. Much of the treatment in chapter 2, in which
we  discuss carefully the Bismut-Witten Laplacian  $\triangle^h$ and its associated
semigroups, was arrived at together with D. Elworthy and S. Rosenberg. In chapter
3 we look at the invariant measure for an h-Brownian motion and its ergodic
 properties.

\bigskip
Inspired by the idea of uniform covers, we introduce the weak uniform cover
method in chapter 4, which allows us to look at many different problems from the
same point of view including the explosion problem, the behaviour at infinity 
of diffusion processes and  probability semigroups, along with the geometry
of the underlying manifold at infinity. As an example we conclude  that linear growth
 gives the $C_0$ property and have a result on the nonexplosion of s.d.e. on an
 open set of $R^n$. The results  on nonexplosion are used in chapter 5 to obtain 
strong completeness.

\bigskip

In chapter 5,  the concept of strong p-completeness is introduced due to the 
complexity of strong completeness. A stochastic differential equation is said 
to be strongly p-complete if it has a version which is smooth on any given smooth 
singular p-simplex.  It is called strongly complete if there is a smooth flow
 of the solution on the whole manifold.   We start with a theorem on strong
 1-completeness which is applied to get $dP_tf=\delta P_t(df)$ in chapter 6.
Here $\delta P_t(df)=Edf(TF_t\chi_{t<\xi})$, and $TF_t$ is the derivative 
flow(c.f. section 1.2). 
Strong 1-completeness is needed in chapter 7 to get homotopy vanishing
results and is also used in chapter 8. Later we give a theorem on 
strong p-completeness, in particular the existence of a continuous flow
and flows of diffeomorphisms.  
 Strong p-completeness leads naturally  to   cohomology vanishing  given
 strong moment stability (see page ~\pageref{th: p-complete cohomology}). 
 These theorems are originally given in terms of integrability  conditions on 
$TF_t$, but those conditions can also be checked in terms of the  bounds on 
the coefficients of the stochastic differential equations.  See 
section ~\ref{applications}. In particular, as a simple example,   there is a 
smooth  flow of  Brownian motions on a submanifold of $R^n$  whose second 
fundamental form  is bounded.

\bigskip

In chapter 6, we look at the probability semigroup for 1-forms, the heat
 semigroup for  1-forms, and the differentiation of the probability semigroups
 for functions and ask  a basic and yet important question: does $dP_tf$ equals
 $\delta P_t(df)$? This question is answered from different approaches, but all
  require conditions on $TF_t$.  We use this knowledge of the semigroups 
extensively in the next two chapters. Following from the discussions, we get
a result on the $L^p$ boundedness and contractivity of heat semigroups for forms. 
As a corollary, 
we give a cohomology vanishing result.

\bigskip

Proceeding to chapter 7, the interplay  between diffusion processes and 
geometric and  topology properties of the underlying manifold becomes clearer.
The main theorems here are  on the vanishing of the first  homotopy group
$\pi_1(M)$   given conditions on the regularity of the diffusions and strong
 moment stability, i.e. the supremum over a compact set of the $p^{th}$ norm of
the derivative flow decays exponentially fast.  In particular we conclude 
that a certain class of manifold cannot have a strongly  moment stable s.d.s.
if its fundamental group is not trivial. We also 
look briefly at the vanishing of harmonic 1-forms and   cohomology vanishing.

\bigskip

Finally in chapter 8, we get a formula for $d(P_tf)$ for elliptic systems in terms
of $f(x_t)$ and a martingale,  extending the formula given in \cite{ELflow}. We
also obtain a   similar type of formula for  $q$ forms. 
In particular we have an explicit formula
for the gradient of the logarithm of the heat kernel extending Bismut's formula.

\bigskip
Bibliographical notes are scattered in every chapter.  But we would like to
 refer to \cite{ELbook} and  \cite{EL-survey} for general references. See 
also \cite{IK-WA}.

\vfill

\input{picture.1}

\section{Probabilistic set-up}

\subsection*{Brownian motions and Brownian systems}

  A path continuous strong Markov process on a Riemannian manifold is called a 
{\it Brownian motion}\index{Brownian motion} (BM\index{BM}) if its  generator  is ${1\over 2 } \triangle$, where $\triangle$ denotes the  Laplace-Beltrami operator. It is called a Brownian motion with drift $Z$ if its  generator is ${1\over 2 } \triangle +Z$. Here $Z$ is a $C^\infty$ smooth vector field on $M$.  The drift $Z$ is called a {\it gradient drift}\index{gradient drift} if $Z=\nabla h$ for some  function $h:  M \to R$. We will always assume $h$ is 
$C^\infty$ smooth for simplicity. But in many cases we only need it to be $C^2$.

A stochastic dynamical system $(X,A)$ is called a 
{\it Brownian system}\index{Brownian system} (with drift $Z$) if its solution is
 a Brownian motion  (with drift $Z$). Equivalently  $X(x): R^m\to T_xM$ is
 an orthogonal projection for 
each $x\in M$. i.e. $X^*X=Id$ and 
$Z=A^X\index{$A^X$}=: A+\half {\rm trace}\nabla X(X(\cdot))(\cdot)$. 

In the text, we often use {\it h-Brownian system} for a Brownian system with drift 
$\nabla h$, and   {\it h-Brownian motion}\index{$h$-Brownian motion} for a Brownian motion with drift $\nabla h$.

\subsection*{Gradient Brownian flow}

 Let $f: M\to R^m$ be an isometric immersion. Define
$X(x): R^m\to T_xM$ as follows: 
$$X(x)(e)=\nabla <f(\cdot),e>(x), \hskip 6pt  e\in R^m.$$
 A stochastic dynamical system $(X, A)$ with $X$ so defined is called a {\it gradient Brownian system}\index{gradient Brownian system} (with drift).  When $A=0$, its 
solution flow  is called a {\it  gradient Brownian flow}. It is called a 
{\it h-gradient Brownian system}\index{h-gradient Brownian system}
 if $A=\nabla h$. For such a system, we can always 
choose, for each $x\in M$, an  orthonormal basis (o.n.b.) $\{e_1, \dots, e_m\}$ of $R^m$, such that  for all $v\in T_xM$ and $i=1,  \dots, m.$ either
$$\nabla X(v)(e_i)=0  $$
or $$X(x)(e_i)=0.$$ 
Here $\nabla X(-)(e_i)$ denotes for the covariant derivative of $X(e_i)$. 
  In particular we have 
\begin{equation}
\sum_{i=1}^m \nabla X^i(X^i)=0.
\label{eq: gradient}
\end{equation}
 Here $X^i(x)=X(x)(e_i)$. See \cite{ELbook} or \cite{ELflour} for the proof.

\section*{The derivative flow}

Let $M_t(x)=\{\omega: t<\xi(x,\omega)\}$.  Denote by $TF_t$ the solution flow to the following covariant  stochastic differential equation:
\begin{equation} 
\label{coveq}
Dv_t=\nabla X(v_t)\circ dB_t +\nabla A(v_t)dt.
\end{equation}
It is in fact the derivative\index{derivative flow} of $F_t$ in measure in the 
following sense as shown in \cite{ELbook}:
let $f: M\to R$ be a $C^1$ map and  $\sigma \colon [0, 1]\to M$   a smooth curve 
with $\sigma(0)=x$, and $\dot{\sigma}(0)=v$.   Then for any $\delta>0$:

$$\lim_{r\to 0}P\{\omega\in M_t(x): 
|{f(F_t(\sigma(r))-f(F_t(x))\over r} -df(TF_t(v))|>\delta\}=0.$$

\noindent    By convention  $f(F_t(\sigma(r)),\omega)=0$  if
$\omega\not \in M_t(\sigma(r))$.          Let $x_0\in M$, and let $ v_0\in T_{x_0}M$. We write $x_t=F_t(x_0)$\index{$x_t$}, $v_t=TF_t(v_0)$\index{$v_t$} for simplicity. Clearly $TF_\cdot$ is a solution to the stochastic differential equation  on the tangent bundle $TM$ corresponding to  equation ~\ref{coveq}. Furthermore it has the same explosion time as $x_\cdot$ according to 
\cite{ELflow}. In a trivialisation we may not distinguish between    the derivative 
flow and its principal part, if there is no confusion caused.

\bigskip

\section*{It\^o  formula for forms}
 
First, we recall {\it It\^o  formula for one forms}\index{Ito's formula}  from \cite{ELflow}:
Let $T$ be a stopping time with $T<\xi$, then 
\begin{equation}\begin{array}{cl}
\phi (v_{t\wedge T})
 =&\phi(v_0)+\int_0^{t\wedge T} \nabla \phi\left (X(x_s) dB_s\right)(v_s)  \\
  &+ \int_0^{t\wedge T}\phi\left(\nabla X(v_s)  dB_s \right)
 +\int_0^{t\wedge T}{\cal L}\phi(v_s)\, ds.
\end{array}
\end{equation}

Here $\cal L$ is  the differential operator on 1-forms associated with a s.d.s. $(X,A)$ defined as follows (for Levi-Civita connection):

 \begin{equation}\begin{array}{cl}
({\cal L} \phi)_x(v)
=&{1\over 2}{\rm trace} \nabla^2\phi(X(x)(\cdot), X(x)(\cdot))(v) \\
 &+{1\over 2}\phi({\rm trace} R(X(x)(\cdot), v)X(x)(\cdot) )    \\
&+L_{A^X} +{\rm trace} \nabla\phi(X(x)(\cdot))\nabla X(v)(\cdot).
\end{array}
\label{sy: ell}
\end{equation}

If $\A={1\over 2}\triangle +L_Z$ and $\phi$ is closed, then as shown in 
\cite{ELflow}: 
$${\cal L}\phi={1\over 2}\triangle^1 +L_Z.$$

\bigskip

For higher order forms, there is  a similar formula from \cite{ELflow}; we
 quote here the formula for gradient systems. But first we need some
 notation (see \cite{AB-MA}):

\bigskip

\noindent
Let  $A$ be a linear map from  a vector space $E$ to $E$.  We define $(d\Lambda)^qA$ from $E\times \dots \times E$ to $E\times \dots \times E$ as follows:

$$(d\Lambda)^q A(v^1,\dots, v^q)
=\sum_{j=1}^q (v^1, \dots, Av^j, \dots, v^q).$$

\bigskip

\noindent
 Let $v_0=(v^1_0, \dots, v_0^q)$, where  $v_0^i\in T_{x_0}M$.  Denote by $v_t$
the $q$ vector induced by $TF_t$: 
 $$v_t=(TF_t(v_0^1), TF_t(v_0^2),\dots, TF_t(v_0^q)).$$ 
 
\bigskip

Here is   {\it It\^o  formula for gradient Brownian systems  for $q$ forms}\index{Ito formula} as given in \cite{ELflow}:

  \begin{equation} \begin{array}{cl}
\label{eq: Ito formula for gradient}
\psi(v_t)=&\psi(v_0)+\int_0^t\nabla \psi(X(x_s)dB_s)(v_s)\\
&+\int_0^t\psi\left((d\Lambda) ^q (\nabla X(\cdot )(dB_s)(v_s)\right)
+\int_0^t \frac 1 2 \triangle ^{h,q}(\psi)(v_s)ds.
\end{array}
\end{equation}

\bigskip

\noindent
{\bf Moment exponents:} For $p\in R$, $K\subset M$ compact, 
 we have {\it moment exponents}\index{moment exponents} $\mu_K(p)$:

$$\mu_K(p)=\limsup_{t\to \infty} \sup_{x\in K}{1\over t}\log E|T_xF_t|^p$$
and point moment exponents $\mu_x(p)$:
$$\mu_x(p)=\limsup_{t\to \infty} {1\over t}\log E|T_xF_t|^p.$$
We will say a flow\label{moment stable} is {\it moment stable}\index{moment stable} if $\mu_x(1)<0$ 
for each x,  and {\it strongly moment stable} if $\mu_K(1)<0$ for each compact 
set $K$. It is called {\it  strongly $p^{th}$-moment stable}\index{strong $p^{th}$-moment stability} 
if $\mu_K(p)<0$ for each compact subset $K$ of $M$. For discussions
on various exponents and related problems, see Arnold\cite{Arnold83}, 
Baxendale and Stroock \cite{BA-ST}, Carverhill, Chappel and Elworthy 
\cite{ca.ch.el},  Chappel\cite{Chappell}, and Elworthy \cite{ELflow} 
\cite{ELflour}.

\bigskip

\noindent
{\bf Invariant measure:} A $\sigma$-finite measure $m$ on $M$ is called an invariant measure\index{invariant measure} for $F_\cdot$, if  the following holds for all $t>0$
and   $L^1$  functions $f$:
$$\int_M P_tf(x)m(dx)= \int_M f(x)m(dx).$$

\bigskip

\noindent{\bf Some notation:}
Let $dx$ be the Riemannian volume element\index{$dx$} on $M$. Then $M$ is said to
 have finite volume if  $\hbox{Vol}(M)=\int_M dx<\infty$. Similarly $M$ is said to 
have {\it finite  h-volume}\index{h-volume} if $\hbox{h-Vol}(M)=
\int_M e^{2h}\, dx<\infty$. 
e.g. $R^n$ has finite volume  for ${1\over (2\pi)^{n/2}}e^{2h}$ the Gaussian
 density given by   $h(x)=-{1\over 4} |x|^2$.

\section{Unbounded operators on Hilbert space}

Let $\h$ be a Hilbert space. An operator $T$ on $\h$ is a linear map from a subspace of $\h$ to $\h$. The subspace is called its domain, which we denote by $D(T)$. We also denote by ker$(T)$ the kernel of $T$, and $Im(T)$ ( or Ran$(T)$) the image of $T$ in $\h$. A  {\it closed} operator $T$ is an operator with its graph 
$\Gamma(T)=\{(\phi, T\phi): \phi\in D(T)\}$ a closed subspace of $\h\times\h$. An operator $S$ is called an $extension$ of $T$ if $\Gamma(S)\supset\Gamma(T)$, i.e. $D(S)\supset D(T)$ and $S=T$ on $D(T)$. If $T$ has a closed extension, it is 
 {\it closable}. In this case the smallest extension of it is called its 
{\it closure}, which we denote as $\bar T$. 

\begin{proposition}\cite{RE-SI80}[p.250]. If $T$ is closable, 
$\Gamma(\bar T)=\overline{\Gamma(T)}$, i.e.
$D(\bar T)$ contains precisely those $\phi$ such that: there exist
 $ \{\phi_n\}\subset D(T) $ with 
$\phi_n\to \phi$ in $\h$ and $T\phi_n\to \eta,$ in $\h$, some  $\eta \in \h$
(and then $\bar T\phi=\eta$).
\end{proposition}

An interesting class of operators are operators with dense domain, called densely defined operators. Let $T$ be such an operator, there is a well defined  adjoint operator $T^*$ defined by:
$$D(T^*)=\{\phi\in \h: <T\psi, \phi>=<\psi, \eta>,\hskip 4pt\hbox{some}\hskip 4pt \eta \in \h,
\hskip 4pt \hbox{all} \hskip 4pt \psi \in D(T)\}.$$
For such $\phi$ and $\eta$ we set: $T^*\phi=\eta$. 

Notice that if $S\subset T$, then $T^*\supset S^*$. If $T\subset T^*$, $T$ is called a {\it symmetric operator}. There are also the following properties of adjoint operators:

\begin{theorem}\cite{RE-SI80}[p.253]. $\label{th: properties of adjoint}$
Assume $D(T)$ is dense in $\h$, then:
\begin{enumerate}
\item 
$T^*$ is closed,  and $T^{**}$ is symmetric.
\item 
$T$ is closable if and only if $D(T^*)$ is dense, in which case $\bar T=T^{**}$.
\item If $T$ is closable , then $(\bar T)^*=T^*$.
\end{enumerate}
\end{theorem}

 A symmetric operator is
$self$-$adjoint$ if $T=T^*$. It is  {\it essentially self-adjont} if its closure is self-adjoint, equivalently if it has only one self-adjoint extension.

Here is another  test for essential self-adjointness:

\begin{theorem}\cite{RE-SI80} [p.257].
 Let $T$ be a symmetric operator on a Hilbert space. Then the following are equivalent:
\begin{enumerate}
\item
$T$ is essentially self-adjoint.
\item
ker$(T^*\pm i)=\{0\}$.
\item
Ran$(T\pm i)$ are dense.
\end{enumerate}
\end{theorem}

\noindent
{\bf The Friedrichs extension}: Let $T$ be a positive symmetric operator. Let 
${\varepsilon}$ be its quadratic form (bilinear form) defined as follows: 
$${ \varepsilon}(\phi, \psi)=:<\phi, T\psi>, \hskip 18pt\hbox{for} \hskip 4pt \phi, \psi\in D(T).$$
Let $|\phi|_\varepsilon^2=|\phi|^2+\varepsilon(\phi,\phi)$. A quadratic form is 
{\it closed} if $\{\phi_n\}$ is a sequence in $D(T)$ with 
$\lim_{n\to\infty}\phi_n=\phi$
 and $lim_{n\to\infty}\varepsilon(\phi_n,\phi_m)=0$, then it follows $\phi\in D(T)$ and $\phi_n$ converges to $\phi$ in $|\, |_\varepsilon$ norm. The {\it closure} of $\varepsilon$ is the least extension of $\varepsilon$ which is closed.  
The closure of ${ \varepsilon}$ is in fact the quadratic form of a unique self adjoint operator $\hat T$, which is a positive extension of $T$ and is called the Friedrichs extension of $T$.

\begin{theorem} {(Von Neumann theorem)}\cite{Yosida} Let $T$ be a closed densely defined operator, then $T^*T$ and $TT^*$ are self-adjoint.
$\label{VonNeu}$
\end{theorem}

Furthermore, If $T$ is symmetric and $T^2$ is densely defined, then $T^*T$ is the Friedrichs extension of $T^2$. See \cite{RE-SI80} [p.181]

Given operators $S$ and $T$, we may define new operators $S+T$, and $ST$. For this, $S+T$ is defined on $D(S)\cap D(T)$, on which $(S+T)\phi=: S\phi+T\phi$,  while  $D(ST)=\{\phi:\phi\in D(T), T\phi\in D(S)\}$, and $(ST)\phi=:S(T\phi)$.

\begin{theorem}[spectral theorem-functional calculus form]
\cite{RE-SI80}

Let $A$ be a self-adjoint operator on $\h$. Then there is a unique map 
$\hat \phi$ from the bounded Borel functions on $R$ into $L(\h)$ so that
 \begin{enumerate}

\item  $\hat \phi$ is an algebraic *-homomorphism. 
\item  $\hat \phi$ is norm continuous, that is, 
	$||\hat \phi(h)||_{L(\h)}\le ||h||_\infty$.
\item Let $\{h_n\}$ be a sequence of bounded Borel functions with 
	$\lim_{n\to \infty}h_n(x)$ converging to $x$ for each $x$ and $|h_n(x)|\le |x|$ for
	all $x$  and $n$. Then, for any $\psi\in D(A)$, 
	$\lim_{n\to \infty}\hat \phi(h_n)\psi =A\psi$.
\item  If $h_n(x)\to h(x)$ pointwise and if the sequence $||h_n||_\infty$
	is bounded, then $\hat \phi(h_n)\to \hat\phi(h)$  strongly.
\item If $A\psi=\lambda\psi$, $\hat \phi(h)\psi =h(\lambda)\psi$.
\item If $h\ge 0$, then $\hat \phi(h)\ge 0$.
\end{enumerate}

\end{theorem}

\section{Semigroups and Generators}

\begin{definition}  A family of operators $\{T_t\}$ on a Banach space 
$(X, | \cdot |)$
is called a one parameter semigroup\index{semigroup: of class $C_0$} (of class $C_0$) if they satisfy the following:
(not  to be confused with the $C_0$ semigroup later in chapter 4.)
\begin{enumerate}
\item $T_{t+s}=T_tT_s, t\ge 0, s\ge 0$.
\item $T_0=I$.
\item For each $t_0\ge 0$, $f\in X$,
       $$\lim_{t\to t_0} |T_tf  -  T_{t_0}f|=0.$$
\end{enumerate}
\end{definition}
\noindent
For such a semigroup, all $T_t$ are bounded operators with $|T_t|\le Me^{\beta t}$, for  $0\le t< \infty$ and some constants $M>0$,  $\beta<\infty$.

The third condition in the definition is equivalent to each of the following if $T_t$ satisfies the semigroup property (1) and (2) in the definition above:\cite{DA80}, \cite{RE-SI80}

\begin{enumerate}
 \item 
 It is weakly continuous at $0$:  $w-\lim_{t \to 0} T_tf=f$, each $f\in X$.  i.e. for any $\phi \in X^*$, the dual space, $\lim_{t\to 0} \phi(T_tf)=\phi(f).$
\item 
The map $(t,f) \mapsto T_tf$ from $[0,\infty)\times X \to X$ is jointly continuous.
\end{enumerate}

\noindent
Denote by $D(\A)=\{f: \lim_{t\to 0} {T_tf -f \over t} \hskip 12pt \hbox{exists}\}$.  Let $\A f=\lim_{t\to 0} {T_tf - f \over t}$, if $f\in D(\A)$. 
The operator $\A$ on $X$ is called the
 {\it infinitesimal generator}\index{generator}, 
which enjoys the following properties:

\begin{theorem}\cite{DA80}
 Let $T_t$ be a semigroup of class $C_0$, $\A$ its generator. Then the following 
hold:

\begin{enumerate}
\item  The operator $\A$ is closed and densely defined.
\item   $T_t\{D(\A)\}\subset D(\A)$.
\item   If $f\in D(\A)$, $T_tf$ is continuously differentiable on $[0, \infty)$ and satisfies:
$$ {{\partial(T_tf)} \over \partial t} =\A(T_tf)=T_t(\A f).$$
\item 
Furthermore if a function $g: [0,a]\to \hbox{Dom}(\A)$ satisifies
 $${{\partial g_t} \over {\partial t}}  =\A g_t$$

\noindent
for all $t\in [0,a]$, then  $g_a=T_ag_0$.
\end{enumerate}
\label{th: properties of semigroups}
\end{theorem}

\begin{definition}
A contraction semigroup\index{semigroup: contraction} is a semigroup of class $C_0$ with $|T_t|\le 1$ for all $t\ge 0$.
\label{de: contraction}
\end{definition}

\noindent
\begin{theorem}[The Reisz-Thorin interpolation\index{interpolation theorem} theorem]\cite{DA-heat}

Let $1\le p_0,$  $ p_1,$   $q_0,$   $q_1\le \infty$. Let $T$ be a linear operator from $L^{p_0}\cap L^{p_1}$ to  $L^{q_0}+L^{q_1}$ which satisfies:
$$|Tf|_{q_i} \le M_i |f|_{p_i}$$
for all $f$ and $i=1,2.$ Let $0<t<1$ and define $p$, $q$ by:
$${1\over p}={t\over p_1}+ {1-t\over p_0}.$$
$${1\over q}={t\over q_1}+{1-t\over q_0}.$$
Then $|T_tf|_q\le (M_1)^t\cdot (M_0)^{1-t}|f|_p$ for all $f\in L^{P_0}\cap L^{p_1}$.
\end{theorem}

In particular if $T$ is a bounded operator both on $L_2$ and $L_\infty$:
then $T$ extends to a operator on $L_p$, for $2<p<\infty$. Furthermore if $T$ is both an $L_2$ and $L_\infty$ contraction, it is an $L_p$ contraction for all such $p$ (restricting to $L^2$ functions).

\noindent
Finally we quote the following theorem on dual semigroups\index{semigroup: dual semigroups} from \cite{DA80}:
\begin{theorem}\cite{DA80}
If $\A$ is a generator of a one parameter semigroup $P_t$ on a reflexive Banach space $X$, then $P_t^*$ is a one parameter semigroup on $X^*$ with generator $\A^*$.
\label{th: dual semigroups}
\end{theorem}

\section{Differential forms\index{differential forms}} \label{differential forms}

Let $M$ be a Riemannian manifold with   Levi-Civita connection $\nabla$ and a
 positive measure $\mu$ given by $e^{2h(x)}dx$, for a smooth function $h$ on $M$, where $dx$ denotes the volume element determined by the metric. Denote by $A^p$ the space of differentiable $p$ forms. Let $\phi$ be a p-form;  there is an element $\phi^\#(x)\in \wedge^p T_xM$ such   that for any $x\in M$, and $v\in \wedge^pT_xM$, $\phi(x)(v)=$

\noindent $<\phi^\#,v>_x$. For two p-forms $\phi$ and $\psi$, we
define:
$$<\phi, \psi>_x=<\phi^\#(x), \psi^\#(x)>_x,$$
$$<\phi,\psi>=\int_M <\phi,\psi>_x \, e^{2h(x)}dx.$$
 
Let $L^2(A^p,\mu)$ be the completion of 
 $\{\phi\in A^p : |\phi|^2=<\phi, \phi><\infty\}$, and 
$L^2(A)=\bigoplus L^2(A^p)$, where $\bigoplus$ denotes direct sums. Then 
$L^2(A)$ is a Hilbert space with the inner product defined above (forms of
 different order  are considered orthogonal).   Those forms in $L^2$ are
 called $L^2$ forms. Let $C_K^\infty(A^p)$ be the space of smooth $p$ forms
 with compact  support, and  $C_K^\infty$  the\index{$C_K^\infty$} space of
 smooth functions with compact    support. There is the usual exterior
 differential operator $d$ on $C_K^\infty$ with adjoint $\delta^h$,
here $\delta^h$\index{$\delta^h$} is the formal adjoint of $d$  for $\mu$:

$$d: C_K^\infty(A^p) \to C_K^\infty(A^{p+1}),$$
$$\delta^h: C_K^\infty(A^p) \to C_K^\infty(A^{p-1}). $$ 

Denote by $d^*$ the adjoint of $d$ in  $L^2(A^p, \mu)$ and $\delta^*$  the
 $L^2$ adjoint of operator $\delta$,  so that $d^*=\delta^h$, $(\delta^h)^*=d$
 when restricted to $C_K^\infty$.  Notice that  $C_K^\infty$ is dense in 
 $L^2$,  so all  the operators    concerned have dense domain.

Denote by   $(\wedge^p T_xM)^*$  the dual space of the antisymmetric tensor
 tangent bundle of order $p$ on $R^1$.  Let
 $(\wedge T_xM)^*=\bigoplus (\wedge^p T_xM)^*$.  For each $e\in T_xM$, there
 is the annihilation operator
 $i_e: (\wedge^{p+1}T_xM)^*\to (\wedge^p T_xM)^*$  given by:

$$i_e(\phi)(v_1, \dots, v_p) =\phi(e_1, v_1, \dots, v_p).$$
Here we have identified $(\wedge^p T_xM)^*$ with antisymmetric multilinear
 forms.

Let  $Y$ be a vector field on $M$ with $S_t$ the corresponding solution flow. 
 There is the interior product of $\phi$ by $Y$   given by:
$$i_Y\phi(v_1, \dots, v_{p-1})=\phi(Y(x), v_1, \dots, v_{p-1})$$
for $v_i\in T_x M$.  If $\psi $ is a 1-form, we write $i_{\psi}\phi$ for
 $i_{\psi^\#}\phi$.  Here $\#$ denotes the adjoint. 
There is also the Lie derivative $L_Y$ of $\phi$ in direction $Y$:
$$L_Y\phi(v_1,\dots, v_p)={d \over dt}\phi(TS_t(v_1), \dots, TS_t(v_p))|_{t=0}.$$
Here $v_i\in T_xM$. They are related by the following formula \cite{AB-MA}:
$$L_Y\phi =d(i_Y\phi) +i_Y(d\phi).$$

Take an orthonormal basis $e_1, \dots, e_n$ of $T_xM$, the formal adjoint 
 $\delta$ of the exterior differential operator on $L^2(M, dx)$ is given by:

$$ (\delta\phi)_x = -\sum_1^n i_{e_k} (\nabla _{e_k}\phi). $$
For a function  $f$ on  $M$, 

	$$\delta(f\phi) = f \delta \phi - i_{\nabla f} \phi.$$
Also for $\phi\in D(\delta)$, $\psi \in C_K^\infty$,

$$\int <d\phi, \psi> e^{2h}\, dx = \int <\phi, \delta(e^{2h} \psi)>\, dx
=\int <\phi, \delta\psi -2i_{\nabla h} \psi> e^{2h} \, dx. $$

\noindent
  Thus $\delta^h=\delta - 2 i_{\nabla h}$. Similar arguments show that:
  $\delta^h=e^{-2h}\delta e^{2h}$ and
\begin{equation}
\delta^h(f\phi)=f\delta^h\phi - i_{\nabla f}\phi.
\end{equation}
Let $\triangle$ be the {\it Hodge-Laplace operator} on forms defined as:
$$\triangle=-(d+\delta)^2,$$
 and $\triangle^h$ the {\it Bismut-Witten Laplacian}\index{Bismut-Witten
 Laplacian}  \label{de: Witten Laplacian}\index{$\triangle^h$}
  defined by:
$$\triangle^h=-(d+\delta^h)^2.$$
Clearly $\triangle^h=\triangle +2L_{\nabla h}$. The restricted operators on
 $q$ forms are denote by $\triangle^q$ and $\triangle^{h,q}$ respectively.
 If $q=0$, we simply write $\triangle$ and $\triangle^h$.

\noindent
{\bf  Divergence theorem }  Let $Y$ be a vector field. Define the divergence
 of $Y$ to be: $\hbox {div} Y=\hbox{trace}\nabla Y$.  It is the formal adjoint
 of $-\nabla$ on $L^2(M, dx)$. The h-divergence of $Y$ is given by:

\begin{equation}
\hbox{div}_h Y=e^{-2h} \hbox{div} (e^{2h} Y).
\label{eq: divergence 1}
\end{equation}

This is the formal adjoint of $-\nabla$ in $L^2(M, e^{2h} dx)$:
$$\int f (\hbox{div}_h Y) e^{2h}\, dx = -\int <Y, \nabla f> e^{2h} \, dx. $$
The divergence theorem holds for h-divergences, i.e. let $U$ be a precompact
 open set in $M$ with piecewise smooth boundary, then:

\begin{equation}
\int_U \hbox{div}_h Ye^{2h}\, dx = \int_{\partial U} <Y, \nu> e^{2h}\, dS.  
\label{eq: divergence theorem}
\end{equation}
Where $\nu$ is the unit outer normal vector to $\partial U$, and $dS$ the
 Riemannian surface area element(corresponding to $dx$). The equation can be
 proved from  equation ~\ref{eq: divergence 1} and  from the usual divergence
 theorem with  $h=0$.

\section{Parabolic regularity  etc.}

Recall  that $h$ is a smooth function on $M$. Let $P_t^h\phi$ be a $L^2$
 solution to  the following equation on $q$ forms:
\[ \left\{\begin{array}{cl}
{\partial \over \partial t}u(x,t)&=\left(\half \triangle^h\right)u(x,t)\\
u(x,0)&=\phi
\end{array}\right. \]

Then $P_t^h\phi$ is in fact a classical solution(i.e. $C^2$ in $x$
and $C^1$ in $t$). Furthermore it is smooth  if $\phi$ and $h$ are.
 See  \cite{NARA} and \cite{Eells}.

The following theorem is given by Strichartz for the Laplacian, but it is
 easy to see it is true for ${1 \over 2}\triangle^h$ (defined on
 page ~\pageref{de: Witten Laplacian}, c.f. theorem
 $~\ref{th: chernorffselfadjoint}$).  Let $\heatsemi$ be 
the heat semigroup defined by functional analysis.

\begin{theorem}  \cite{STRIana}
There is a heat kernel\index{heat kernel}\index{$p_t^h(x,y)$} $p_t^h(x,y)$ 
satisfying:
\begin{enumerate}
\item
The function $p_t^h(x,y)$ is smooth on $R^+\times M \times M$, symmetric in
 $x$ and $y$, and is strictly positive for $t>0$.
\item
$\int_M p_t^h(x,y)e^{2h(x)} dx \le 1$ for all $x$ and $t>0$, and
\begin{equation}   \heatsemi f(x)=\int p_t^h(x,y)f(y)e^{2h(y)} \, dy
\label{eq: integral kernel}
\end{equation}
for all $f\in L_2$.
\item
$|\heatsemi f|_p \le |f|_p$, for all $t>0$, and $f\in L^2\cap L^p$, 
$1\le p\le \infty$ with $\lim_{t\to 0}|\heatsemi f-f|_p=0$,  for
 $p\not = \infty$, and
\item
${\partial \over \partial t} \heatsemi f
={1\over 2}\triangle^h \heatsemi f$ for all $f\in L^2$ and $t>0$ , and these
 properties continue to hold for all $f\in L^p$,  $1\le p\le \infty$ if we
 define  $\heatsemi f$ by equation~\ref{eq: integral kernel}. Moreover we
 have uniqueness of the semigroup for $1<p<\infty$ in the following sense:
if $Q_t$ is any strongly continuous contractive semigroup on $L^p$ for fixed
 $p$ such that $Q_tf$ is a solution to
 ${\partial \over \partial t}=\half \triangle^h$,  then $Q_tf=\heatsemi f$.
\end{enumerate}
\label{th: STRana}
\end{theorem}

\chapter {On the Bismut-Witten  Laplacian and its semigroups}

{\bf Main result:} From the essential self-adjointness of $d+\delta^h$, we get 
the Hodge decomposition theorem and a result on $dP_t^h=P_t^hd$, which we use 
extensively later. These results are generally considered well known, but there
does not appear to be a suitable reference which gives these in exactly the
 context we want, especially for the Bismut-Witten Laplacian $\triangle^h$. 
This problem arose when I was trying to work out $dP_t^h=P_t^hd$, and occurred
again when D. Elworthy and S. Rosenberg were working on their paper
 \cite{EL-RO92}. 
This chapter is the result of discussions among the  three of us, part of the
 result is in fact used  in \cite{EL-RO92}.

\section{The self-adjointness of weighted Laplacians}

Chernoff \cite{CH73} has proved that $d+\delta^h$ and all its powers are
 essentially  self-adjoint on $C_K^\infty$ (with $h=0$, but exactly the same
 proof gives the result for $h\not = 0$ making use of the h-divergence theorem
 stated on page ~\pageref{eq: divergence theorem} as seen later in the
 context).    In particular the Hodge-Laplacian operator
 $\triangle^h=-(d+\delta^h)^2$ is essentially self-adjoint.
 From this we obtain $\overline{d+\delta^h}=\bar d +\overline{\delta^h}$  and
 the Hodge-Kodaira decomposition, and therefore : $d^*=\bar \delta^h$ and 
$(\delta^h)^*=\bar d$.  Consequently we get: $\overline{\triangle^h} =-(\bar d
\overline{\delta^h} +\overline{\delta^h} \bar d)$.  Thus all the operations on
 $\triangle^h$ are exactly as for  $\triangle$.    We will give Chernoff's
 proof for the self-adjointness, adapted to our context.
 
\bigskip

First we introduce the terminology. Let $L$ be a first order differential
 operator 
on the complexified cotangent bundle $(\wedge T_xM)^*_\C$. Consider 
${\partial \over \partial t}=L$, which is called a symmetric hyperbolic system
 if $L^* +L$ is a zero order operator, i.e. it is given
locally by multiplication by a matrix valued function of $x$.  The symbol
 $\sigma$ for $L$ is a linear map  from 
$(\wedge T_xM)^*_\C \to (\wedge T_xM)^*_\C$ for each $x\in M$ and 
$v\in T_x^*M$:
\begin{equation}
\sigma(v)(e) = L(f\phi)(x) -f L\phi(x).
\end{equation}
Here $f\in C^\infty (M)$ with $df(x)=v$, and $\phi: M\to (\wedge TM)^*$  is a
 smooth form on $M$ with $\phi(x)=e$. The local velocity of propagation $c$
 associated to $L$ is defined as a function on  $M$:
$$c(x)=\sup_v \{||\sigma(v)(x)|| : v\in T_x^*M, ||v||=1\}.$$
We are mainly interested in the following operator $L=i(d+\delta^h)$ for the 
proof of the next two theorems.
It indeed gives rise to a hyperbolic system since $L+L^*=0$, and has  symbol:
$\sigma(v)(e)= v\wedge e -i_{v^\#} e$. Here $v^\#$ is defined by: 
$v(w)=<v^\#,w>$ for $w\in T_xM$. The local velocity for $L$ is constantly $1$, 
from the following:
 
For each fixed $v\not = 0$, we may write $e$ into the sum of  two orthogonal
 parts: $e=e_0+e_1$ satisfying: $v\wedge e_0=0$, and $|v\wedge e_1|=|v||e_1|$,
 by letting  $e_0$ be the component of $e$ along $v$. So $i_{v^\#}e_1=0$. Thus

$$ |\sigma(v)| = \sup_{|e|=1} (|v\wedge e  -  i_{v^\#} e|)
	=\sup_{|e|=1} (|v\wedge e_1| +|i_{v^\#} e_0|)
	=|v|.$$
 
\begin{lemma}[Energy inequality] \cite{CH73} Let $M$ be a Riemannian manifold
 with Riemannian distance $\rho$. The geodesic ball  radius $R$ centered at 
$x_0$  is denoted by  $B(x_0,R)$.  Let $x_0\in M$, $u$ a smooth solution to 
 ${\partial \over \partial t}=L$ on $[0,a]\times B(x_0, R)$. Let $r$ and $a$ be
 real numbers with   $r+a\le R$, then the  following inequality holds:

$$\int_{B(x_0,r)} |u(a)|^2 e^{2h}\, dx
 \le \int_{B(x_0, r+a)} |u(0)|^2 e^{2h}\, dx. $$

In particular if $u(0)$ vanishes on $B(x_0, r+a)$, then $u$ vanishes throughout
 $K$ for $K=\{(t,x): 0\le t\le a , \rho(x,x_0)\le r+(a-t)\}$.
\end{lemma}

\noindent
Proof: 
 We first introduce a (h-)divergence free vector field $Z$ on $[0,a]\times M$ 
as follows.  Let $f$ be a real valued function on $[0,a]\times M$,  denote by
${\partial f \over \partial t}$  the differential of $f$ in the time component,
 and  $df$ the space differential. Then define $Z$ by:

$$df(Z)=<Z, \nabla f>  
=|u_t|^2{\partial f \over\partial t} -<u_t, \sigma(df)u_t>.$$

Write $Y(f)=<u_t, \sigma(df)u_t>$. Notice $Z$ is given in the form of the sum 
of time and space parts, so therefore is the divergence:

$$\hbox{div}_h Z
= 2<u_t, {\partial \over \partial t} u_t> - \tilde{\hbox{div}}_h Y,$$

\noindent
where $\tilde{\hbox{div}}_h Y$ denotes the h-divergence of $Y$ in its space
 variables. 

But $\tilde{\hbox{div}}_h Y=< (L-L^*)u_t, u_t>$, since for any function $f\in C_K^\infty$:

\begin{eqnarray*}
\int_M f\tilde {\hbox {div}}_h Y\, e^{2h}  dx
 &=& -\int_M <Y, \nabla f> \,e^{2h} dx \\
	&=&-\int <u_t, L(f u_t) -fLu_t> \,e^{2h} dx\\
	&=&-\int_M \{<L^* u_t, u_t>- <u_t, Lu_t>\}f\,e^{2h} dx\\
	&=&-\int_M f<u_t, (L^* -L)u_t>\, e^{2h} dx. 
\end{eqnarray*}
%\nopagebreak

\noindent
Thus $\hbox{div}_h Z=<u_t, (L+L^*)u_t> =0$. Applying divergence theorem for the
 h-divergence, we get:
\begin{eqnarray*}
0&=&\int_K \hbox{div}_h Z \,e^{2h} dt dx \\
&=&\int_{\partial K} <Z, \nu> e^{2h} dS\\
	&=&\int_{B(x_0,r)} |u(a)|^2\, e^{2h}dx 
-\int_{B(x_0,r+a)} |u(0)|^2\, e^{2h}dx
	+ \int_\Sigma <Z, \nu> \,e^{2h} dS.
\end{eqnarray*}

%
%
%\begin{figure}
%\vspace{50mm}
%\end{figure}
%

 \noindent
Here $\Sigma$ is the mantle part of the boundary of $K$, and $\nu$ the unit
 outer normal vector to $\Sigma$, which is in fact given by:  
$\nu =\alpha\nabla \phi = \alpha (1, \nabla \rho)$, 
for  $\phi(t,x)=t +\rho(x_0,x)$, and some positive constant $\alpha$. 
\begin{eqnarray*}
<Z, \nu>&= &\alpha|u_t|^2 -<u_t, \sigma(\alpha\, d\rho)u_t>\\ 
	&\ge& \alpha |u_t|^2- \alpha |d\rho||u_t|^2 =0.
\end{eqnarray*}
 The required inequality follows.  \hfill \rule{3mm}{3mm}

Let $L=i(d+\delta^h)$.  We shall view the equation 
${\partial u \over \partial t}=i(d+\delta^h)u$ as a symmetric hyperbolic
system of partial differential equations as follows:

\noindent
Let $\{\omega_1,\dots, \omega_{2n}\}$ be local sections of $A=\bigoplus_{p=0}^n
A_p$ which are pointwise orthonormal. Then we may wite
$v=\sum_1^n v_\alpha(x)\omega_\alpha$ and
$$(d+\delta^h)v=\sum_{i,\alpha,\beta}A_{i\alpha\beta}(x)
{\partial v_\alpha \over \partial x_i}\omega_\beta+
\sum_{\alpha, \beta}B_{\alpha \beta}v_\alpha(x)\omega_\beta.$$
Let $v, w\in C^\infty(A)$, the space of smooth forms in $A$. We shall assume
that one of them has compact support. The identity
$$<(d+\delta^h)v,w>=<v,(d+\delta^h)w>$$
and an integration by parts show that
$$A_{i\alpha\beta}(x)=-A_{i\beta\alpha}(x), \hskip6pt
\hbox{any} \hskip3pt i,\alpha,\beta, \hskip 3pt \hbox{and}\hskip 3pt x.$$
We shall let $A_i(x)$ denote the matrix whose entries are
 $A_{i\alpha \beta}(x)$. On setting $u=v+iw$, the equation:
${\partial u\over \partial t}=i(d+\delta^h)u$ can be written as
$${\partial \over \partial t}\left(\begin{array}{c}\vec v\\  \vec w
\end{array}\right)
=\sum_1^n\bar A_i(x){\partial \over \partial x^i} \left(\begin{array}{c}
\vec v\\ \vec w\end{array}\right)
 +\bar B(x)\left(\begin{array}{c}\vec v\\ \vec w\end{array}\right),$$
where 
$\vec v=(v_1, \dots,  v_{2n})$ and 
$\vec w=(w_1,\dots,  w_{2n})$ and
$\bar A_i(x)=\left(\begin{array}{cc}0&A_i(x)\\-A_i(x)&0\end{array}\right)$ is  
 symmetric.

 Thus locally there is a unique smooth solution to the Cauchy problem for 
this  hyperbolic equation by standard theory. See e.g. \cite{John}.

\begin{theorem}\cite{CH73}
Let $M$ be a complete Riemannian manifold, and let $L={\sqrt  {-1}}(d+\delta^h)$. 
There is a unique  smooth solution to the  hyperbolic equation: 
${\partial \over \partial t} =L$, for each initial value $u_0\in C_K^\infty$. 
Moreover $u_t\in C_K^\infty$, each $t$. 
\end{theorem}

\noindent
Proof: Take a local trivialization of the tangent bundle $TM$.  For  $y\in M$, 
there is a constant $r(y)>0$ such that $B(y,r)$ is contained in a single chart.
Then by the local existence and uniqueness theorem: 
for each $u(0)$ smooth on $B(y,r)$, there is a unique smooth solution on 
$K_y=\{(t,x): 0\le t \le r, \rho(x,y)\le r-t\}$, from the energy inequality.
 Notice the propagation speed is identically $1$ for $L$.

 For $R>0$, $x_0\in M$, the geodesic ball $B(x_0, R)$ is compact on a complete
Riemannian manifold. Let $u(0)$ be smooth on $B(x_0, R)$. 
 We can find a common $r>0$ for all points in the ball $B(x_0, R-r)$ such that the
 statement above holds, given initial data on $B(x_0, R)$. Moreover the
 solutions coincide on overlapping areas from the local uniqueness. Altogether they
 define a solution on the truncated cone:

$$\{(t,x): 0\le t\le {r}, \rho(x, x_0)\le R-t\}.$$
The solution at $t={r}$ in turn serves as initial data on 
$B(x_0, R-{r})$, then we have a solution on:
$\{{r} \le t \le 2r, \rho(x, x_0)\le R- {r} -t\}$.
Continuing  with the procedure, we obtain a unique solution 
$\{(t,x): 0\le t\le 1, \rho(x, x_0)\le R-t\}$.   
 
In particular this show that the solution up to time $N$ on $B(x_0, R-N)$ is
determined by solution at time $N-1$ on the ball $B(x_0, R-N+1)$, and 
therefore by initial value on $B(x_0, R)$.

Let $u(0)\equiv 0$ on $M$. For each number $N$ and $R$,    $u\equiv 0$ on
 $[0,N]\times B(x_0, R)$ from the above argument. So the uniqueness holds.

We show next that there is a globally defined solution for each smooth data
 $u_0$ with compact support inside the ball $B(x_0, a)$, some $a$. As has been 
shown there is a unique solution on $K=[0,1]\times B(x_0, a+3-t)$. Moreover the
solution vanishes outside $B(x_0, a+1)$.  Extend the solution to $[0,1]\times M
$ by setting  it to be $0$ outside $K$. The solution so defined is smooth and h
as support in $B(x_0, a+1)$  at time $t=1$, which serves in turn as initial 
data and gives us a solution on  $[0,2]\times M$. We can by this means 
propagate the solution for all time.
It is clear that the solution $u_t$ thus obtained has compact support for each
 $t$.  \hfill\rule{3mm}{3mm}

\begin{theorem} \cite{CH73}
Denote by $V_tu$ the solution given in the theorem above for $u\in C_K^\infty$.
Then $V_t$ is a unitary group from $C_K^\infty \to C_K^\infty$ with
$L(V_tu)=V_t(Lu)$, and  extends to $L^2$ as a unitary group.
 \end{theorem}

 \noindent
Proof:
We only need to prove the unitary part. Take $u$, $w$ from $C_K^\infty$.Than
\begin{eqnarray*}
{d \over dt} <V_tu, V_tw> &=& <LV_tu, V_tw> +<V_tu, LV_tw> 	\\
	&=&<(L+L^*) V_tu, V_tw>=0.
\end{eqnarray*}
The fact that $V_t$ is a group comes from  the uniqueness of the solution,
 and $L(V_tu)=V_t(Lu)$  follows from the standard semigroup result. 
\hfill \rule{3mm}{3mm}

\begin{theorem}\cite{CH73} $\label{th: chernorffselfadjoint}$
The operator $T=d+\delta^h$ and all its powers are essentially self-adjoint.
\end{theorem}
\noindent
Proof:
Let $A=T^n$, $n>0$. It is a symmetric operator. According to theorem 1.4 we only need to show $\psi=0$ if $A^*\psi=\pm i\psi$, $i=\sqrt {-1}$. Suppose $A^* \psi = i\psi$. Let 
$\phi \in C_K^\infty$. Consider $<V_t\phi, \psi>$ which is bounded in $t$ since
$V_t$ is a unitary operator.

\begin{eqnarray*}
{d^n \over dt^n}<V_t \phi, \psi>&=&i^n <T^n V_t\phi, \psi>\\
	&=& i^n<V_t\phi, A^*\psi> =i^n<V_t\phi, i\psi>\\
	&=& -i^{n+1}<V_t\phi, \psi>.
\end{eqnarray*}

\noindent
So $<V_t\phi, \psi>=Ce^{\alpha t}$, for some constant $C$. Where $\alpha$ satisfies: $\alpha^n +i^{n+1}=0$, and so is not a pure imaginary. Thus 
 $<V_t\phi, \psi>=0$, since it is bounded in $t$.
In particular $<\phi, \psi>=0$. Thus $\psi=0$. A similar argument works with 
$A^*\psi =-i\psi$.  \hfill \rule{3mm}{3mm}

\section{The Hodge decomposition theorem}

  To explore the relationship of these operators, let us first notice that
$\bar d\subset \delta^*$, and $\bar \delta \subset d^*$ before giving the following lemmas:

\begin{lemma} On $L^2(M, e^{2h}dx)$, 
 $D(\bar d^2)=D(\bar d)$, $D((\bar \delta^h)^2)=D(\bar \delta^h)$,
and $(\bar d)^2=(\bar \delta^h)^2=0$. It is also true that: $ D(\bar d)=
 D((\delta^h)^*\bar d)$, and $(\delta^h)^*\bar d=0$. Similarly for $d^*\overline {\delta^h}$.
\end{lemma}

Proof: Take $\phi \in D(\bar d)$, there are  therefore $\phi_n$ in $C_K^\infty$ 
converging  to $\phi$ in $L^2$, and $d\phi_n\to \bar d\phi$ in $L^2$. But $d^2\phi_n=0$, so $\bar d\phi\in D(\bar d)$, and $\bar d^2\phi=0$. 

Let $\phi\in D(\bar d)$, then for all $\psi \in D(\delta^h)=C_K^\infty$, we have:
 $$<\bar d \phi, \delta^h\psi>= <\phi, d^*\delta^h \psi> =0.$$
So $\bar d \phi\in D(\delta^*)$ and $(\delta^h)^*(\bar d\phi)=0$.
The rest  can be proved analogously. \hfill \rule{3mm}{3mm}

\begin{proposition} Let $M$ be a complete Riemannian manifold.  Then:

 $\bar d+{\overline \delta^h}={\overline{d+\delta^h}}=d^*+(\delta^h)^*$.
\end{proposition}
 
\noindent
Proof:  We write $\delta$ for $\delta^h$ in the proof. Since both $d^*$, and $\delta^*$ are closed from theorem 1.2, so is $d^*+\delta^*$ since $d^*$ and $\delta^*$ map to different spaces. For the same reason 
$\bar d+\bar \delta$ is also a closed operator. Noticing $\bar d\subset \delta^*$, $\bar \delta\subset d^*$, we have:

\begin{equation}
{\overline{d+\delta}}\subset \bar d+ \bar \delta \subset d^*+\delta^*.
\label{eq:operator1}
\end{equation}
We also have:
\begin{equation}
d^*+\delta^*\subset (\bar d +\bar\delta)^*
\label{eq:operator2}
\end{equation}
from the following argument:

\noindent
Let $\phi\in D(d^*+\delta^*)$.  For any $\psi\in D(\bar d+\bar \delta)$, there is the following:

$$<(\bar d+\bar\delta)(\psi),\phi>
=<\psi, (\bar d)^*\phi>+<\psi,(\bar\delta)^*\phi>=<\psi, (d^*+\delta^*)\phi>.$$

The last equality comes from theorem ~\ref{th: properties of adjoint}.  We now
 take adjoint of both sides of $~\ref{eq:operator1}$  and $~\ref{eq:operator2}$
 and get:

\begin{equation}
(d^*+\delta^*)^*\subset (\bar d+\bar \delta)^*\subset (\overline{d+\delta})^*,
\label{eq:operator3}
\end{equation}
\begin{equation}
\bar d+\bar \delta=(\bar d +\bar \delta)^{**}\subset (d^*+\delta^*)^*.
\label{eq:operator4}
\end{equation}

\noindent
But by theorem $~\ref{th: chernorffselfadjoint}$,  
$(\overline{d+\delta})^*={\overline{d+\delta}}.$

So combine  $~\ref{eq:operator1}$ with $~\ref{eq:operator3}$ 
and $~\ref{eq:operator4}$, we have:
$${\overline{d+\delta}} \subset \bar d+\bar \delta    \subset (d^*+\delta^*)^*
\subset \overline {d+\delta}.$$
Therefore $\bar d+\bar \delta=(d^*+\delta^*)^*={\overline{d+\delta}}$, and each
is a self-adjoint operator. Thus we have proved that:
 $$\bar d+\bar \delta=d^*+\delta^*={\overline{d+\delta}}. $$ 
\hfill \rule{3mm}{3mm}
 
\noindent
{\bf Remark:} Let $\overline{\triangle^h}=-\overline{(d+\delta^h)^2}$
\label{Witten Laplacian}, it is now clear that:
  $$\overline {\triangle^h}=-(\bar d +\bar \delta^h)^2=-(\bar d \bar{\delta^h}+\bar{\delta^h} \bar d).$$
\noindent
Proof: First we know that $\bar d+\bar \delta^h$ is self adjoint from
 $~\ref{th: chernorffselfadjoint}$.  By theorem $~\ref{VonNeu}$, 
$\bar\triangle^h$ is self adjoint. But there is only one self adjoint extension
 for $\triangle^h$. The result follows.
 Furthermore $(\bar d +\bar \delta^h)^2$ is  in fact the Friedrichs extension 
of $\triangle^h$ from the remark after theorem $~\ref{VonNeu}$. 

The following is a standard result:

\begin{proposition}
Let $T$ be a self-adjoint operator on Hilbert space $\h$ with dense domain. Then
$$H=\overline {Im(T)}\bigoplus ker(T).$$
\end{proposition}

\noindent{\bf Proof:}
Let $\phi\in D(T)$, $\psi\in ker(T)$. Then
$$<T\phi, \psi>=  <\phi, T\psi> =0.$$
So $\overline{Im(T)}$ is orthogonal to $ker(T)$. 

\noindent
Let $\phi\in (\overline{Im(T)})^\bot$. Then by definition for any $\psi\in D(T)$,
 we have:
$<\phi, T\psi>=0$. However this shows that: $\phi\in D(T^*)=D(T)$, and $T^*\phi=0$.
Thus $T\phi=0$. Therefore $\overline{Im(T)}^\bot \subset Ker(T)$. This finishes the
proof. \hfill \rule{3mm}{3mm}

We are ready to prove the following Hodge decomposition theorem:
\begin{theorem} Let $M$ be a complete Riemannian manifold with measure 
 $e^{2h}dx$,  here $h\in C^\infty(M)$.  Let
 $L^2(\h)=\{\phi: \bar d \phi=\bar \delta^h \phi=0\}$. Then:
$$L^2\Omega^p
=\overline {Im(\delta^h)}\bigoplus \overline {Im(d)}\bigoplus L^2(\h)$$
and
$$L^2(\h)=Ker\overline {\triangle^h}.$$
$\label{Hodge decomposition}$
\end{theorem}

\noindent 
{\bf Proof:} 
From the propositions above, we know:
$$L^2=\overline{Im(\bar d +\bar \delta)} \bigoplus Ker(\bar d +\bar \delta ).$$
We only need to prove the following:
$$\overline{Im(\bar d +\bar \delta)}
=\overline{Im(\bar d)}\bigoplus \overline{Im(\bar \delta)}.$$

It is clear the two spaces on the right hand side are orthogonal to each other.
Take $\phi \in \overline{Im(\bar d +\bar \delta)}$. By definition there are
 $\phi_n\in D(\bar d +\bar \delta)$, such that 
$(\bar d +\bar \delta)\phi_n\to \phi$. Since $\bar d\phi$ and $\bar \delta \phi$
are in different spaces, we may write $\phi=\phi_1+\phi_2$ such that:
$\bar d\phi_n\to \phi_1$, and $\bar \delta \phi_n \to \phi_2$. 
 Thus follows: 
$$\overline{Im(\bar d +\bar \delta)}
\subset \overline{Im(\bar d)}\bigoplus \overline{Im(\bar \delta)}.$$

Next let $\phi\in D(\bar d)$. Let
 $\psi\in Ker(\bar d +\bar \delta)=Ker(\bar d)\cap Ker(\bar \delta)$.
 Take $\psi_n\in C_K^\infty$ converging to $\psi$ such that
 $\delta\psi_n\to \bar \delta \psi =0$. Then there is the following:
  $$<\bar d \phi, \psi>=\lim_{n\to \infty}<\bar d \phi, \psi_n>
=\lim_n<\phi, \delta \psi_n>=<\phi, \bar \delta\psi>=0.$$
Therefore we get:  
$\overline{Im(\bar d)}\subset ker(\bar d +\bar \delta)^\bot $. Similarly
$\overline{Im(\bar \delta)}\subset ker(\bar d +\bar \delta)^\bot$.
Thus 
$$L^2\Omega^p
=\overline {Im(\delta^h)}\bigoplus \overline {Im(d)}\bigoplus L^2(\h).$$

\noindent
from $\overline{Im(\bar \delta^h)}=\overline{Im( \delta^h)}$, and
$\overline{Im(\bar d)}= \overline{Im(d)}$.

For the proof of the second part, recall  $\delta^*\bar d\phi=0$. 
For  $\phi\in D(\bar \triangle)$,
 $$-<\bar \triangle\phi, \phi> = <(\bar d+\bar \delta)\phi, (\bar d+\bar \delta)\phi>=|\bar d\phi|^2 +|\bar\delta \phi|^2.$$
So  $L^2(\h)=Ker(\bar \triangle).$
The proof is complete. \hfill \rule{3mm}{3mm}

\begin{proposition}   Let $M$ be a complete Riemannian manifold, then:

  $\bar d=(\delta^h)^*$,  and $\bar \delta^h=d^*$.
\end{proposition}

\noindent
{\bf Proof}  First notice $\delta^{h*}$ is  an extension of $\bar d$.  Let
 $\phi\in D(\delta^{h*})$, we may write $\phi=\phi_1+\phi_2+h$.  Here
$h\in L^2(\h)$, $\phi_1\in \overline {Im (d)}$,  and
$\phi_2\in \overline {Im(\delta^h)}$.
However $\overline{Im(d)} \subset D(\bar d)$, $\overline{Im(\delta^h)}\subset D(\bar \delta^h)$. So $h+\phi_1\in D(\bar d) \subset D(\delta^{h^*})$, and $\phi_2$ is in the domain of  $(\delta^h)^*$.  But then $\phi_2\in D(\delta^{h*})\cap D(\bar \delta^h)\subset D(\delta^{h*})\cap D(d^*)=D(\bar d+ \bar \delta^h)$, by lemma 1.
So $\phi_2\in D(\bar d)$, and therefore $\phi\in D(\bar d)$. 

\noindent
This gives: $D(\delta^{h*})\subset D(\bar d)$, i.e. $\delta^*=\bar d$. Similarly we can prove $\bar\delta^h=d^*$. \hfill \rule{3mm}{3mm}

\noindent
{\bf Remark:} Restricted to $C^1$ forms, the above lemma can be obtained by an
approximation method as suggested  by Gaffney \cite{GA51}.

With these established, we are happy to use $d$, $\delta^h$, and  $\triangle^h$ for their closure without causing any confusion.

\section{The semigroup associated with $\triangle^h$}

Since $\triangle^h$ is non-positive, there is a semigroup 
 $e^{{1\over 2} t \triangle^h}$\index{$\heatsemi$} with generator ${1\over 2} \triangle^h$  defined by the spectral theorem. Furthermore $\heatsemi \alpha$ solves: 

\begin{equation}\left\{
\begin{array}{clc}
{\partial \over \partial t} g_t &=&{1\over 2} \triangle^h g_t\\
g_0&=&\alpha.
\end{array}\right.
\label{eq: heat with h}
\end{equation}
\noindent Here $\alpha\in L^2(\Omega)$ is a $L^2$ form.

\begin{proposition}\label{pr: exchangeability}
Let $M$ be a complete Riemannian manifold. Let  $\alpha\in D(d)\cap D(\delta^h)$, then $d\heatsemi\alpha =\heatsemi(d\alpha)$, and
 $\delta^h \heatsemi\alpha  =\heatsemi(\delta^h \alpha)$.
\end{proposition}

\noindent
{\bf Proof:}  See Gaffney\cite{GA54} for a proof  of the lemma for $C^1$ 
forms (h=0). 
Let $h(y)=ye^{-{1\over 2} ty^2}$, $h_n(y)=\chi _{[-n,n]}(y)(ye^{- {1\over 2} ty^2})$, and  $g_n(y)=\chi _{[-n,n]}(y)y.$  
Denote by $h_n(\A)$ and $g_n(\A)$ the operators defined by the spectral theorem corresponding to  an operator $\A$. Then 
$$h_n(\A)=g_n(\A)e^{-{1\over 2} t\A^2}=e^{- {1\over 2} t\A^2}g_n(\A).$$
Now $|g_n(y)|\le |y|$, all $n$ and $g_n(y)\to y$, so the spectral theorem gives the following convergence result:
$$\lim_{n\to \infty} g_n(\A)\psi =\A \psi,  \hskip 6pt \psi\in D(\A).$$
Thus  for $\psi \in D(\A)$,  we have (  $P_t\psi$ is automatically in  $D(\A)$):

$$\A e^{- {1\over 2} t\A^2}\psi  
 =\lim_{n \to \infty} g_n(\A)e^{- {1\over 2} t\A^2}\psi
 =\lim_{n \to \infty} e^{- {1\over 2} t\A^2}g_n(\A)
 =e^{- {1\over 2} t\A^2}\A.$$

Now let $\A=d+\delta^h$. Let $\psi \in D(d+\delta^h)$, then:
$$(d+\delta^h)\heatsemi \psi
=\heatsemi(d+\delta^h) \psi.$$
Thus the proof is finished noticing both sides of the equality above consist of
orthogonal forms.  \hfill \rule{3mm}{3mm}

\chapter{Invariant measures and ergodic properties of BM on manifolds of finite
 h-volume}

Let $M$ be a complete Riemannian manifold. There is an  invariant measure 
for h-Brownian motion, i.e. for the diffusion
 process with generator $\half \triangle^h$.  In this chapter we give a direct
 proof of the existence of invariant measures  and deduce some ergodic properties, 
which are used in chapter 6 and 7.  These ergodic properties are essentially 
 known and well treated in \cite{Kunitabook} and \cite{skorohod},  at least when the
 weight $h$ is zero. However in their treatment, the processes are required to have
 the $C_0$ property if the manifold is not compact.  Our treatment avoids the 
problem of having to assume the $C_0$ property for the Brownian motion concerned. 
\bigskip

The following lemma is a standard result from semigroup theory:

\begin{lemma}
Let $M$ be a complete Riemannian manifold given the measure $e^{2h(x)}dx$.
 Let $f\in L^2$. Denote by $P_t^h$ the heat semigroup. Then 

\noindent
$\int_0^t P_s^h fds \in D(\bar \triangle^h)$ and moreover:
$$\half\bar \triangle^h(\int_0^t P_s^hf ds)=P_t^hf-f.$$
\end{lemma}
To prove this, we only need to show for all $g\in C_K^\infty$:
$$\int_M \left(\int_0^t P_s^hf ds\right) \left(\half\bar \triangle^h g\right)
 e^{2h}  dx
=\int_M\left(P_t^hf-f\right)g e^{2h} dx.$$
using the fact  that 
$\half(\overline \triangle^h)^*=\half\overline \triangle^h$. See chapter 2.
Noticing  the dual 
semigroup $(P_s^h)^*$ equals $P_s^h$, we obtain:

\[\begin{array}{ll}
&\half\int_M \left(\int_0^t P_s^hf ds\right)\left(\bar \triangle^h g\right)\,
e^{2h}dx 	 \\
=&\half\int_0^t\int_M P_s^hf(\triangle^h g)\ e^{2h}dx \,ds, \\
=&\half\int_0^t\int_M f P_s^h(\triangle^h g)\ e^{2h} dx \, ds, \\
=&\half \int_0^t\int_M f\bar\triangle^h(P_s^hg)\, e^{2h}dx\, ds  \\
=&\int_0^t\int_M f{\partial \over \partial s}(P_s^h g)\, e^{2h}dx\, ds  \\
=&\int_M f\left (\int_0^t {\partial \over \partial s} (P_s^h g)\, ds\right)\,
 e^{2h} dx,  \\
=&\int_M f(P_t^hg-g)\, e^{2h} dx \\
=&\int_M(P_t^hf-f)g\, e^{2h}dx. 
\end{array} \]
The proof is finished. \hfill \rule{3mm}{3mm}

\bigskip
Next we notice if $M$ is a complete Riemannian manifold of h-finite volume, then $1$ is in the domain of $\bar{\triangle^h}$ (let $h_n$ be the sequence in $C_K^\infty$ approximating the constant function $1$ as in the appendix, then $h_n\to 1$ in $L^2$). 
Thus we have nonexplosion for a h-Brownian motion in this case as is well known (see \cite{GA59} for the case of $h=0$). A quick proof\index{complete} \label{nonexplosion finite volume} is as follows:
$${\partial P_t^h 1\over \partial t}=\half\bar\triangle^h (P_t^h 1) 
=\half P_t^h(\bar \triangle^h 1)=0.$$

 We also have:
\begin{theorem}\label{th: invariant measure}
Let $M$ be a complete\index{invariant measure} Riemannian manifold, then
 $e^{2h}dx$ is an invariant measure\index{invariant measure} for $P_t^h$, i.e.
\begin{equation}
\int P_t^h f(x) e^{2h}dx =\int f(x)e^{2h} dx         \label{eq:invariant}
\end{equation}
for each $L^1$ function $f$.
\end{theorem}

\noindent{\bf Proof:}
First notice if the invariance property $~\ref{eq:invariant}$ holds for
 functions in $C_K^\infty$, it holds for all $L^1$ functions since $P_t^h$ is
 continuous on  $L^1$.

Next there is the following divergence theorem: let $\phi$ be a $L^1$ 1-form
with $\delta^h \phi$ also in $L^1$, then

$$\int_M \left(\delta^h\phi\right) e^{2h} dx=0$$

\noindent
which can be proved by taking approximations of $L^1$ forms by smooth
 compactly supported forms and using the Green theorem on
 page \pageref{eq: divergence theorem}, as proved in \cite{Gromov} for $h=0$. 

Let  $f\in C_K^\infty$, then $\int_0^t P_s^hf ds$
is in the domain of the h-Laplacian $\bar \triangle^h$ by the previous lemma
 and  we have:

\begin{eqnarray*}
\int_M (P_t^h f-f)e^{2h} dx
 &=& \half\int_M \bar \triangle^h (\int_0^t P_s^h f ds) e^{2h} dx=0.
\end{eqnarray*}
 \hfill \rule{3mm}{3mm}

\bigskip

\vspace{6.0mm}

Recall for an elliptic system our semigroup  $P_t$  has the strong Feller property, i.e. it sends $B(M)$ to $C(M)$.  See \cite{MolFeller}. For  such processes there are several notions of recurrence. The basic definition is as follows:
\begin{definition} Denote by $P^x$ the probability law of a process $X$ starting from $x$. The process $X$   is called recurrent if for each $x\in M$, the trajectories of $X$ return  $P^x$ a.s. infinitely often to any given open set in $M$. It is called transient if for each $x\in M$, the trajectories $X$ tend  $P^x$ 
a.s. to infinity as $t\to \infty$. 
\end{definition}

Let $G$ be the potential  kernel of the differential operator
 $\A=\half \triangle +L_Z$:
\begin{equation}
Gf(x)=\int_0^\infty P_tf(x) dt.  \label{eq: potential kernel}
\end{equation}
Then in potential language, $X$ is transient if and only if $G$ is everywhere finite on compact sets, i.e. when $G$ applied to $\chi_K$ for $K$ compact is finite. It 
is recurrent if and only if $G$ is identically infinite on open sets. So a strong 
Feller process is either transient or recurrent  as proved in \cite{AZ74} following 
\cite{AZ-DU-RE}.

\bigskip
\begin{theorem}
Let $X_t$ be an h-Brownian motion on a complete Riemannian manifold  of finite
 h-volume, then $X_t$ is recurrent.
\end{theorem}
We only need to show that it is not transient. Assume $x_t$ is transient. Let
 $K$ be a compact set, then 
$\lim_{t\to \infty} \chi_K(x_t)=0$ a.s. So
$\lim_t E\chi_K(x_t)=0$ by the dominated convergence theorem and
 $\int_M E\chi_K(x_t) e^{2h} dx \to 0$. But
 $\int_M E\chi_K(x_t) e^{2h} dx=$h-vol(K). This gives a contradiction.
 \hfill \rule{3mm}{3mm}

\bigskip
\begin{proposition} \label{pr: ergodic property}
Let $M$ be a complete Riemannian manifold of h-finite volume. Let $\mu$ be the invariant probability measure, then we have for any compact set $K$:
$$\lim_{t\to \infty} P_t^h(\chi_K(x))
={\hbox{h-vol}(K) \over \hbox{h-vol}(M) }=\mu(K).$$
Here "vol" denotes the h-volume.
\end{proposition}

{\bf Proof:}
First notice $\chi_K \in L^2$, so $P_t^h\,\chi_K$ converges in $L^2$ to a harmonic function. The convergence is also in $L^1$ since $\hbox{h-vol}(M)<\infty$, and  also 
the limit function is a constant. So

\begin{eqnarray*}
\lim_{t \to \infty} \int_M P_t^h\chi_K \, e^{2h}dx 
&=& \int_M \lim_{t\to \infty} P_t^h\chi_K e^{2h}dx \\
&=& (\lim_{t\to \infty} P_t^h\chi_K) \hbox{h-vol}(M).
\end{eqnarray*}
But 
$$\int_M P_t^h\chi_K e^{2h} dx =\int_M\chi_K e^{2h} dx=\hbox{h-vol}(K).$$
Thus
$$\lim_{t\to\infty} P_t^h\chi_K ={ \hbox{h-vol}(K) \over \hbox{h-vol}(M)}.$$

\part{Completeness and properties at infinity}

\chapter{Properties at infinity of diffusion semigroups and stochastic flows via  weak uniform covers}
\section{Introduction}

{\bf I.  Background}

A diffusion process is said to be a $C_0$ diffusion\index{$C_0$ diffusion} if 
its semigroup leaves invariant $C_0(M)$\index{$C_0(M)$}, the space of 
continuous functions vanishing at infinity, in which case the semigroup is said
to have the $C_0$ property\index{$C_0$ property}. A Riemannian manifold is said
 to be stochastically complete if the Brownian motion on it is complete, it is
 also said to have the $C_0$ property if the Brownian motion on it does. One 
example of a  Riemannian manifold which is stochastically complete is a 
complete manifold  with finite volume. See Gaffney \cite {GA59}.  More 
generally a complete Riemannian manifold with Ricci curvature bounded from 
below  is stochastically complete and has the $C_0$ property as proved by Yau 
\cite{Yau78}. See also  Ichihara \cite{IC82II}, Dodziuk \cite{Dodz}, Karp-P. Li
 \cite{KA-LI}, Bakry \cite{BA86},  Grigory\'an \cite{GR87}, Hsu\cite{Hsu89}, 
Davies \cite{DA89}, and Takeda \cite{TA91}  for further discussions in terms of
volume growth and bounds on Ricci curvature. For discussions on the behaviour
 at infinity of diffusion processes, and the $C_0$ property, we refer the
 reader to Azencott \cite{AZ74} and Elworthy \cite{EL82}. But we would like to
 mention that a flow consisting of diffeomorphisms (c.f. page 
~\pageref{diffeomorphism}) has the $C_0$ property by arguing by contradiction 
as in \cite{EL82}.

Those papers above are  on a Riemannian manifold except for the last reference. 
For a manifold without a Riemannian structure, Elworthy \cite{ELbook} following
 It$\hat o$ \cite{ITO50}  showed that the diffusion  solution to
 (\ref{eq: basic}) does not explode if there is a uniform cover for the 
coefficients of the equation.
 See also Clark\cite{Clark73}. In particular this shows that the s.d.e
 (\ref{eq: basic}) does not explode on a
 compact manifold if the coefficients are reasonably smooth. See  
\cite{CA-EL83}, \cite{ELflour}.  To apply this method to check whether a 
Riemannian manifold is stochastically complete, we usually construct a 
stochastic differential equation whose solution is Brownian motion.

\noindent
{\bf II.  Main results} 

 The main aim \cite{LI89}  of this chapter is to give unified treatment to some
 of the results from H. Donnelly-P. Li and L. Schwartz.  It gives the first a 
probabilistic interpretation and extends part of the latter. We first introduce
 weak uniform covers in an analogous way to  uniform covers, which gives 
nonexplosion test by using estimations on exit times of the diffusion
 considered. As a corollary this gives the known result on nonexplosion of a 
Brownian motion on a complete Riemannian manifold with Ricci curvature going to
negative  infinity at most quadratically in the distance function\cite{IC82II}.

One interesting example is that a solution to a stochastic differential equation on $R^n$ whose coefficients have linear growth has no explosion and has the $C_0$ property. Notice under this condition,
its associated generator has quadratic growth. On the other hand let $M=R^n$,
 and let $L$ be an elliptic  differential operator :
$$L=\sum_{i,j} a_{ij}{{\partial^2}\over {\partial x_i\partial x_j}}
       +\sum_i b_i{\partial \over {\partial x_i}},$$
where $a_{ij}$ and $b_i$ are $C^2$. Let $(s_{ij})$ be the positive square root of 
the matrix $(a_{ij})$. Let
  $X^i=\sum_j s_{ij}{\partial \over \partial x_j}$,
$A=\sum_j b_j {\partial \over {\partial x_j}}$.  Then 
the s.d.e. defined by:
$$(It\widehat o)  \hskip 20pt    dx_t=\sum_iX^i(x_t)dB_t^i +A(x_t)dt$$
has generator $L$. Furthermore if $|(a_{ij})|$ has quadratic growth and $b_i$ has 
linear growth,  then both $X$ and $A$ in the s.d.e. above have linear growth.
In this case any solution $u_t$ to the following partial differential equation
\label{pde connection}:

$${{\partial u_t}\over {\partial t}} =Lu_t $$ 

\noindent 
satisfies : $u_t\in C_0(M)$, if $u_0\in C_0(M)$ (see next part).

\noindent

\noindent
{\bf III.  Heat equations, semigroups, and flows}

Let $\bar M$ be a compactification of $M$, i.e. a compact Hausdorff space 
which contains $M$ as a dense subset. We assume $\bar M$ is first countable. 
Let $h$ be a  continuous function  on $\bar M$. Consider the following  heat equation with initial boundary conditions:  
\begin{equation}
 \frac {\partial f} {\partial t}={1 \over 2} \Delta f, \hskip 12pt  x\in M , t>0
\label{heat equation}
\end{equation}
\begin{equation}  f (x,0)=h(x), x\in\bar M
\end{equation}
\begin{equation}   f(x,t)=h(x), x\in \partial M.
\label{boundary}
\end{equation}

It is known that there is a unique minimal solution to the first two equations on a 
stochastically complete manifold, the solution is in fact given by the semigroup
associated with Brownian motion on the manifold. So the above equation is not
solvable in general. However with a condition imposed on the boundary of the
 compactification, Donnelly-Li \cite{DO-LI84} showed that the heat semigroup satisfies (~\ref{boundary}). Here is the condition and the theorem:

\bigskip

\noindent {\bf The ball convergence\index{ball convergence criterion} criterion:} 
Let $\{x_n\}$ be a sequence in $M$ converging to a point $\bar x$ on the 
boundary,  then the geodesic balls $\B_r(x_n)$, centered at $x_n$ radius $r$,
  converge to  $\bar x$  as $n$ goes to infinity for each fixed $r$.  

An example of manifolds which satisfies the ball convergence criterion is $R^n$
with sphere at infinity. However this is not true if $R^n$ is given  the 
compactification of  a cylinder with a circle at infinity  added at each end. 
The one point compactification also satisfies the 
ball convergence  criterion. Another class of examples is manifolds with their
geometric compactifications and with the cone topology, i.e. the boundary of
 the manifold are equivalent classes of geodesic rays. See
 \cite{Eberlein-Oneil}. Recall  that two geodesic rays 
$\{\gamma_1(t), t\ge 0\}$ and $\{\gamma_2(t), t\ge 0\}$  are said to be 
equivalent if the the  Riemannian  distance $d(\gamma_1(t),\gamma_2(t))$  
between the two points $\gamma_1(t)$  and $\gamma_2(t)$ is  smaller than a 
constant for each $t$. 

\begin{theorem}[H.Donnelly-P. Li] 
Let $M$ be a complete Riemannian manifold with  Ricci curvature bounded from 
below. The over determined equation (~\ref{heat equation})-(~\ref{boundary}) is
 solvable for any given continuous function  $h$ on $\bar M$,
if and only if the ball convergence criterion holds for $\bar M$.
\end{theorem}

Notice that if the Brownian motion  starting  from $x$ converges to the same
point on the boundary to which  $x$ converges, then
 (~\ref{heat equation})-(~\ref{boundary}) is clearly solvable. See section 4 
for details. We would also like to consider the opposite question: Do we
 get any information on the diffusion processes if we know the behaviour at 
infinity of the associated semigroups? This is true for many cases. In 
particular, for the one point compactification, Schwartz has the following 
theorem\cite{SC89}, which provides a partial converse to \cite{EL82}:

\begin{theorem}[L. Schwartz]
 Let $\bar F_t$ be the standard extension of $F_t$ to $\bar M=M\cup \{\infty\}$, the 
one point compactification. Then the map $(t,x)\mapsto \bar F_t(x)$ is continuous 
from $R_+\times \bar M $ to $L^0(\Omega,P, \bar M)$, the space of measurable maps 
with topology induced from convergence in probability,   if and only if the Semigroup $P_t$  has the $C_0$ property and the map $t\mapsto P_tf$ is continuous from $R_+$ to $C_0(M)$. 
\end{theorem}

\section{Weak uniform covers and non-explosion}

\begin{definition} \cite{ELbook} A stochastic dynamical system (1) is said to admit a
 uniform\index{uniform cover} cover( radius $r>0$, bound $k$), if there are charts $\{\phi_i, U_i\}$ of 
diffeomorphisms from open sets $U_i$ of the manifold onto open sets $\phi_i(U_i)$ of
 $R^n$,  such that:
\begin{enumerate}
\item $B_{3r}\subset \phi_i(U_i)$, each $i$. ($B_\alpha$ denotes the open ball about 
$0$, radius $\alpha$).

\item The open sets $\{\phi_i^{-1}(B_r)\}$ cover the manifold.

\item  If $(\phi_i)_*(X)$ is given by:
$$(\phi_i)_*(X)(v)(e)=(D\phi_i)_{\phi_i^{-1}(v)}X(\phi_i^{-1}v)(e)$$
with  $(\phi_i)_*(A)$ similarly  defined, then both  $(\phi_i)_*(X)$ 
and  $\A(\phi_i)$   are bounded by $k$ on $B_{2r}$. Here $\A$ is the generator 
for the  dynamical system. 
\end{enumerate}
\end{definition}

\begin{figure}
\vspace{40mm}
\end{figure}

Let $M=R^n$. Equation (\ref{eq: basic}) can be interpreted as It$\widehat o$ 
 integral.

\begin{definition} 
	A diffusion process $\{F_t,\xi\}$ is said to have a 
weak\index{weak uniform cover} uniform cover if 
there are pairs of connected open sets $\{U_n^0, U_n\}_{n=1}^\infty$,  and a 
non-increasing sequence  $\{\delta_n\}$ with $\delta_n>0$, such that:
\begin{enumerate}
\item 
$U_n^0\subset U_n$, and the open sets  $\{U_n^0\}$ cover the manifold.
 For $x\in U_n^0$ denote by  $\tau^n(x)$ the first exit time of $F_t(x)$ from 
		the open set $U_n$. Assume $\tau^n<\xi$ a.s. unless $\tau^n=\infty$ almost surely.
\item  
  There exists $\{K_n\}_{n=1}^\infty$, a family of increasing open subsets of
		 $M$ with  $\cup K_n=M$, such that each $U_n$ is contained in one of
 these sets and  intersects at most one boundary from
 $\{\partial K_m\}_{m=1}^\infty$.
\item 
 Let $x\in U_n^0$ and  $U_n\subset K_m$, then for $t<\delta_m$: 
\begin{equation}
 P\{\omega\colon \tau^{n}(x)< t\}\le Ct^2.
\label{in the definition of the weak cover}
\end{equation}
\item 
$\sum_{n=1}^\infty  \delta_n=\infty$. 
\end{enumerate}
\end{definition}

\begin{figure}
\vspace{60mm}
\end{figure}

Notice the introduction of $\{K_n\}$ is only for giving an order to the open 
sets $\{U_n\}$. This is  quite natural when looking at concrete manifolds. In a
 sense the condition says  the geometry of the manifold under consideration 
changes slowly as far as the diffusion process is concerned. In particular if 
the  $\delta_n$ can be taken all  equal, we take $K_n=M$;  e.g. when the number
 of open sets in the cover is  finite. On a Riemannian manifold the open sets
 are often taken as geodesic balls.

\noindent
\begin{lemma} Assume there is a uniform cover for the stochastic dynamic system(1), 
then the solution has a weak uniform cover with $\delta_n=1$, all $n$. 
\end{lemma}
{\bf Proof:} This comes directly from lemma 5 , Page 127 in  \cite{ELbook}.
\hfill \rule{3mm}{3mm}

An example of a stochastic differential equation which has a weak uniform cover
 but not a uniform cover is given by  example 0 on page 50 and the example on
 page  $~\pageref{ex: nonexplosion fo Ricci}$.

\noindent
{\bf Remarks:}
\begin{enumerate}
\item
 If $T$ is a stopping time, then the inequality
 (~\ref{in the definition of the weak cover}) gives the following  from the strong 
Markov property of the process.: 

let $V\subset U_n^0$, and $V\subset K_m$, then when $t<\delta_m$:
\begin{equation}
P\{\tau^n(F_T(x))< t| F_T(x)\in V\}\leq Ct^2 ,
\label{for the weak cover}
\end{equation}
since 
\begin{eqnarray*}
&P&\{\tau^n(F_T(x))<t, F_T(x)\in V\}\\
&=&\int_{V} P(\tau^n(y)<t)P_T(x,dy)
\leq Ct^2   P\{F_T(x)\in V\},
\end{eqnarray*}

here $P_T(x,dy)$ denotes the distribution of $F_T(x)$.
\item
Denote by $P_t^{U_n}$ the heat solution on $U_n$ with Dirichlet boundary condition, then (~\ref{in the definition of the weak cover}) is equivalent to the 
following: when $x\in U_n^0$, 
\begin{equation}
1-P_t^{U_n}(1)(x)\leq Ct^2.
\end{equation}
\item
The methods in the article work in infinite dimensions to give analogous results
(when a Riemannian metric in not needed).
 \end{enumerate}

\noindent
{\bf Exit times:\, }  Given such a cover, let $x\in U_n^0$. We define stopping times $\{T_k(x)\}$ as
follows:  Let $T_0=0$. Let $T_1(x)=\inf_t\{F_t(x,\omega)\not \in U_n\}$ be the 
first exit time of $F_t(x)$ from the set $U_n$. Then  $F_{T_1}(x,\omega)$ must be 
in one of the open sets $\{U_k^0\}$. Let
\begin{eqnarray*}
&\Omega_1^1&=\{\omega: F_{T_1}(x)(\omega)\in U_1^0, T_1(x,\omega)<\infty\},\\
&\Omega_k^1&=\{ \omega: F_{T_1}(x)\in U_k^0- \bigcup_{j=1}^{k-1} U_j^0, T_1<\infty\}.
\end{eqnarray*}
Then $\{\Omega_k^1\}$ are disjoint sets such that $\cup \Omega_k^1=\{T_1<\infty\}$.
In general we only need to consider the nonempty sets of such.  Define further 
the following:
 Let $T_2=\infty$, if $T_1=\infty$. Otherwise if $\omega\in\Omega_k^1$, let: 
\begin{equation}
T_2(x,\omega)=T_1(x,\omega)+\tau^k(F_{T_1}(x,\omega)).
\end{equation}

In a similar way, the whole sequence of stopping times $\{T_j(x)\}$ and sets 
$\{\Omega_k^j\}_{k=1}^\infty$ are defined for $j=3, 4, \dots$. Clearly $\Omega_k^j$
 so defined is measurable with respect to the sub-algebra $\F_j$.

\begin{lemma}
Given a weak uniform cover as above. Let $x\in U_n^0$ and $U_n\subset K_m$. Let 
 $t< \delta_{m+k}$. Then
\begin{equation}
P\{\omega: T_k(x,\omega)-T_{k-1}(x,\omega)< t, T_{k-1}<\infty\}\leq Ct^2.
 \end{equation}  
\end{lemma}

\noindent {\bf Proof:} Notice for such  $x$, $F_{T_{k}}(x)\in K_{m+k-1}$. Therefore for $t< \delta_{m+k}$ we have: 
\begin{eqnarray*}
&P&\{\omega: T_k(x,\omega)-T_{k-1}(x,\omega)< t, T_{k-1}<\infty\}\\
&=& \sum_{j=1}^\infty P(\{\omega: T_k(x,\omega)-T_{k-1}(x,\omega)< t\}
\cap \Omega^{k-1}_j)\\
&=& \sum_{j=1}^\infty P\{\tau^j(F_{T_{k-1}}(x))<t, \Omega^{k-1}_j\}\\
&\leq& Ct^2\sum_{j=1}^\infty P(\Omega^{k-1}_j)\leq Ct^2, 
\end{eqnarray*}           
as in remark 1. Here $\chi_A$ is the characteristic function for a measurable set $A$, and  $E$ denotes taking expectation.  \hfill \rule{3mm}{3mm}

\begin{lemma} If $\sum_n t_n=\infty$, $t_n>0$  non-increasing. Then there is a
non-increasing sequence $\{s_n\}$, such that $0<s_n\le t_n$:

(i) $\hskip 12pt \sum s_n=\infty$

(ii) $\hskip 12pt \sum s_n^2 <\infty$
\end{lemma}
{\bf Proof:}
 Assume $t_n\le 1$, all $n$. Group the sequence $\{t_n\}$ in the following way: 
$$t_1;\ t_2,\ \dots ,\ t_{k_2};\ t_{k_2+1},\ \dots , \ t_{k_3}; \ t_{k_3+1}\ \dots$$
such that $1\le t_2+\dots+t_{k_2}\le 2$,\ $1\le t_{k_i+1} + \dots +t_{k_{i+1}}\le 2$, \ $i\ge 2$.
\ Let $s_1=t_1$, $s_2=\frac {t_2} 2, \dots $, $s_{k_2}=\frac {t_{k_2}} 2$, $s_{k_2+1}=\frac {t_{k_2+1}} 3,\dots$, $s_{k_3}=\frac {t_{k_3}} 3$, $s_{k_3+1}=\frac {t_{k_3+1}} 4, \dots$.  \ Clearly the $s_n$'s so defined satisfy the requirements. $ \hfill \rule{3mm}{3mm}$

So without losing generality,  we may assume from now on that the constants 
$\{\delta_n\}$ for a weak uniform cover fulfill the two conditions in the 
above lemma. With these established, we can now state the nonexplosion result.
The proof is  analogous  to that  of theorem 6  on Page 129 
 in \cite{ELbook}.

\begin{theorem} 
If\index{nonexplosion: theorem} the\index{complete: theorem}
 solution $F_t(x)$ of the equation $(1)$ has a weak uniform  cover,  then it is 
 complete(nonexplosion).
\label{th: complete from weak}
\end{theorem}

\noindent{\bf Proof: \ }
Let $x\in K_n$, $t>0$, $0<\epsilon< 1$. Pick up a number $p$ (possibly  depending
 on $\epsilon$ and $n$), such that $\sum_{i=n+1}^{n+p}\epsilon\delta_i>t$. This is
 possible since  $\sum_{i=1}^\infty \delta_i=\infty$.
So
\begin{eqnarray*}
 P\{\xi(x)< t\}&\leq&     P\{ T_p(x)< t, T_{p-1}<\infty\}\\
&=&    P\{\sum_{k=1}^p (T_k(x)-T_{k-1}(x))<t, T_{p-1}<\infty\}\\
&\leq& \sum_{k=1}^p P\{ T_k(x)-T_{k-1}(x)< \epsilon \delta_{n+k}, T_{k-1}<\infty\} \\
&\leq& Ct^2\epsilon^2\sum_{k=n}^{n+p} \delta_k^2   
\leq Ct^2\epsilon^2\sum_{k=1}^{\infty} \delta_k^2.
\end{eqnarray*}
Let $\epsilon\to 0$, we get:  $P\{ \xi(x)<t\}=0$. 
\hfill \rule{3mm}{3mm}

\noindent
{\bf Remark:}
The argument in the above proof is valid if the definition of a weak uniform
cover is changed slightly, i.e. replacing the constant $C$ by $C_n$ (with
some slow growth condition, say ${1\over p^2}\sum_{j=1}^p C_{n+j}$  is bounded 
for all n) but keep all $\delta_n$ equal.

As a corollary, we have the following known result:

\begin{corollary}\label{co: wuc 2}
Let $\xi$ be the explosion time. If $P\{\xi<t_0\}=0$, for some $t_0>0$. Then there is no explosion.
\end{corollary}
Prove by induction.
In the following we look into some examples:

\noindent{\bf Example 0:}
The flow $F_t(x)=x+B_t$ on $R^n-\{0\}$ does not have a uniform cover. The
 problem occurs at  the origin. But it does have a weak uniform cover as
 constructed below. First note  that we only need to worry about the  origin. 
Take
$U_n=\{x: |x|\le a_n\}$, for $a_n=\left({1\over n+1}\right)^{1\over 8}$.
 Let $U_n^0=U_{n-1}$,
$K_n=\bar U_n$, and $C_n=k \sqrt n$. Here $k$ is a constant.  Now
\begin{eqnarray*}
P\{\tau_n<t\}&\le& P\{\sup_{s\le t}|B_s|\le a_{n}-a_{n+1}\}\\
&\le &{kt^2 \over [a_{n}-a_{n+1}]^4}\le C_nt^2
\end{eqnarray*}
by the maximal inequality. Then
$$\lim_{p\to \infty} {1\over p^2} \sum_1^{p} C_{n+j}=0.$$
See the remark above to see $\{U_n\}$ is a "weak uniform cover".

\bigskip

\noindent
{\bf Example 1 :}  Let $\{U_n\}$ be a family of relatively compact open (proper) 
subsets of $M$ such that $U_n\subset U_{n+1}$ and $\cup_{n=1}^\infty U_n=M$. Assume
 there is a 
sequence of numbers $\{\delta_n\} $ with $\sum_n \delta_n=\infty$, such that the 
following inequality holds  when $t<\delta_{n-1}$ and $x\in U_{n-1}$:
$P\{\tau^{U_n}(x)<t\}\leq ct^2$. 
Then the diffusion concerned does not explode by taking
 $\{U_{n+1}-\overline{U_{n-1}},U_n-\overline{U_{n-1}}\}$ to be a weak uniform cover and $K_n=\bar U_n$.

\section{Boundary behaviour of diffusion processes}

To consider the boundary behaviour of  diffusion processes, we introduce the 
following concept:

\begin{definition} A weak uniform  cover\index{weak uniform cover: regular} $\{U_n^0, U_n\}$ is said to be regular (at 
infinity for $\bar M$), if the following holds: let $\{x_j\}$ be a sequence in $M$ converging to $\bar x\in \partial M$, and
$x_j\in U_{n_j}^0\in\{U_n^0\}_{n=1}^\infty$, then the corresponding open sets 
$\{U_{n_j}\}_{j=1}^\infty  \subset \{U_n\}_1^\infty$ converges to $\bar x$ as well. 
A regular uniform cover can be  defined in a similar way. 
\end{definition}
 
For a point $x$ in $M$, there are a succession of related open sets 
$\{W_x^p\}_{p=1}^\infty$,  which are defined as follows:
Let $W_x^1$ be the union of all open sets from $\{U_n\}$ such that  $U_n^0$
contains $x$, and $W_k^2$ be the union of all open sets from 
$\{U_n\}$ such that $U_n^0$ intersects one of the small balls $U_{n_j}^0$ defining 
$W_x^1$. Similarly $\{W_{k}^{p}\}$ are defined. These sets are well defined
and in fact form an increasing sequence.

\begin{lemma} Assume $F_t$ has a regular (weak) uniform cover. Let $\{x_n\}$ be
 a sequence of $M$ which converges to a point $\bar x \in \partial M$. Then 
$W_{x_n}^p$ converges to $\bar x$ as well for each fixed $p$.
\end{lemma}
{\bf Proof:} Note that by arguing by contradiction, we only need to prove the
 following:
let  $\{z_k\}$ be a sequence  $z_k\in W_{x_k}^p$, then $z_k \to \bar x$, as $n\to \infty$.  First let $p=2$. 

By definition, for each $x_k$, $z_k$, there are open
 sets $U_{n_k}^0$ and $U_{m_k}^0$ such that $x_k\in U_{n_k}^0$, 
  $z_k\in U_{m_k}^0$ and $U^0_{n_k}\cap U^0_{m_k}\not = \emptyset$. Furthermore $U_{n_k}\to \bar x$ as $k\to \infty$. 
Let $\{y_k\}$ be a sequence of points with $y_k\in U_{n_k}^0\cap U_{m_k}^0$.  But $y_k\to \bar x$ since $U_{n_k}$ does.  So $U_{m_k}\to \bar x$ again 
from the definition of a regular weak uniform cover. Therefore $z_k$ converges to $\bar x$ as $k\to \infty$, which is what we want. 
The rest can be proved by induction. \hfill\rule{3mm}{3mm}

\begin{theorem} If the diffusion $F_t$ admits a regular weak uniform cover, 
with $\delta_n=\delta$, all $n$, then the map $F_t(-):M\to M$ can be extended
to the compactification $\bar M$ continuously in probability with the
restriction to the boundary to be the identity map,  uniformly in $t$ in finite
intervals. (We will say $F_t$ extends.)
\label{th: extends}
\end{theorem}

\noindent {\bf Proof:} 
Take $\bar x\in \partial M$ and a sequence $\{x_n\}$ in $M$ converging to
 $\bar x$.  Let $U$ be a neighbourhood of $\bar x$ in $\partial M$. We want to
 prove for each $t$: 
$$\lim_{n\to\infty}P\{\omega:F_s(x_n,\omega)\not\in U, s<t \}=  0.$$
 Since $x_n$ converges to $\bar x$, there is a number $N(p)$ for each $p$, such
 that if $n>N$, $W_{x_n}^p\subset U$. Let  $t>0$, Choose $p$ such that 
${2t \over p}<\delta$.  For such a number $n>N(p)$ fixed,   we have:
\begin{eqnarray*}
 &P&\{\omega: F_s(x_n,\omega)\not \in U, s< t\}\\
&\le& P\{\omega: F_s(x,\omega)\not\in W_{x_n}^p,s< t\}\\
	&\leq& P\{\omega:T_p(x_n)(\omega)< t,T_{p-1}(x_n)<\infty\}\\
	&\leq& \sum_{k=1}^{p} P\{T_k(x_n)-T_{k-1}(x_n)<  {t \over p} , 
		T_{k-1}(x_n)<\infty\}.\\
	&\leq& {Ct^2 \over p}.	\end{eqnarray*}
Here $C$ is the constant in the definition of the weak uniform cover. Let $p$ go 
to infinity to complete the proof. \hfill \rule{3mm}{3mm}

\noindent {\bf Remark:}
If the $\delta_n$ can be taken all equal, theorem ~\ref{th: complete from weak},  theorem ~\ref{th: extends}  hold if (\ref{in the definition of the weak 
cover}) is relaxed to:         $$P\{T^n(x)<t\}\le f(t),$$
for some nonnegative function $f$ satisfying $\lim_{t\to 0} {f(t) \over t}=0$.

\noindent
{\bf Example 2:\ }
Let $M=R^n$ with the one point compactification. Consider the following s.d.e.:
$$ (It\widehat o) \hskip 12pt  dx_t=X(x_t)dB_t+A(x_t)dt.$$ 
Then if both $X$ and $A$ have linear growth, the solution has the $C_0$ property.

\noindent
{\bf Proof:} There is a well known uniform cover for this system. See \cite{Clark73},
 or \cite{ELbook}. 
 A slight change gives us the following regular uniform cover:

Take a countable set of points $\{p_n\}_{n\ge 0}\subset M$ such that the open sets 
 $U_n^0=\{z:|z-p_n|<{ |p_n| \over 3}\}, n=1,2, \dots$  and  $U_0^0=\{ z: |z-p_0|<2\}$
cover $R^n$. Let $U_0=\{ z: |z-p_0|<6\}$; and $U_n=\{z:|z-p_n|<{|p_n|\over 2}\}$, for $n\not = 0$.  Let  $\phi_n$ be the chart map on $U_n$:
$$\phi _n(z)={z-p_n\over |p_n|}.$$
This certainly defines a uniform cover(for details see Example 3 below). Furthermore if $z_n\to \infty$ and
$z_n\in U_n^0$, then any $y\in U_n$ satisfies the following:
 $$|y|>{|p_n| \over 2} >{1\over 3}|z_n| \to \infty,$$
since $|p_n|\ge {{3 |z_n|} \over 4}$.
Thus we have a regular uniform cover which gives the required $C_0$ property. 
$\hfill \rule{3mm}{3mm}$
\noindent

\noindent {\bf Example 3:}  
Let $M=R^n$, compactified with the sphere at infinity: 
 $\bar M=R^n\cup S^{n-1}$. Consider the same s.d.e as in the example above. 
Suppose both $X$ and $A$ have sublinear growth of power $\alpha<1$: 
\begin{eqnarray*}
|X(x)|&\le&   K(|x|^\alpha+1)\\
 |A(x)|& \le& K(|x|^\alpha+1)
\end{eqnarray*}
for a constant $K$. Then there is no explosion. Moreover the solution  $F_t$ extends.

\noindent {\bf Proof:} 
 The proof is as in example 2, we only need to construct a regular
 uniform cover for the s.d.e.:

Take points  $p_0, p_1, p_2, \dots$ in $R^n$ (with $|p_0|=1$, $|p_n|>1$) ,  such that  the open sets $\{U_n^0\}$ defined  by : 
 $ U_0=\{ z: |z-p_0|<2\}$, $U_i^0=\{z:|z-p_i|<{|p_i|^\alpha \over 6}\}$ cover $R^n$. 

 Let $U_0=\{ z: |z-p_0|<6\}$, $U_i=\{z:|z-p_i|<{|p_i|^\alpha\over 2}\}$, and
let   $\phi_i$ be the chart map from $U_i$ to $R^n$:
$$\phi _i(z)={6(z-p_i)\over |p_i|^\alpha}.$$

Then $\{\phi_i, U_i\}$ is a uniform cover for the stochastic dynamical system. 
In fact, for $i\not = 0$, and $y\in B_3\subset R^n$: 
$$ |(\phi_i)_\ast (X)(y)|\le {K(1+|\phi_i^{-1}(y)|^\alpha) \over |p_i|^\alpha}
\leq {K \over |p_i|^\alpha} (1+2|p_i|^\alpha)<18K.$$

\noindent
Similarly $|(\phi_i)_\ast (A)(y)|\le 18K$, and $D^2\phi_i=0$.

Next we show this cover is regular. Take a sequence  $x_k$ converging to  $\bar x$
in $\partial R^n$. Assume $x_k\in U_k^0$. Let $z_k\in U_k$. We want to prove $\{z_k\}$ converge to $\bar x$.  First the norm of $z_k$ converges to infinity as $k \to \infty$, since $|p_k|>{ {2|x_k|} \over 3}$ and $|z_k|>|p_k|-{1\over 2} |p_k|^\alpha$.

\begin{figure}
\vspace{55mm}
\end{figure}

Let $\theta$ be the biggest angle between points in $U_k$, then 

$$     \hbox{tan}{\theta\over 2} \leq \sup_{z\in U_n} { |z-p_n| \over |p_n|}
\le {|p_n|^\alpha \over 2|p_n|}
\le {|p_n|^{\alpha -1}\over 2} \to 0.         $$
Thus $\{U_n, \phi_n\}$ is a uniform cover satisfying the convergence criterion for
the sphere compactification. The required result holds from the theorem. 
\hfill \rule{3mm}{3mm}

This result is sharp in the sense there is a s.d.e. with coefficients having
linear growth but the solution to it does not extend to the sphere at infinity to be identity:

\noindent {\bf Example 4:} 
 Let $B$ be a one dimensional Brownian motion. Consider the 
following s.d.e on the complex plane $\cal C$:
		$$dx_t=ix_tdB_t.$$
The solution starting from $x$ is in fact $xe^{iB_t+{t\over 2}}$, which does not 
continuously extend to be the identity on the sphere at infinity.

\bigskip

\noindent
{\bf Example 5:}
Let $U$ be a  bounded open set of $R^n$. Let $(X,A)$ be a s.d.s. (in It\^o form)
 on $U$  satisfying
$$|X(x)|\le k\ d(x,\partial U),$$
and
$$|A(x)| \le k\ d(x,\partial U)$$
for some constant $k$, then there is no explosion. Here $d$ denotes the distance
function  on $R^n$ and $\partial U$ denotes the boundary of $U$.

\bigskip
\noindent
{\bf Proof:}  Choose points $\{x_n\}$ such that the balls 
 $B(x_n, \half\, d(x_n,\partial U))$ centered at $x_n$ radius
 $\half\, d(x_n,\partial U)$  cover  $U$. 
Define a map:
 $$\phi_n\colon  B(x_n, \half\, d(x_n,\partial U))\to B_3$$
 by
$$\phi_n(x)= 6\times {x-x_n \over d(x_n,\partial U)}.$$
Then 
$$(T\phi_n)_*(X)\le {6|X(\phi_n^{-1}(x))| \over d(x_n,\partial U)}
\le {6kd(\phi_n^{-1}(x),\partial U) \over d(x_n,\partial U)}
\le 12k$$
by the triangle inequality of the distance function. Thus we have a uniform cover and
so nonexplosion.  \hfill \rule{3mm}{3mm}

This result is sharp in the sense that there is an example (in It\^o form) given in
\cite{TANI89} which satisfies:
$${|X(x)|+|A(x)|\over d(x,\partial U)^\epsilon } \le 2$$
for $0<\epsilon <1$, but  has explosion. Here is the example (on $\{|x|<1\}\subset
R^2$):

\begin{eqnarray*}
dx_t^1&=&(1-|x_t|^2)^\epsilon x_t^1 dB_t -\half (1-|x_t|^2)^{2\epsilon} x_t^1\,dt\\
dx_t^2&=&(1-|x_t|^2)^\epsilon x_t^2 dB_t -\half (1-|x_t|^2)^{2\epsilon} x_t^2\,dt.
\end{eqnarray*}
Here $0<\epsilon<1$. See \cite{TANI89} for more discussions on nonexplosion on 
open sets of $R^n$.

\section{Boundary behaviour continued}

A diffusion process is a $C_0$ diffusion\index{$C_0$ diffusion} if its semigroup has the $C_0$ property. This is equivalent to the following \cite{AZ74}: let $K$ be a
 compact set, and $T_K(x)$ the first entrance time to $K$ of the diffusion starting 
from $x$, then  $\lim_{x\to\infty} P\{T_K(x)<t\}=0$ for each $t>0$, and each compact set $K$. 

The following theorem follows from theorem 3.2 when $\delta_n$ in the definition
of  weak uniform cover can be taken all equal:

\begin{theorem} Let\index{$C_0$ diffusion: theorem} $\bar M$ be the one point compactification. Then if the 
diffusion process $F_t(x)$ admits a regular weak uniform cover, it is a 
$C_0$ diffusion.
\end{theorem}

\noindent
{\bf Proof:}
Let $K$ be a compact set with $K\subset K_j$; here $\{K_j\}$ is as in definition 2.2.
Let $\epsilon>0$, $t>0$, then there is a number $N=N(\epsilon, t)$ such that:

$$\delta_{j+2} +\delta_{j+4}+\dots+ \delta_{j+2N-2}>{t\over \epsilon}.$$

Take $x\not \in K_{j+2N}$.
Assume $x\in K_m$, some $m>j+2N$. Let $T_0$ be the first entrance time of 
$F_{t}(x)$ to $K_{j+2N-1}$,   $T_1$ be the first entrance time  of 
$F_{t}(x)$ to $K_{j+2N-3}$ after $T_0$, (if $T_0<\infty$), and so on. 
 But $P\{T_i<t, T_{i-1}<\infty\}\le Ct^2$ 
for $t<\delta_{j+2N-2i}, i>0$, since any open sets from the cover  intersects at most one boundary of sets from $\{K_n\}$. Thus
\begin{eqnarray*}
P\{T_K(x)<t\}&\leq&   P\{\sum_1^{N-1} T_{i}(x)<t, T_{N-2}<\infty\}   \\
&\leq& \sum_{i=1}^{N-1} P\{T_i(x)<\epsilon \delta_{j+2N-2i}, T_{i-1}(x)<\infty\}\\
&\le& C\epsilon^2\sum_1^{N-1}\delta_{j+2N-2i}^2 \, 
\le C\epsilon^2\sum_1^\infty\delta_{j}^2.
\end{eqnarray*}
The proof is complete by letting $\epsilon \to 0$. \hfill \rule{3mm}{3mm}

\noindent
{\bf Example 6:\,} \label{ex: infinity Ricci}
 $\label{ex: nonexplosion fo Ricci}$
 Let $M$ be a complete Riemannian manifold, $p$ a fixed point in $M$.
Denote by $\rho(x)$ the distance between $x$ and $p$,  ${\cal B}_r(x)$ the 
geodesic ball centered at $x$ radius $r$, and Ricci$(x)$ the Ricci curvature at $x$.

\noindent {\bf Assumption A:}  
\begin{equation}
\int_1^\infty {1\over \sqrt{K(r)}} dr=\infty .  
\label{Assumption A} \end{equation}

\noindent
Here $K$ is defined as follows:

$$K(r)=-\{\inf_{B_r(p)} \hbox{Ricci}(x)\wedge 0\}.$$

Let $X_t(x)$ be a Brownian motion on $M$ with $X_0(x)=x$. Consider the first exit 
time of $X_t(x)$ from $B_1(x)$:

 $$T=\inf_t\{t\geq 0: \rho(x ,X_t(x))=1\}.$$ 

Then we have the following  estimate on $T$ from \cite{Hsu89}:

If $L(x)>\sqrt{K(\rho(x)+1)}$, then
$$P\{T(x)\leq {c_1 \over L(x)} \}\leq \exp^{-c_2L(x)}$$
for all $x\in M$. Here $c_1, c_2$ are positive constants independent of $L$. 

\noindent
This can be rephrased into the form we are familiar with:
 when $0\le t<{c_1 \over \sqrt{K(\rho(x)+1)}}$,
$$P\{T(x)\leq t\}\leq \exp^{-{c_1c_2   \over t}}.$$
 But  $\lim_{t\to 0} {e^{-{c_1c_2\over t}} \over t^2}=0$. So there is a 
 $\delta_0>0$, such that:
$e^{-{c_1 c_2 \over t}} \leq t^2$, when $t<\delta_0$.
Thus:
 
\noindent
{\bf Estimation on exit times:  \,} 
when $t<{c_1 \over \sqrt{K(\rho(x)+1)}}\wedge \delta_0 $,
\begin{equation}
P\{T(x)<t\}\leq t^2.
\end{equation}

\noindent
Let $\delta_n={1\over \sqrt{K(3n+1)}}\wedge \delta_0$, then we also have the
 following :
\begin{equation}
\sum_1^\infty \delta_n      \ge   \sum_1^\infty{1\over \sqrt{K(3n+1)}}
\ge \int_1^\infty {1\over \sqrt{K(3r)}}\, dr=\infty.
\end{equation}

With this we may proceed to prove the following from \cite{Hsu89}:

\noindent {\bf Corollary:} [Hsu] \,\,
A complete Riemannian manifold $M$ with Ricci curvature satisfies assumption A is 
stochastically complete and has the $C_0$ property.

\noindent {\bf Proof:} 
There is a regular weak uniform cover as follows:

 First take any $p\in M$, and let $K_n=\overline{B_{3n}(p)}$. 
Take points $p_i$ such that $U_i^0=B_1(p_i)$ covers the manifold. Let $U_i=B_2(p_i)$. Then  $\{U_i^0,U_i\}$ is a regular weak uniform cover for $M\cup \Delta$.  

\noindent{\bf Remark:} Grigory\'an has the following volume growth test on nonexplosion.  The Brownian motion does not explode on a manifold if  
$$\int^\infty {r\over \hbox{Ln}(\hbox{Vol}(B_R))}\, dr=\infty.$$
Here $\hbox{Vol}(B_R)$ denotes the volume of a geodesic ball centered at a point $p$ in  $M$. 
This result is stronger than the corollary obtained above by the following 
comparison theorem on a $n$ dimensional manifold: let $\omega_{n-1}$ denote the 
volume of  the $n-1$ sphere of radius $1$,
 
$$\hbox{Vol}(B_R) \le \omega_{n-1} \int_0^R\{ \sqrt{{(n-1) \over K(R)}} 
\hbox{Sinh}(\sqrt{{K(R)\over (n-1)}} r)\}^{(n-1)}\, dr.$$
Notice $K(R)$ is  positive when $R$ is sufficiently big provided the  Ricci curvature is not nonnegative everywhere. So Grigory\'an's result is stronger than the one obtained above.

\bigskip

The definition of weak uniform cover is especially suitable for the one point 
compactification. For general compactification the following definition explores
more of the geometry of the manifold and gives better result:

\begin{definition}
Let\index{uniform cover: boundary} $\bar M$ be a compactification of $M$, $\bar x \in \partial M$. A diffusion 
process $F_t$ is said to have a uniform cover at point $\bar x$, if there is a 
sequence  $A_n$ of open neighbourhoods of $\bar x$ in $\bar M$ and positive numbers $\delta_n$  and a constant $c>0$, such that:
\begin{enumerate}
\item 
The sequence of $A_n$ is strictly decreasing, with $\cap A_n=\bar x$, and 
$A_n\supset \partial A_{n+1}$. 
\item  
The sequence of numbers $\delta_n$ is non-increasing  with $\sum \delta_n=\infty$
and $\sum_n\delta^2_n<\infty$.
\item When $t<\delta_{n}$, and $x\in A_{n}-A_{n+1}$,
$$P\{\tau^{A_{n-1}}(x)<t\}\leq ct^2.$$
Here $\tau^{A_n}(x)$ denotes the first exit time of $F_t(x)$ from the set $A_n$.
\end{enumerate}
\end{definition}

\begin{figure}
\vspace{60mm}
\end{figure}

\begin{proposition}
If there is a uniform cover for $\bar x \in \partial M$, then $F_t(x)$ converges
to $\bar x$ continuously in probability as $x \to\bar x$ uniformly in $t$ in finite interval.
\end{proposition}

\noindent
{\bf Proof:}
 The existence of $\{A_n\}$ will ensure $F_{\tau^{A_{n}}}(x)\subset A_{n-1}$, 
which allows us to apply a similar argument as in the case of the one point 
compactification. Here we denote by $\tau^A$ the first exit time of the process 
$F_t(x)$ from a set $A$.

Let $U$ be a neighbourhood of $\bar x$. For this  $U$, by compactness of $\bar M$, there is a number $m$ such that $A_m\subset U$, since $\cap_{k=1}^\infty A_k=\bar x$.
Let $0<\epsilon<1$, $\bar \epsilon=({\epsilon \over c\sum_k \delta_k^2})^{1\over 2}$.
we may assume $\bar \epsilon<1$.  Choose $p=p(\epsilon)>0$ such that: 

$$\delta_{m}+\delta_{m+1}+\dots +\delta_{m+p-1}>{t \over\bar \epsilon}.$$

Let $x\in A_{m+p+2}$. Denote by $T_0(x)$ the first exit time of $F_t(x)$ from
 $A_{m+p+1}$, $T_1(x)$ the first exit time of $F_{T_0}(x)$ from $A_{m+p}$ where 
defined.  Similarly $T_i, i>1$ are defined.

Notice  if $T_i(x)<\infty$, then $F_{T_i(x)}\in A_{m+p-i}-A_{m+p+1-i}$, for $i=0, 1, 2,\dots $. Thus for $i>0$ there is the following inequality from the definition and the Markov property:
$$P\{T_i(x)<\bar\epsilon \delta_{m+p-i}\}\leq c\bar\epsilon^2 \delta^2_{m+p-i}.$$

Therefore we have:
\begin{eqnarray*}
P\{\tau^U(x)<t\}    &\leq&  P\{\tau^{A_m}(x)<t\}     \\
&\leq&P\{T_{p}+\dots +T_1<t, T_{p-1}<\infty\}\\
&\leq& \sum_{i=1}^{p} P\{T_i<\bar \epsilon\delta_{m+p-i}, T_{i-1}<\infty\}\\
&\leq& c\bar\epsilon^2\sum_{i=1}^{p}\delta_{m+p-i}^2<\epsilon .    
\end{eqnarray*}
The proof is finished.   \hfill \rule{3mm}{3mm}

\section{Properties at infinity of semigroups}

Recall a semigroup is said to have  the $C_0$ property, if it sends 
$C_0(M)$, the space of continuous functions on $M$ vanishing at infinity, to itself. 
Let $\bar M$ be a compactification of $M$. Denote by $\Delta$ the  point at infinity
for the one point compactification.  Corresponding to the $C_0$ property of
 semigroups we consider the following $C_*$ property for
 $\bar M$: 

\begin{definition}
A semigroup $P_t$ is said to have the $C_*\hskip 4pt$ property\index{$C_*$ property} for $\bar M$,  if for each continuous function $f$ on $\bar M$, the following holds: let $\{x_n\}$ be a sequence converging to $\bar x$ in $\partial M$, then
\begin{equation}
	\lim_{n\to \infty} P_tf(x_n)=f(\bar x), 
\end{equation}
\end{definition}

To justify the definition, we notice if $\bar M$ is the one point compactification,
condition $C_*\hskip 4pt$ will imply the $C_0\hskip 4pt$ property of the semigroup.
On the other hand if $P_t$ has the $C_0$ property, it has the $C_*$ property for
 $M\cup \Delta$ assuming nonexplosion. This is observed by subtracting a constant
function from a continuous function $f$ on $M\cup\Delta$: Let $g(x)=f(x)-f(\Delta)$,
then $g\in C_0(M)$. So $P_tg(x)=P_tf(x)-f(\Delta)$. Thus 
$$\lim_{n \to \infty} P_tf(x_n)=\lim_{n \to \infty} P_tg(x_n) +f(\Delta)=f(\Delta),$$ 
\noindent
if $\lim_{n \to \infty} x_n=\Delta$. 

In fact the $C_*$ property holds for the one point compactification if and only if
there is no explosion and the $C_0$ property holds. These properties are often 
possessed by processes, e.g. a Brownian motion on a Riemannian manifold with
 Ricci curvature which satisfies (~\ref{Assumption A}) has this property.

Before proving this claim, we observe first that:
\begin{lemma} If $P_t$ has the $C_*$ property for any compactification $\bar M$, 
it must have the $C_*$ property for the one point compactification. 
\end{lemma}

\noindent {\bf Proof:\ }
Let $f\in C(M\cup \Delta)$. Define a map $\beta$ from $\bar M$ to $M\cup \Delta$:
$\beta(x)=x $ on  the interior of $M$, and $\beta(x)=\Delta$, if $x$ belongs to 
the boundary. Then $\beta$ is a continuous map from $\bar M$ to $M\cup \Delta$,
since for any compact set $K$, the inverse set $\bar M-K=\beta^{-1}(M\cup\Delta-K)$ 
is open in $\bar M$.

Let $g$ be the composition map of $f$ with $\beta$:
$g=f\circ \beta: \bar M\to R$.
 Thus $g(x)|M=f(x)|M$,  and $g(x)|_{\partial M}=f(\Delta)$.
So for a sequence $\{x_n\}$ converging to $\bar x\in \partial M$,
$\lim_n P_tf(x_n)=\lim_n P_tg(x_n)=g(\bar x)=f(\Delta)$.                  $\hfill \rule{3mm}{3mm}$   

We are ready to prove the following theorem:

\begin{theorem}
If\index{complete: theorem} a semigroup $P_t$ has the $C_* $ property, the associated diffusion process $F_t$ is complete\index{complete: theorem}. 
\end{theorem}

\noindent{\bf Proof:} 
  We may assume the compactification under consideration is the one point
compactification from the the lemma above. Take $f\equiv 1$, $P_tf(x)=P\{t<\xi(x)\}$.
But $P\{t<\xi(x)\}\to 1$ as $x\to \Delta$ from the assumption. 
More precisely for any $\epsilon>0$, there is a compact set $K_\epsilon$ such that
if $x\not\in K_\epsilon$, $P\{t<\xi(x)\}>1-\epsilon$.

Let $K$ be a compact set containing $K_\epsilon$. Denote by $\tau $ the first
exit time of $F_t(x)$ from $K$. So $F_\tau(x)\not \in K_\epsilon$ on the set 
$\{\tau<\infty\}$. Thus:
\begin{eqnarray*}
P\{t<\xi(x)\}
 &\geq& P\{\tau<\infty, t<\xi(F_\tau(x))\}+P\{\tau=\infty\}\\
 &=& E\{\chi_{\tau<\infty} E\{\chi_{t<\xi(F_\tau(x))}|\F_\tau\}\}+  P\{\tau=\infty\}.
\end{eqnarray*}
Here $\chi_A$ denotes the characteristic function of set $A$. Applying the strong Markov property of the diffusion we have:
$$E\{\chi_{\tau<\infty} E\{\chi_{t<\xi(F_\tau(x))}|\F_\tau\}\}
=E\{\chi_{\tau<\infty} E\{t<\xi(y)|F_\tau=y\}\}.$$
However $$E\{\chi_{t<\xi(y)}|F_\tau=y\}> 1-\epsilon,$$
So
\begin{eqnarray*}
P\{t<\xi(x)\}&\ge& P\{\tau=\infty\}+E\{\chi_{\tau<\infty}(1-\epsilon)\}  \\
&=& 1-\epsilon P\{\tau<\infty\}.
\end{eqnarray*}
Therefore  $P\{t<\xi(x)\}=1$, since $\epsilon$ is arbitrary. $ \hfill \rule{3mm}{3mm}$

In the following we examine the  the relation between the behaviour at $\infty$
of a diffusion process  and the diffusion semigroup.

\begin{theorem}The semigroup $P_t$ has the $C_*$ property if and only if the
 diffusion process
$F_t$ is complete and can be extended to $\bar M$ continuously in probability
with $F_t(x)|_{\partial M}=x$.
\end{theorem}

\noindent{\bf Proof:}
        Assume $F_t$ is complete and extends. 
	Take a point $\bar x \in \bar M$, and sequence $\{x_n\}$ converging to
 $\bar x$. Thus 
$$\lim_{n\to \infty} P_tf(x_n)=\lim_{n\to \infty} Ef(F_t(x_n))=Ef(\bar x)
=f(\bar x)$$
 for any continuous function on $\bar M$, by the dominated convergence theorem.
	
	On the other hand $P_t$ does not have the $C_*$ property if the 
	assumption above is not true.
	In fact let $x_n$ be a sequence converging to $\bar x$ , such that
 for some  neighbourhood $U$ of $\bar x$, and a number $\delta>0$:
	$$\liminf_{n\to \infty} P\{F_t(x_n)\not \in U\}=\delta.$$
	There is therefore a subsequence  $\{x_{n_i
}\}$ such that:
	$$\lim_{i\to \infty} P\{F_t(x_{n_i})\not \in U\}=\delta.$$
	Thus there exists $N>0$, such that if $i>N$:
	$$P\{F_t(x_{n_i}) \in \bar M-U\}>{\delta \over 2}.$$
	But since $\bar M$ is a compact Hausdorff space, there is a continuous 
function $f$
	from $M$ to $[0,1]$ such that $f|_{\bar M-U}=1$, and $f|_G=0$, for any
 open set $G$ in $U$. Therefore
	
\begin{eqnarray*}
P_tf(x_n)&=&Ef(F_t(x_n))  \\
&\geq& \int_{\{\omega: F_t(x_n)\in \bar M-U\}}
	 f(F_t(x_n))\, P(d\omega)\\
&=&P\{F_t(x_n)\in \bar M-U\} >  {\delta \over 2}.  
\end{eqnarray*}
So $\lim P_tf(x_n)\not =f(\bar x)=0.$ 
 $\hfill \rule{3mm}{3mm}$

As is known, a flow consisting of diffeomorphisms has the $C_0$ property. But
this is , in general, not true for the $C_*$ property. See example 4.
	
\begin{corollary}
Assume the diffusion process $F_t$ admits a weak uniform cover regular for 
$M\cup \Delta$,
then its diffusion semigroup $p_t$ has the $C_*$ property. The same is true 
 for a general compactification if all $\delta_n$ in the weak uniform cover 
can be taken equal.
\end{corollary}

\noindent {\bf Example 7:}\cite{DO-LI84}\,
Let $M$ be a complete connected Riemannian manifold with Ricci curvature bounded from below. Let  $\bar M$ be  a compactification such that the ball convergence criterion holds (ref. section 1).  In particular the over determined equation 
(~\ref{heat equation})-(~\ref{boundary}) is 
solvable  for any continuous function $f$ on $\bar M$ if  the ball convergence 
criterion holds.

\noindent {\bf Proof:} 
We keep the notation of example 6. 
Let $K=-\{\inf_x Ricci(x)\wedge 0\}$, $\delta={c_1\over k}$, where $c_1$ is
the constant in example 6. Let $p\in M$ be a fixed point, and $K_n=\overline{B_{3n}(p)}$ be compact sets in $M$.  Take points $\{p_i\}$ in $M$ such that $\{B_1(p_i)\}$ cover the manifold. 
Then $\{B_1(p_i), B_2(p_i)\}$ is a weak uniform cover from (4.10) and (4.11). 
Moreover 
this is a regular cover if the ball convergence criterion holds for the compactification.

\part{Strong completeness, derivative semigroup, and moment stability}

\chapter{ On the existence of  flows: strong completeness}
\section{Introduction}
A stochastic dynamical system $(X,A)$ is said to be {\it complete} if its explosion
 time  $\xi(x)$ is infinite almost surely for each $x$ in $M$. It is called 
{\it strongly  complete}\index{strong completeness: definition} if there is a
 version of the solution which is jointly continuous in time and space for all time. In this case the solution is called a {\it continuous flow}\index{flow: continuous}. 
Examples of s.d.s. which are complete but not strongly complete can be found
in \cite{ELbook}, \cite{Jakel}, and \cite{Kunitabook}.

The known results on the existence of a continuous flow are concentrated on $R^n$ 
and compact manifolds. On $R^n$ results are given (for It\^o equations) in terms of
 global Lipschitz  or similar conditions. See  Blagovescenskii and Friedlin 
\cite{BL-FR}. The  problems concerning the diffeomorphism property of flows 
have been discussed by e.g. Kunita \cite{KUNITA80}, Carverhill and Elworthy 
\cite{CA-EL83}. See  Taniguchi \cite{TANI89} for discussions on the strong
 completeness of a stochastic dynamical system on an open set of $R^n$.  For 
discussions of higher derivatives 
of solution flows on $R^n$, see Krylov \cite{Krylov} and  Norris \cite{Norr-smc}.

On a compact manifold, a stochastic differential equation  with $C^2$ coefficients 
is strongly complete.  In fact the solution flow is almost as smooth as the
 coefficients of the stochastic differential equation.  More  precisely  the 
solution flow is $C^{r-1}$  if the coefficients are $C^r$. Moreover the flow 
consists of diffeomorphisms. See Kunita   \cite{KUNITA80},  Elworthy \cite{EL78}, 
and Carverhill and Elworthy \cite{CA-EL83}. For discussions in the framework of
diffeomorphism groups see Baxendale \cite{Baxendale80} and Elworthy \cite{ELbook}.

\bigskip

In general  we know very little about strong 
completeness. Our aim is to prove  strong completeness given nonexplosion and
 certain regularity properties of the solution.

\bigskip

To begin with we quote the following theorem on the  existence of a
 {\it partial flow}\index{flow: partial} from
 \cite{ELbook} (first proved in \cite{KUNITA80},  extended later in \cite{ELbook},
 and  \cite{CA-EL83}):

\begin{theorem}\cite{ELbook} $\label{th: partial flow}$
Suppose $X$, and $A$ are  $C^{r}$, for $r\ge 2$. Then there is a partially defined flow $(F_t(\cdot),\xi(\cdot))$ which 
is a maximal solution  to $(\ref{eq: basic})$ such that if 

$$M_t(\omega)=\{x\in M, t<\xi(x,\omega)\},$$

\noindent
then there is a set $\Omega_0$ of full measure such that for all
 $\omega\in \Omega_0$:
\begin{enumerate}
\item
$M_t(\omega)$ is open in $M$ for each $t>0$, i.e. $\xi(\cdot,\omega)$ is lower semicontinuous.
\item
$F_t(\cdot,\omega): M_t(\omega)\to M$ is in $C^{r-1}$ and is a diffeomorphism onto an open
subset of $M$. Moreover the map : $t\mapsto F_t(\cdot,\omega)$ is continuous into
 $C^{r-1}(M_t(\omega))$, with the topology of uniform convergence on compacta of the first r-1 derivatives.
\item
Let $K$ be a compact set and $\xi^K=\inf_{x\in K} \xi(x)$. Then

\begin{equation}
\lim_{t\nearrow \xi^K(\omega)} \sup_{x\in K} d(x_0, F_t(x))=\infty
\end{equation}

\noindent
almost surely on the set $\{\xi^K<\infty\}$. (Here $x_0$ is a fixed point of $M$ and $d$ is any complete  metric on $M$.)

\end{enumerate}
\end{theorem}

\noindent
{\bf Remark:}  (1). For each compact set $K$, $\xi^K>0$ almost surely. This is 
easily seen from $\xi(x)>0$ a.s. for each $x$ and the fact that a lower 
semicontinuous function on a compact set is bounded from below and assumes its 
minimum.

\bigskip

\noindent
(2).  \label{re: partial flow 2}As pointed out in \cite{EL82}, if there are two 
partial flows $(F_t^1, \xi_1)$, and $(F_t^2,\xi_2)$ which satisfy conditions 1-3
 of theorem $\ref{th: partial flow}$ for the $C^0$ topology, then  we have 
 uniqueness: for all $x$, $\xi_1(x)=\xi_2(x)$ almost surely.  Consequently 
$\inf_{x\in M} \xi_1=\inf_{x\in M}\xi_2$ almostly surely, and for each compact set 
$K$, $\xi^K_1=\xi_2^K$ almost surely. In particular if we have a version of the 
solution which is jointly continuous, then we actually have a version  smooth in
 the $C^{r-1}$ topology by part 3 of the theorem above.

\bigskip

From the theorem we see that starting from a compact set, the solution can be chosen to be continuous until part of it explodes. This and the following example suggests 
that strong completeness is a very demanding property, and there is a rich
 layer between being complete and being strongly complete.

\bigskip
\noindent
{\bf Example:} \cite{EL78}, \cite{ELbook} \label{ex: -0 1}
   Let $X(x)(e)=e$, and $A=0$. Consider the following 
stochastic differential equation  $dx_t=dB_t$ on $R^n$-$\{0\}$ for $n>1$. The solution is: $F_t(x)=x+B_t$, which is complete since for a fixed  starting point $x$, $F_t(x)$ almost  surely never hits $0$ . But it is not  strongly complete.  However for any $n$-$2$ dimensional hyperplane (or a submanifold) $H$ in the manifold, $\inf_{x\in H}\xi(x,\omega)=\infty$ a.s., since a Brownian motion does not charge a set of codimension $2$. To get an example on a complete metric space, apply the inversion map $z\mapsto {1\over z}$ in complex form  as in \cite{CA-EL83}. The resulting system on $R^2$ is $(\hat X, B)$ where  
\[      \hat X(x,y)=\left[ \begin{array}{cc} y^2-x^2 &2xy\\
-2xy &y^2-x^2\end{array} \right].        \]

\noindent
The corresponding solution is in fact (in complex notation): ${z\over 1+zB_t}$. We'll
continue this example on page ~\pageref{ex: -0 2}.

\bigskip
This leads to the following definition, which originally uses $p$-dimensional sub-manifolds, 
 those $p$-dimensional sub-manifolds are replaced by $p$-simplices at the suggestion of  D. Elworthy.

\begin{definition}  A\index{strongly $p$-complete} stochastic dynamical system on a  manifold is called strongly 
$p$-complete if  $\xi^K=\infty$  a.s. for every  $K\in S_p$. 
$\label{def:complete}$
\end{definition}
Here $S_p$ is the space of images of all smooth (smooth in the sense of extending
 over an open neighborhood) singular p-simplices. Recall a 
singular p-simplex in $M$ is a map from the standard p-simplex to $M$. For
 convenience we also use the term singular p-simplex for the image of a singular
 p-simplex map.

\bigskip

The example above on $R^n$-$\{0\}$ (for $n>2$) gives us a s.d.s. which is strongly
 n-2 complete, but not strongly (n-1)-complete. It is strongly  (n-2)-complete
since  a singular (n-2)-simplex
has finite Hausdorff (n-2) measure and is thus not charged by Brownian motion. It is
not strongly (n-1) complete from proposition ~\ref{pr: n-1 complete} on page
~\pageref{pr: n-1 complete}.

\bigskip

Of all these "completeness" notions, we are particularly interested in strong 
1-completeness, which helps us to get a result  on $d(P_tf)=(\delta P_t)f$
 (see page \pageref{th: differentiate semigroups}) and is used to get a homotopy 
vanishing result in theorem $~\ref{th: homotopy}$  replacing the obvious 
requirement of strong completeness. It turns out on most occasions, 
we only need   strong 1-completeness and  this follows from natural 
assumptions.

\bigskip
Once we get strong completeness, naturally we would like to know  when does the flow {\it consist of diffeomorphisms}, i.e.  there is a version of the flow\index{flow: diffeomorphism} such that except for a set of probability zero, $F_t(\cdot, \omega)$ is a diffeomorphism from $M$ to $M$ for each $t$ and $\omega$. This is basically the "onto" property of the flow, since the flow is always injective as showed in \cite{Kunitabook} and by part 2 of theorem 
~\ref{th: partial flow}. We will discuss this  at the end of the next section.

\noindent{Note:} Results in this chapter remain true when  (\ref{eq: basic}) 
is changed to a time dependent  equations.

\section{Main Results}

If not specified, by $(F_t,\xi)$ we mean the partial flow defined in 
theorem $~\ref{th: partial flow}$.

\begin{proposition} If the stochastic differential equation considered is

\noindent
 strongly p-complete, then  $\xi^N=\infty$ a.s.  for any $p$ dimensional submanifold 
$N$ of $M$.
\label{pr: on definition of strong}
\end{proposition}

\noindent
{\bf Proof:}   Let $N$ be a $p$ dimensional submanifold.   Since all smooth 
differential manifolds have a  smooth triangulation \cite{Munkres},
 we can write: $N=\cup V_i$. Here $V_i$ are smooth singular p-simplexes.
But  $\xi^{V_i}=\infty$ a.s. for each $i$ from the assumption. Thus 
 $F_\cdot(\cdot)|_{V_i}$ is continuous a.s. and thus so is $F|_N$ itself. Thus
 $\xi^N=\infty$  a.s. from remark 2 on page ~\pageref{re: partial flow 2}.
\hfill \rule{3mm}{3mm}

\bigskip
Note if $p$ equals the dimension of $M$, p-completeness gives back the usual
 definition of strong completeness, i.e. the partial flow defined in theorem 
~\ref{th: partial flow} satisfies $\inf_{x\in M}\xi(x)=\infty$ almost surely
 as  showed above. See also remark 2  after theorem $~\ref{th: partial flow}$. 
In this case we will continue to use strong completeness for strong n-completeness.

\bigskip

We need the following cocycle property from \cite{ELbook} for
the next proposition: For almost all $\omega$, for all $s>0$, $t>0$,

\begin{equation}
F_{t+s}(x,\omega)=F_t(F_s(x,\omega), \theta_s(\omega)).
\label{eq: cocycle}
\end{equation}

\noindent
Note that the exceptional set  for (\ref{eq: cocycle}) can be taken independent of
 $s$ and $t$, according to a recent survey by L. Arnold (manuscript).
However we do not need this refinement here.

\bigskip

\begin{proposition}
Let $\xi=\inf_{x\in M} \xi(x)$. If there is a number $\delta>0$ such that: 
$P\{\xi\ge \delta\}=1$,  then $\xi=\infty$ almost surely.
\end{proposition}

\noindent{\bf Proof:} 
There is no explosion by corollary ~\ref{co: wuc 2}. Let 
$\Omega_0$ be the set of $\omega$ such that if $\omega\not \in \Omega_0$, then 
$(t,x) \mapsto F_t(x,\omega)$ is continuous on $[0, \delta]\times M$. Consider
 $F_\cdot(\cdot, \theta_\delta(\omega))$ which is continuous on 
$[0, \delta]\times M$ if $\theta_\delta(\omega)\not \in \Omega_0$. Let
$\Omega_1=\Omega_0\cup\{\omega: \theta_\delta(\omega)\in \Omega_0\}$. Thus
 $\xi\ge 2\delta$ for $\omega \not \in \Omega_1$ and $\Omega_1$ has measure zero. 
Inductively we get $\xi=\infty$ almostly surely. \hfill \rule{3mm}{3mm}

\bigskip

The following proof was suggested to me by D. Elworthy, improving an earlier
result of myself proved in more restrictive situation:

\begin{proposition}\label{pr: n-1 complete}
A stochastic dynamical system on a n-dimensional 

\noindent
 manifold is strongly complete if  strongly (n-1)-complete. 
\end{proposition}

\noindent
{\bf Proof:}
Since strong n-completeness holds for compact manifolds, we shall assume $M$ is
not compact. Let $B$ be a geodesic ball centered at some point $p$ in $M$ with
 radius smaller than the injectivity radius at $p$. Since $M$ can be covered by
 a countable number of such balls, we only need to prove $\xi^B=\infty$ almost 
surely.

Let $B$ be such a ball. It clearly divides $M$ into two parts, one bounded and
 the other unbounded.  Write $M-\partial B=K_0\cup N_0$. Here $K_0$ is the
 bounded piece.  Fix $T>0$. By the ambient isotopy theorem there is a 
diffeomorphism from $[0,T]\times M$ to $[0,T]\times M$ given by:
$(t,x)\mapsto (t, h_t(x))$ for $h_t$ some diffeomorphism from $M$ to its image,
 and satisfying:
$$h_t|_{\partial B}=F_t|_{\partial B}.$$
Set $K_t=h_t(K_0)$, $N_t=h_t(N_0)$. Then
$$M=K_t\cup F_t(\partial B) \cup h_t(N_0),$$
 and
\begin{equation}
F_t(\stackrel {\circ}{B})\subset K_t
\label{eq: n-complete 1}
\end{equation}
 on $\{\omega: t<\xi^{B}(\omega)\}$.

Now
$$\cup_{0\le t\le T}\bar K_t={Proj}^1\left[H(\bar K_0 \times [0,T])\right],$$
here $\hbox{Proj}^1$ denotes the projection to $M$. Thus 
$\cup_{0\le t\le T} \bar K_t$ is compact. 
By $(~\ref{eq: n-complete 1})$, 
$F_t(B)=F_t(K_0)\cup F_t(\partial B)$, for $0\le t\le T\wedge \xi^B$, stays in a 
compact region. So $\xi^B\ge T$ almost surely from part 3 of theorem
 ~\ref{th: partial flow}.  \hfill \rule{3mm}{3mm}

\bigskip

 Take a sequence of nested relatively 
compact open sets  $\{U_i\}$ such that it is a cover for  $M$ and 
$\bar U_i\subset U_{i+1}$.  Let $\lambda^i$ be a standard  smooth cut off 
function such that:

\[\lambda^i =\left\{ \begin{array} {cl}1  &x\in U_{i+1}  \\
0, & x\not \in U_{i+2}. \end{array} \right.\]
Let $X^i=\lambda^i X$, $A^i=\lambda^i A$, and $F^i_\cdot$ the solution flow to the 
s.d.s. $(X^i,A^i)$. Then $F^i$ can be taken smooth since both $X^i$ and $A^i$ have compact support.  Let $S_i(x)$ \label{stopping times}
be the first exit time of $F^i_t(x)$ from $\bar U_i$ and $S_i^K=\inf_{x\in K} S_i(x)$ for a compact set $K$.  Thus 
$S^K_i$ is a stopping time.  Furthermore $F_t^i(x)=F_t(x)$ before  $S_i^K$.

\noindent
 Clearly  $S_i^K\le \xi^K$,  and in fact $\lim_{i\to\infty}S_i^K=\xi^K$ as proved in \cite{CA-EL83}.

\bigskip

\noindent  Let 
$$K_1^1=\{\hbox{Image}(\sigma)| \sigma: [0,\ell]\to M \hskip 3pt \hbox{is $C^1$, $\ell<\infty$}\}.$$

\begin{theorem}
 Let\index{strong 1-completeness: theorem} $M$ be a complete  connected 
Riemannian manifold. 
Suppose all the coefficients of the stochastic differential equation are  $C^{2}$, and assume there is a point $\bar x\in M$ with $\xi(\bar x)=\infty$ almost surely. Then we have $\xi^H=\infty$ for all $H\in K_1^1$,  if 

\begin{equation}
\liminf_{j\to \infty}\sup_{x\in K}  E\left(|T_xF_{S^K_j}|\chi_{S^K_j< t}\right)
<\infty
\label{eq: dim1 1}
\end{equation}

\noindent
for every compact set  $K\in K_1^1$ and each $t>0$. In particular when 
$(~\ref{eq: dim1 1})$ holds  we have strong 1-completeness, and strong completeness
 if dimension of  $M$ is less or equal to  $2$.
 \label{th: dim1}\end{theorem}

\bigskip
\noindent {\bf proof:}
Let $y_0\in M$.  We shall show $\xi(y_0)=\infty$.  Take a piecewise $C^1$ curve $\sigma_0$ connecting the two points $\bar x$ and $y_0$.  Suppose it is parametrized by arc length $\sigma_0 \colon [0, \ell_0]\to M$ with $\sigma_0(0)=\bar x$.

Denote by  $K_0$  the image set of the curve.
 Let $K_t=\{F_t(x): x\in K_0\}$, and $\sigma_t=F_t\circ \sigma_0$.  Then
on $\{\omega: t<\xi^{K_0}(\omega)\}$, $\sigma_t(\omega)$ is a piecewise $C^1$ curve.
Denote by $\ell(\sigma_t)$ the length of $\sigma_t$. 

Let $T$ be a stopping time such that  $T<\xi^{K_0}$, then: 

\begin{eqnarray*}
\ell(\sigma_{T(\omega)}(\omega))
&\le& \int_0^{\ell_0} |{d\over ds} \left(F_{T(\omega)}\left (\sigma(s),\omega\right)\right)|\, ds  \\
&\le& \int_0^{\ell_0} |T_{\sigma(s)}F_T\left(\omega\right)|\,ds.
\end{eqnarray*}

\noindent Thus for each $t>0$:

\begin{eqnarray}
E\ell(\sigma_T)\chi_{T<t}
&\le& \int_0^{\ell_0} E(\chi_{T<t}|T_{\sigma(s)}F_T|)\,ds\\
&\le& \ell_0 \sup_{x\in K_0} E\left(|T_xF_T|\chi_{T<t}\right). 
\label{eq: dim1 2}
\end{eqnarray}

\noindent
 Assume $\xi^{K_0}<\infty$.  Take $T_0$ with $P\{\xi^{K_0}<T_0\}>0$. 
Now $\cup_{0\le t\le T_0} F_t(\bar x,\omega)$ is a
 bounded set a.s. since $F_t(\bar x)$ is sample continuous in $t$ and 
$\xi(\bar x)=\infty$. 
Thus there is  $R(\omega)<\infty$ a.s. such that:

\begin{equation}
\sup_{0\le t\le T_0} d\left(F_t(\bar x,\omega), \bar x\right)\le R(\omega)<\infty.
\label{flow1} \end{equation}

\noindent  
But by theorem$~\ref{th: partial flow}$, almostly surely on  $\{\xi^{K_0}<\infty\}$

\begin{equation}
\lim_{t\nearrow \xi^{K_0}} \sup_{x\in K_0}  d\left(\bar x, F_t(x,\omega)\right)
=\infty .
\label{flow2}  \end{equation} 

\noindent
So by the triangle inequality we get:

$$\sup_{x\in K_0} d(F_t(x,\omega), F_t(\bar x,\omega))
\ge \sup_{x\in K_0} d\left(F_t(x, \omega), \bar x\right)-d\left(\bar x, F_t(\bar x, \omega)\right).$$

\noindent
Combining $(\ref{flow1})$ with $(\ref{flow2})$, we get, almostly surely 
 on $\{\omega: \xi^{K_0}<T_0\}$:

\begin{eqnarray*}
& &\liminf_{t\nearrow \xi^{K_0}} \sup_{x\in K_0} d\left(F_t(x,\omega),F_t(\bar x,\omega)\right)\\
&\ge& \liminf_{t\nearrow \xi^{K_0}} \sup_{x\in K_0} d\left(F_t(x,\omega),\bar x\right)
-\sup_{0\le t\le T_0} d\left(\bar x,F_t(\bar x,\omega)\right)\\
&=&\infty.
\end{eqnarray*}

\noindent
Therefore 
$\liminf_{t\nearrow\xi^{K_0}}\ell\left(\sigma_t(\omega)\right)=\infty$ almostly surely on $\{\xi^{K_0}<T_0\}$.

\noindent
Let $T_j=: S_j^{K_0}$ be the stopping times defined immediately before this theorem, which  converge to $\xi^{K_0}$. Then  there is a subsequence of 
$\{T_j\}$,  still denote by $\{T_j\}$, such that:

$$\lim_{j\to \infty} \ell(\sigma_{T_j})\chi_{\xi^{K_0}<T_0}
=\infty, $$

\noindent
Thus

$$E\lim_{j\to \infty} \ell(\sigma_{T_j})\chi_{\xi^{K_0}<T_0} =\infty.$$

\noindent
However from equation (~\ref{eq: dim1 2}) and our hypothesis $(~\ref{eq: dim1 1})$
, we have:

$$\liminf_{j\to \infty} E\ell(\sigma_{T_j(\omega)}(\omega))\chi_{\xi^{K_0}<T_0}
  \le  \ell_0 \liminf_{j\to\infty} \sup_{x\in K_0} E|T_xF_{T_j}|\chi_{T_j<T_0}
<\infty$$

\noindent
since $T_j<\xi^{K_0}$ almost surely.  Applying  Fatou's lemma, we have:

$$E\liminf_{j\to \infty} \ell(\sigma_{T_j})\chi_{\xi^{K_0}<T_0}
\le  \liminf_{j\nearrow\infty}E \ell(\sigma_{T_j})\chi_{\xi^{K_0}<T_0}<\infty.$$

\noindent 
This gives a contradiction. Thus $\xi^{K_0}=\infty$. In particular $\xi(y)=\infty$
for all $y\in M$.
 
Next take $K\in K_1^1$. Replacing $K_0$ by $K$ in the proof above we get 
$\xi^K=\infty.$  This is because in the proof above we only used the fact that 
there is a point     $\bar x$  in $K_0$ with $\xi(\bar x)=\infty$  and $|TF_t|$ 
satisfies $(\ref{eq: dim1 1})$.  

To see strong 1-completeness, just notice the set of smooth singular  1-simplexes 
$S_1$ is contained in $K_1^1$.
The proof is finished. \hfill \rule{3mm}{3mm}

\bigskip

\noindent {\bf Remark:} 
\bigskip

(1). We only need  inequality $~\ref{eq: dim1 1}$ to hold for 
one sequence of  exhausting open sets  $\{U_j\}$. In particular they may be different for different compact sets $K$. 

\bigskip

(2). The second inequality holds if we assume:
 
$$\sup_{x\in K} E\left(\sup_{s\le t}|T_xF_s|\chi_{s<\xi(x)}\right)<\infty$$

\noindent  for each number $t$ and compact set $K$. However we keep the inequality
with stopping times since it is sometimes easier to calculate from the original
stochastic differential equation, and gives sharper result as will be shown later in
the examples. See lemma ~\ref{le: applications 2} on page 
~\pageref{le: applications 2}.

(3).  The requirement on the connectness can be removed by assuming that
the s.d.e. is complete at  one point of each component of the manifold.

\bigskip

\noindent
{\bf Example: } \label{ex: -0 2}
(1).  The requirement for the manifold to   be  complete is necessary.
e.g. the example on $R^2-\{0\}$ on page ~\pageref{ex: -0 1}  satisfies
equation (\ref{eq: dim1 1}) but is  not strongly complete (see also the
 example on page ~\pageref{ex: Langevin}). In fact the transformed flow 
  $F_t(z)={z\over 1+zB_t}$  on $R^2$ by inverting does not satisfy the 
condition of the theorem on its derivative and it is not strongly 1-complete.

\noindent 
   (2).    Theorem ~\ref{th: dim1} does not work with
 equation (\ref{eq: dim1 1})  replaced  by 

\noindent
$\sup_xE|T_xF_t|<\infty$. This can be seen by using the above example
on $M=R^2-\{0\}$ but with the following  Riemannian metric:
$$|v|^{\tilde { } }={|v|\over |x|}, \hskip 6pt v\in T_xM,$$
 since the  metric  is complete and for each compact set $K$,
$$\sup_{x\in K}E|T_xF_t|^{\tilde { } } =\sup_{x\in K} E{1\over |x+B_t|}<\infty.$$

\bigskip

We say a s.d.e. is {\it complete at one point}  if there is a point $x_0$ in $M$ 
with $\xi(x_0)=\infty$. From the theorem we have the following corollary, which 
is known for elliptic diffusions.

\begin{corollary}
A stochastic differential equation with $C^2$ coefficients and satisfying hypothesis
$(\ref{eq: dim1 1})$ in theorem $~\ref{th: dim1}$
 is complete if it is complete at one point.
\end{corollary}

\bigskip

With strong 1-completeness we may apply the following  Kolmogorov's criterion on regularity \cite{ELbook} to get strong p-completeness:

\bigskip

\noindent
 {\bf Kolmogorov's criterion:}
 Let $M$ be a complete Riemannian manifold with $d( , )$ denoting the distance between two points.  Let $F$ be  a set of $M$-valued random variables indexed by  $[0,1]^p$ for which there exist positive numbers $\gamma$, $c$, and $\epsilon$ such that:

$$E d(F(\underline{s_1}), F(\underline{s_2} ))^\gamma 
\le c |\underline{s_1}-\underline{s_2}|^{p+\epsilon}$$

\noindent
hods for all $\underline{s_1}, \underline{s_2}$ in $[0,1]^p$. 
Then there is a  modification $\tilde F$ of $F$ such that the paths of
 $\tilde F$  are continuous. Here $|\underline{s_1}-\underline{s_2}|$ is the 
distance between  these two points (induced from $R^p$).

There is a corresponding result for  pathwise continuous stochastic processes 
$\{F_t(\cdot), t\ge 0\}$ parametrized by $[0,1]^p$.  
There is a version  which is jointly continuous in $t$ and $x$ on
 $[0,T]\times [0,1]^p$ for each $T$, if the following is satisfied: 

 $$E\sup_{s\le t}d( F_s(\underline{s_1}),  F_s(\underline{s_2}))^\gamma 
\le c |\underline{s_1} -\underline{s_2}|^{p+\epsilon}$$

\noindent
for all $\underline{s_1}, \underline{s_2}$ in $[0,1]^p$.  
This comes naturally by letting $N=C([0,T], M)$ 
with the following metric:
$\bar d(f,g)=\sup_{0\le s\le T}d(f(s), g(s))$.

\bigskip
\noindent
Recall a map $\alpha: [0,1]^p\to M$ is called Lipschitz continuous if there is a constant $c$ such that:

\begin{equation}
d\left(\alpha(\underline s), \alpha(\underline t) \right)
\le c|\underline s -\underline t|
\end{equation}

\noindent
 for all $\underline s$, $\underline t$ in $[0,1]^p$. 

\noindent
Denote by $L_p$ the space of all the image sets of such a Lipschitz map. This space contains $K_p$.
 
\bigskip

\begin{theorem}\label{th: strong completeness}
Let\index{strong p-complete: theorem} $M$ be a complete  connected 
Riemannian manifold. 
Consider a s.d.e. which is complete at one point and with  $C^2$ coefficients. 
Let $1\le d\le n$. Then we have $\xi^K=\infty$  for each $K\in L_d$
if  for each positive number $t$ and compact set $K$ there is a number 
$\delta>0$ such that:

$$\sup_{x\in K} E\left(\sup_{s\le t}|T_xF_s|^{d+\delta}\chi_{s<\xi}\right)<\infty.$$ 

In particular this implies strong d-completeness.
\end{theorem}

\noindent {\bf Proof:}  Let $\sigma$ be a Lipschitz map from $[0,1]^d$ to $M$
 with image set $K$. Take a compact set  $\hat K $ with the following property:
 for  any two points of $K$, there is a piecewise $C^1$ curve lying in $\hat K $
 connecting them.

For example the set $\tilde K$ can be taken in the following way:
Let  $\sigma$ be a minimum length curve in $M$ connecting two points of $K$;
 its length will be smaller or equal to  $\hbox{dia}(K)$.  Let  $\hat K $  be the 
closure of  the union of the image sets of  such curves.

Let $x=\sigma(\underline s)$ and  $y=\sigma(\underline t)$ be two points from $K$.  
Let $\alpha$ be a piecewise $C^1$ curve in $\hat K $ connecting them. 
Denote by $H_\alpha$ the image set of $\alpha$ and $\ell$ its length. 
By proposition $~\ref{th: dim1}$, $\xi^{H_\alpha}=\infty$. Thus for any $T_0>0$
 we have:

\begin{eqnarray*}
 E\sup_{t\le T_0}\left[d\left(F_{t}(x), F_{t}(y)\right)\right]^{d+\delta}
&\le&  E\left(\int_0^\ell \sup_{t\le T_0}|T_{\alpha(s)}F_{t}|\, ds\right)^{d+\delta}\\
&\le& \ell^{d+\delta-1}
 E\int_0^\ell \left(\sup_{t\le T_0}|T_{\alpha(s)}F_{t}|^{d+\delta}\right)\, ds\\
&\le&  \ell^{d+\delta} 
\sup_{x\in \hat K }\,\left(E\sup_{t \le T_0} |T_xF_{t}|^{d+\delta}\right).
 \end{eqnarray*}

\noindent
Taking infimum over a sequence of such curves which minimizing the distance 
between $x$ and $y$ we get:

$$ E\left(\sup_{t\le T_0}d(F_{t}(x), F_{t}(y))^{d+\delta}\right)
\le d(x,y)^{d+\delta}\sup_{x\in \hat K }\, E\left(\sup_{t\le T_0} |T_xF_t|^{d+\delta}\right).$$

The Lipschitz property of the map $\sigma$ gives

$$ E\left(\sup_{t\le T_0}d(F_{t}(\sigma(\underline s), F_{t}(\sigma(\underline t)))^{d+\delta}\right)   \le c|\underline s-\underline t|^{d+\delta}
\sup_{x\in \hat K }\, E\left(\sup_{t\le T_0} |T_xF_t|^{d+\delta}\right).$$

Thus we have a modification of   $\tilde F_\cdot(\sigma(-))$ of $F_\cdot(\sigma(-))$ which is jointly  continuous from $[0, T_0]\times [0,1]^d\to M$, according to the Kolmogorov's criterion. So for a fixed point $x_0$ in $M$:

$$\sup_{t\in [0, T_0]}\sup_{\underline s\in [0,1]^d} 
d(F_t(\sigma(\underline s), \omega), x_0)<\infty.$$

On the other hand on $\{\xi^K<\infty\}$, 

$$\lim_{t\nearrow \xi^K} \sup_{x\in K} d(F_t(x,\omega), x_0)=\infty$$

\noindent
almost surely. This give a contradiction. So $\xi^K=\infty$ for all $K\in L_d$. 

Notice every singular d-simplex has a representation of a Lipschitz map
from the cube $[0,1]^d$ to $M$ (by squashing  one half of the cube to the diagonal).
This gives the required strong p-completeness.  \hfill \rule{3mm}{3mm}

\bigskip
This theorem is used in section 7.3 to get a cohomology vanishing result.

\bigskip

\noindent{\bf Flows of diffeomorphisms}

To look at the diffeomorphism \label{diffeomorphism} property, we first quote the following theorem (first proved by Kunita) from \cite{CA-EL83}:

\bigskip\noindent
Let $M=R^n$. Assume both $X$ and $A$ are  $C^2$. If the s.d.s. is strongly complete, then  it has a flow  which  is surjective for each $T>0$  with probability one if and only if the adjoint system:
\begin{equation}  dy_t=X(y_t)\circ dB_t -A(y_t)dt
\label{eq: adjoint}
 \end{equation}
is strongly complete.
\bigskip

For a manifold this works equally well since equation (\ref{eq: adjoint})
 does give the inverse map up to distribution. Thus if both equations 
(\ref{eq: basic}) and 
(\ref{eq: adjoint})  satisfy the conditions of the theorem for strong completeness, 
the solution flow $F_t$ consists of diffeomorphisms. When there is a uniform cover 
for $(X,A)$, there is no explosion for  both $(\ref{eq: basic})$ and 
$(\ref{eq: adjoint})$ . In this case, the  solution consists of diffeomorphisms\index{flow of diffeomorphisms: theorem}  if  for $K$ compact:

$$\sup_{x\in K}E\sup_{s\le t} \left(|T_xF_s|^{n-1+\delta} 
+\left(|T_{F_s^{-1}(x)}F_s|^{-1}\right)^{n-1+\delta} \right)<\infty,$$
since both (\ref{eq: basic}) and (\ref{eq: adjoint}) are strongly complete
by theorem \ref{th: strong completeness}.

\section{Applications} $\label{applications}$

In this section we look at the stochastic differential equation to see when the conditions for the theorems above  are fulfilled. 

	As before let $M$ be a  Riemannian manifold. Let $x_0\in M$, $v_0\in T_{x_0}M$.  Noticing $|v_t|$ is almost surely nonzero for all $t$, \cite{Kunitabook},
we have the following formula for the $p^{\rm th}$ power of the norm of  $v_t$ from \cite{ELflow} for all $p$:  

\begin{equation}
\begin{array}{ll}
|v_t|^p=&|v_0|^p +p\sum_{i=1}^m\int_0^t|v_s|^{p-2} <\nabla X^i(v_s), v_s> dB_s^i\\
&+ p\int_0^t|v_s|^{p-2}<\nabla A(v_s), v_s>ds \\
&+{p\over 2} \sum_1^m\int_0^t |v_s|^{p-2}<\nabla^2 X^i(X^i, v_s), v_s> ds\\
&+{p\over 2} \sum_1^m\int_0^t |v_s|^{p-2} <\nabla X^i(\nabla X^i(v_s)), v_s> ds\\
&+{p\over 2} \sum_1^m\int_0^t |v_s|^{p-2}<\nabla X^i(v_s), \nabla X^i(v_s)>ds  \\
&+ \half p(p-2)\sum_1^m\int_0^t |v_s|^{p-4} <\nabla X^i(v_s), v_s>^2 ds.
\end{array}
\end{equation}

\noindent
on $\{t<\xi\}$. Let $v\in T_xM$. Define $H_p(v,v)$ as follows:

\[ \begin{array}{ll}
H_p(v,v)=& 2<\nabla A(x)(v), v> +\sum_{i=1}^m <\nabla^2 X^i(X^i,v), v>\\
&+\sum_1^m <\nabla X^i(\nabla X^i(v)),v> +\sum_1^m <\nabla X^i(v),\nabla X^i(v)>\\
&+(p-2)\sum_1^m {1\over |v|^2} <\nabla X^i(v), v>^2.
\end{array}\]

Let $\tau$ be a stopping time, then

\begin{equation}\begin{array}{cl}
|v_{t\wedge \tau}|^p=&|v_0|^p +p\sum_{i=1}^m\int_0^{t\wedge \tau} |v_s|^{p-2} <\nabla X^i(v_s), v_s> dB_s^i\\
&+{p\over 2}\int_0^{t\wedge \tau}|v_s|^{p-2}H_p(v_s, v_s)ds.
\label{eq: application 1}
\end{array}
\end{equation}

\noindent
This gives:

\begin{eqnarray*}
|v_{t\wedge\tau}|^{2p}  &\le&  2|v_0|^{2p} +4p^2\left [
	\sum_1^m \int_0^{t\wedge \tau} |v_s|^{p-2}
	 <\nabla X^i(v_s), v_s> dB_s^i\right]^2\\
	&    &+ 4p^2 \left[\int_0^{t\wedge \tau}|v_s|^{p-2}H_p(v_s, v_s)ds\right]^2\\
&\le&  2|v_0|^{2p} +4p^2 2^{m-1}\sum_1^m \left [
	\int_0^{t\wedge \tau} |v_s|^{p-2} <\nabla X^i(v_s), v_s> dB_s^i\right]^2\\
	&    &+ 4p^2 \left[\int_0^{t\wedge \tau}|v_s|^{p-2}H_p(v_s, v_s)ds\right]^2.
\end{eqnarray*}

\noindent
Let $T$ be a positive  number, then:

\[ \begin{array}{ll}
E\sup_{t\le T}|v_{t\wedge \tau}|^{2p} \,& \le 2|v_0|^{2p}\\
 &+  2^{m+1}p^2 \sum_1^m E\sup_{t\le T}  \left [
	\int_0^{t\wedge \tau} |v_s|^{p-2} <\nabla X^i(v_s), v_s> dB_s^i\right]^2\\
&+4p^2 E\sup_{t\le T}\left[
\int_0^{t\wedge \tau}|v_s|^{p-2}H_p(v_s, v_s)ds\right]^2.
\end{array}\]

\noindent
Applying Burkholder-Davies-Gundy inequality and H\"older's inequality we get:

\begin{eqnarray*}
E\sup_{t\le T}|v_{t\wedge\tau}|^{2p}  
&\le& 2|v_0|^{2p}
	 +4Tp^2 E\int_0^{T} \chi_{s \le\tau}|v_s|^{2p-4}H_p^2(v_s, v_s)ds\\
	&  & +2^{m+1}mp^2c_0 E\int_0^{T} \chi_{s\le \tau}|v_s|^{2p-4} 
	<\nabla X^i(v_s), v_s>^2 ds.\\
\end{eqnarray*}

\noindent
Here $c_0$ is the constant in Burkholder's inequality. 

Let $U$ be a relatively compact open set. Denote by $\tau(x)$
the first exit time of $F_t(x)$ from $U$.  For simplicity we write $\tau$ instead of $\tau(x_0)$. 

Since $X$ and $A$ are $C^2$, there is a constant $c$ such that:
$|\nabla X^i(x_s)|^2<c$ and $|H_p(v,v)|<c|v|^2$ on the set 
$\{s\le\tau(x_0)\}$.
Let $k_p=4(cp)^2(T+2^{m-1}mc_0)$, we have:

\begin{eqnarray*}
E\sup_{t\le T}|v_{t\wedge\tau}|^{2p}    
&\le& 2|v_0|^{2p}  +2^{m+1}mp^2c_0 E\int_0^{T} c^2|v_s|^{2p} \chi_{s\le \tau} ds \\
	& &  +  4Tp^2 E\int_0^{T} c^2|v_s|^{2p}\chi_{s \le\tau}ds\\
&=& 2|v_0|^{2p}
	+k_p E\int_0^{T} |v_s|^{2p} \chi_{s\le \tau} ds\\
&\le&  2|v_0|^{2p}
	+ k_p E\int_0^{T}
	 E\left(\sup_{u\le s}|v_u|^{2p} \chi_{u\le \tau} \right)ds\\
 &\le& 2|v_0|^{2p}
	+ k_p E\int_0^{T}
	 E\left(\sup_{u\le s}|v_{u\wedge \tau}|^{2p} \chi_{u\le \tau} \right)ds\\
 &\le&  2|v_0|^{2p}
	+ k_p E\int_0^{T}
	 E\left(\sup_{u\le s}|v_{u\wedge \tau}|^{2p}  \right)ds.\\
\end{eqnarray*}

\noindent
By Gronwall's lemma:

$$E\sup_{t\le T} |v_{t\wedge \tau}|^{2p} \le 2|v_0|^{2p}
{\rm e}^{k_pT}.$$

\noindent
On the other hand, taking an orthornormal basis  $\{e_i\}_1^n$ of $T_{x_0}M$,
we have:
$$E\left(\sup_{t\le T} |T_{x_0}F_{t\wedge \tau}|^{2p}\right) \le
 c\sum_1^n E\left(\sup_{t\le T} |T_{x_0}F_{t\wedge \tau}(e_i)|^{2p}
\right)\le c\exp^{k_pT}.$$
Here $c$ denotes  some constant depending only on $p$ and $n$.  
Thus we arrived at the following useful lemma, which is fairly well known. 

\begin{lemma}\label{le: applications 1}
(1). For each relatively compact open set $U$, there is a constant $c$ depending only on the bounds of the coefficients of  the stochastic differential equation on $U$ such that for all $p$:

\begin{equation}
E\sup_{s\le t} |T_xF_{s\wedge \tau(x)}|^{2p} \le 2c_1{\rm e}^{cp^2t}.
\label{eq: applications 2}
\end{equation}
Here $c_1$ is a constant depends on $n=dim(M)$.

(2). Assume both $|\nabla X|$ and $|H_p|$ are bounded. By the latter we mean

\noindent
 $|H_p(v,v)|\le k|v|^2$ for some constant $k$ and all $v\in T_xM$.  Then

	\begin{equation}
	\sup_{x\in M}E\left(\sup_{s\le t}|T_xF_{s\wedge \xi}|\right)^{2p}
\le 2c_1{\rm e}^{cp^2t}
\label{eq: applications 3}
\end{equation}

\noindent
for all $p$. Here  $c$ is a constant and $c_1$ depends only on $n=dim(M)$.
 In particular this is the case if $X$, $\nabla X$, $\nabla^2 X$, and $\nabla A$ 
are all bounded.
\end{lemma}
 
For the proof of $(~\ref{eq: applications 3})$, let $\tau=\tau^{U_n}$ in the above 
calculation, for
$\{U_n\}$ a sequence of nested relatively compact open set exhausting $M$. Then take the limit.  \hfill \rule{3mm}{3mm}

As a corollary, we have the following result on strong completeness:

\begin{corollary}\label{co: complete 2}
Let $M$ be a complete connected  Riemannian manifold. 
 Assume $|\nabla X|$ is bounded and
there is a constant $k$ such that   $|H_1(v,v)|<k|v|^2$.   Then the s.d.s.  is 
strongly complete if complete for one point.  Note the last condition is
satisfied  if $|X|+|\nabla X|+|\nabla^2 X|+|\nabla A|$ is bounded.
\end{corollary}

Note we do not use any sort of nondegeneracy in the above corollary.

\bigskip
 Let $M=R^n$. Assume the s.d.e. is given in It\^ o form. We have as a corollary the following known result: 

 \begin{corollary} A  stochastic differential equation on $R^n$ (in It\^o form ) is
 strongly complete if all coefficients are $C^2$ and globally Lipschitz continuous.
\label{co: strong completeness Euclidean}
\end{corollary}
\noindent{\bf Proof:} 
First we have completeness as is well known. Write:

\begin{eqnarray*}
dx_t&=&X(x_t)dB_t + A(x_t)dt\\
&=& X(x_t)\circ dB_t +\bar A(x_t)dt.
\end{eqnarray*}

Here $\bar A=A-\half \sum_1^m \nabla X^i(X^i)$. So

$$\nabla \bar A=\nabla A-\half \sum_1^m \nabla^2 X^i(v,X^i)-\half\sum_1^m \nabla X^i(\nabla X^i(v)).$$

 Note also on $R^n$:
 $\nabla^2 X^i(X^i, v)=\nabla^2 X^i(v, X^i)$ on $R^n$. Substituting  these in
equation $~\ref{eq: application 1}$ the second derivatives of $X^i$ disappear. Thus
\[\begin{array}{ll}
H_p(v,v)=&2<\nabla A(v),v> + \sum_1^m <\nabla X^i(v),\nabla  X^i(v)>\\
 &+(p-2) \sum_1^m {1\over |v|^2} <\nabla X^i(v), v>^2.
\end{array}\]

So the boundedness of $\nabla X$ and $\nabla A$ give us  strong completeness.
 \hfill \rule{3mm}{3mm}

\bigskip

\begin{lemma} \label{le: applications 2}
Let $\{U_j\}$ be a sequence of relatively compact open sets exhausting $M$,
and $S_j^K$  the stopping times  defined   before theorem $~\ref{th: dim1}$ on page
 ~\pageref{stopping times}. Here $K$ is  a compact set.

Assume $H_p(v,v)\le k|v|^2$ for some constant $k$.  Then for all $j$
 \begin{equation}
\sup_M E(T_xF_{t\wedge S_j^K}|^p)\le {\rm e}^{{p\over 2} kt}.
\label{eq:  applications 4}
\end{equation}

\end{lemma}

\noindent {\bf Proof: }
First we have: $E\sup_{s\le t}|TF_{s\wedge S_j^K}|^{2p}<\infty$ for $x\in K$ by a 
 similar proof  as for   $(~\ref{eq: applications 2})$.                      
So from $(\ref{eq: application 1})$ on page ~\pageref{eq: application 1}, we
obtain:

$$E|v_{t\wedge S_j^K}|^p=|v_0|^p +{p\over 2} \, E\int_0^{t\wedge S_j^K}
|v_s|^{p-2} H_p(v_s,v_s)\, ds, $$
   since the martingale part disappears.
 Thus, just as before,  there is the following estimate:

$$E|T_xF_{t\wedge S_j^K}|^p\le {\rm e}^{{p\over 2} kt}$$

\noindent
from $H_p(v,v)\le k|v|^2$ and  Gronwall's lemma. 
 \hfill \rule{3mm}{3mm}

\bigskip

\noindent {\bf Remark:}
In fact  $(\ref{eq:  applications 4})$ holds if $S_j^K$ is replaced by $S_j(x)$. 
Here $S_j(x)$ is the first exit time of $F_t(x)$ from $U_j$.

\bigskip

From the proof, we have the following corollary of theorem ~\ref{th: dim1}:
\begin{corollary}\label{co: applications 3}
 A s.d.e. is strongly 1-complete if it is  complete at one point and satisfies:
$$H_1(v,v)\le k|v|^2.$$
Here $k$ is a constant. 
\end{corollary}

\subsection*{The case of Brownian system}

Next we consider special cases. First we assume the s.d.s. considered is a
 Brownian motion with drift $Z$. Then $Z=\half \sum_1^m\nabla X^i(X^i)+A$, and

$$\nabla Z(v)=\half \sum_1^m \nabla^2 X^i(v,X^i)+\half \sum_1^m\nabla X^i\left(\nabla X^i(v)\right)+\nabla A(v).$$ 

\noindent
On the other hand,

\begin{eqnarray*}
\left <\nabla^2 X^i(X^i,v),v\right>-\left <\nabla^2X^i(v,X^i),v\right>
  &=& \left<R(X^i,v)(X^i),v \right>\\
  &=&-{\rm Ric}(v,v).
\end{eqnarray*}

\noindent
Here $R$ is the curvature tensor and Ric is the Ricci curvature. 
Thus:
\begin{equation}\begin{array}{ll} \label{eq: H for BM}
H_p(v,v)=&2<\nabla Z(x)(v),v>-Ric_x(v,v)+\sum_1^m |\nabla X^i(v)|^2 \\
	&+(p-2)\sum_1^m {1\over |v|^2} <\nabla X^i(v),v>^2.
\end{array}\end{equation}

And therefore we have the following theorem:

\begin{theorem} \label{th: strong completeness Hessian}
Assume the s.d.s. is a Brownian motion with gradient drift $\nabla h$. Then if
$\half{\rm Ric}-\hbox{Hess}(h)$ is bounded from below with $|\nabla X|$ bounded the  Brownian motion is strongly complete.  Here  Hess$(h)=\nabla^2 h$.
\end{theorem}

\noindent
{\bf proof:} By a result in \cite{BA86}, we have completeness if $\half$ Ricc-Hess(h) is bounded from below. The strong completeness follows from theorem 
~\ref{th: strong completeness} and lemma ~\ref{le: applications 1}. 
\hfill \rule{3mm}{3mm}

In the above theorem, the s.d.s. may be a Brownian motion with a general drift if 
we know the system  does not explode a priori. The nonexplosion problem is discussed in chapter 4 and chapter 7.

\subsection*{The case of gradient Brownian system\index{gradient Brownian system}}

Next we consider gradient Brownian motion as in the introduction and follow the ideas of \cite{EL-homotopy}. Let $\nu_x$ be the space of normal vectors to $M$ at $x$. There is the second fundamental form:

$$\alpha_x: T_xM\times T_xM \to \nu_x$$

\noindent
and the shape operator:

$$A_x: T_xM\times \nu_x \to T_xM$$

\noindent
related by $\left<\alpha_x(v_1,v_2), w\right> =\left<A_x(v_1,w), v_2\right>$.
If $Z(x): R^m\to \nu_x$ is the orthogonal projection, then

$$\nabla X^i(v)=A_x \left(v,Z(x)e_i\right)$$

\noindent
as showed in \cite{ELbook} and \cite{ELflour}.

Let $f_1,\dots f_n$ be an o.n.b. for $T_xM$. Consider $\alpha_x(v,\cdot)$ as a linear map from $T_xM$ to $\nu_x$. Denote by $|\alpha_x(v,\cdot)|_{H,S}$ the corresponding Hilbert Schmidt norm, and $|\cdot|_{\nu_x}$ the norm of a vector in $\nu_x$. Accordingly we have:

\begin{eqnarray*}
\sum_1^m \left<\nabla X^i(v), \nabla X^i(v)\right>
&=&\sum_{i=1}^m\sum_{j=1}^n \left<A_x(v, Z(x)e_i), f_j\right>^2\\
&=&\sum_{i=1}^m\sum_{j=1}^n \left<\alpha_x(v,f_j), Z(x)e_i\right>^2\\
&=& \sum_{j=1}^m |\alpha_x(v,f_j)|_{\nu_x}^2\\
&=&|\alpha_x(v,\cdot)|_{H,S}^2.
\end{eqnarray*}

\noindent
There is also:
$$\sum_1^m \left <\nabla X^i(v), v\right>^2=|\alpha_x(v,v)|_{\nu_x}^2,$$

\noindent
giving
\begin{equation}\begin{array}{cl}
H_p(v,v)=& -{\rm Ric}(v,v)+2<\nabla Z(v),v>+|\alpha_x(v,\cdot)|_{H,S}^2\\
 &+ {(p-2) \over |v|^2}|\alpha_x(v,v)|^2_{\nu_x}.
\end{array}
\end{equation}

Thus the corollary:

\begin{corollary}
Assume the second fundamental form is bounded. Then the gradient Brownian
 motion with drift $\nabla h$ is strongly complete if
 $\half \hbox{Ri}c-\hbox{Hess} (h)$ is  bounded from  below.  It consists of
 diffeomorphisms if both $\half \hbox{Ric}-\hbox{Hess}(h)$ and
 $\half \hbox{Ricci} +\hbox{Hess}(h)$ are bounded from below.
\end{corollary}

\noindent{\bf Proof:} The strong completeness is clear from the previous
theorem. The diffeomorphism property comes from the fact that its adjoint 
equation is also a gradient Brownian system (with drift -$\nabla h$).

\bigskip

Further,  there is the following Gauss's theorem:

$${\rm Ric}(v,v)
=\left<\alpha(v,v), \hbox{trace } \alpha\right> -|\alpha(v,\cdot)|^2_{H,S}.$$

\noindent
Giving:
\begin{equation}\begin{array}{cl}
H_p(v,v)=& -<\alpha(v,v), \hbox{trace } \alpha> 
	+ 2|\alpha_x(v,\cdot)|_{H,S}^2\\
&+{1\over |v|^2} (p-2)|\alpha_x(v,v)|^2_{\nu_x}
	+2<\hbox{Hess}(h)(v), v>.
\end{array}
\end{equation}

\bigskip

Thus the completeness and strongly completeness of a gradient Brownian motion rely 
only on the bound on the second fundamental form and the bound on the drift:
there is no explosion if $h=0$ and $|\alpha_x|\le k( |x|)$ with $k$ a function on
 $R_+$  satisfying:
$$\int {1\over k(r)} dr =\infty$$
from Gauss theorem  and the example on page ~\pageref{ex: infinity Ricci} in
 chapter 4. 
Furthermore a gradient Brownian motion is strongly complete
 if  the second fundamental form is bounded.

\chapter{Derivative semigroups}

\section{Introduction}

 Assume the derivative of the  solution flow of equation $(\ref{eq: basic})$
 has first moment:  $E|T_xF_s\chi_{s<\xi(x)}|<\infty$,    we may define a
 semigroup (formally) of linear operators $\delta P_t$\index{$\delta P_t$} on
 1-forms as follows: for $v\in T_xM$ and $\phi$ a 1-form

\begin{equation}
\delta P_t\phi(v) =E\phi\left(T_xF_t(v)\right)\chi_{t<\xi(x)}
\label{eq: two1}
\end{equation}
It is in fact a semigroup on $L^\infty(\Omega)$, the space of bounded 1-forms, if 
 \begin{equation}
\sup_{x} E|T_xF_t\chi_{t<\xi(x)}|<\infty. \label{eq: two2}
\end{equation}

\vspace{5mm}\noindent 
 We are interested in three problems:
\begin{enumerate}
\item
When is $\delta P_t$ well defined, as a strongly continuous semigroup?
\item
When is $dP_tf=\delta P_t(df)$?
\item
When  is $\delta P_t \phi =\heatsemif \phi$, if $\A=\half \triangle^h$?
\end{enumerate}

If all answers to the questions are yes, we can obtain informations of heat 
semigroups and answer the question wether $\delta P_t$ sends closed forms to closed
 forms. These problems are also the basis for the next two chapters and 
will be discussed in detail and in great generality in this chapter. For related 
discussions, see Vauthier\cite{Vauthier79} and  Elworthy\cite{ELflour}.

However  first let $\A=\half\triangle^h$.  We have nonexplosion if 
$dP_tf=\delta P_t(df)$ for $f\in C_K^\infty$ and  if $E|T_xF_t|<\infty$ for all $x$, 
as will be shown 
in proposition $~\ref{pr: nonexplosion}$ in chapter 7.  So it is natural to assume
 completeness. For a complete stochastic dynamical system on a complete Riemannian
 manifold,  two basic assumptions:
\begin{equation}
E\sup_{s\le t}|T_xF_s|<\infty  
\end{equation}
and
\begin{equation}
\sup_{x\in M} E|T_xF_t|^{1+\delta}<\infty
\end{equation}
will give  everything we need:
\begin{enumerate}
\item
	$\delta P_t \phi =\heatsemif \phi$, for  $\phi\in L^\infty$,
\item
	$dP_tf=\delta P_t(df)$,  for $f\in C_K^\infty$
 \item   
and  strong 1-completeness.
\end{enumerate}
as shown later. c.f. theorem $~\ref{th: dim1}$ on page \pageref{th: dim1}, 
theorem $~\ref{th: differentiate semigroups}$ on 
page ~\pageref{th: differentiate semigroups} and proposition
 $~\ref{pr: nonexplosion}$ on page ~\pageref{pr: nonexplosion}.
 This basic assumption is satisfied by solutions of a s.d.e. with all the coefficients and their first two derivatives bounded as shown in section ~\ref{applications}.

\bigskip

But first we recall the  properties of probability semigroups for functions.

\section{Semigroups for functions}
\label{semigroup for functions}

Let $P_t$ be the probability semigroup  on bounded measurable   functions determined 
by our stochastic dynamical system.  Let $\A$ be its generator.  Then $\A={1\over 2}\sum X^iX^i +A$ on $C_k^\infty$, the space of smooth compactly supported functions. If further we assume completeness, then it sends bounded continuous functions to bounded continuous functions as showed in 
\cite{ELbook}.   
\bigskip

Let  $\A={1\over 2} \triangle^h$, then $P_t$ is a strongly continuous $L^2$ 
semigroup restricting to $L^2\cap L^\infty$(see \cite{EL-RO88}). Associated
 with $\half \triangle^h$, there is also the functional analytic semigroup
$\heatsemi$. These two semigroups agree as in \cite{EL-RO88}. See  also the 
proposition below. 
Thus from theorem ~\ref{th: STRana} on page  ~\pageref{th: STRana}, 
$P_tf$ is smooth on $L^2\cap L^\infty$. Moreover $P_t$ is $L^p$ contractive on 
 $L^2\cap L^p\cap  L^\infty$ for all $t$, $1\le p\le \infty$. 
See \cite{FU-ST} for more discussions on the $L^p$  contractivity of probability
 semigroups. 

\bigskip

Finally we have the following   known result :

\begin{proposition}\label{pr: semigroups for functions 2}
Let $M$ be a complete Riemannian manifold, then
 
$$P_t1(x)=\heatsemi 1(x)=\int p_t^h(x,y)e^{2h(y)} dy.$$
\end{proposition}

\noindent{\bf Proof:}  
First we show $P_tf=\heatsemi f$ for $f\in L^2\cap L^\infty$. 
Since $P_t$ is a strongly continuous $L^2$ semigroup on $L^2\cap L^\infty$, 
it extends  to a strongly continuous
$L^2$ semigroup $\tilde P_t$ on the whole $L^2$ space. Let $\tilde \A$ be the
 generator of
$\tilde P_t$. Then $\tilde \A$  is a closed operator by theorem
 ~\ref{th: properties of semigroups} and agrees with $\half \triangle^h$ on 
$C_K^\infty$. Thus 
$\tilde \A=\half \bar\triangle^h$  since there is only one closed extension for
$\half \triangle^h$ from the essential self-adjointness of $\triangle^h$ obtained
 in chapter 2. 

Applying the uniqueness theorem for the semigroup of class $C_0$(theorem 
~\ref{th: properties of semigroups}),
we get $\tilde P_tf=\heatsemi f$. Thus $P_tf=\heatsemif$ on $L^2\cap L^\infty$.

Next let $\{g_n\}$ be an increasing sequence of functions in 
$C_K^\infty$ approaching $1$ with $0\le g_n\le 1$. Such a sequence exists as shown 
in the appendix. Then $\heatsemi g_n \to \heatsemi 1$  since $p_t^h(x,-)$ is 
in $L^1$ (c.f. theorem $~\ref{th: STRana}$). But
 $\heatsemi g_n(x)=P_t g_n(x)\to P_t 1(x)$ for each $x$. So the limits must be
 equal: $P_t1=\heatsemi 1$. \hfill \rule{3mm}{3mm}

\section{ $\delta P_t$ and $dP_t$}

With  the help of results on strong 1-completeness, we have the following theorem which improves a theorem in \cite{ELbook}, where strong completeness and bounds on
curvatures are assumed.

\begin{theorem}  $\label{th: differentiate semigroups}$
Assume strong 1-completeness   and suppose for each compact set $K$, there is a constant $\delta>0$ such that:
$$\sup_{x\in K}E |T_xF_t|^{1+\delta} <\infty.$$ 
Then $P_tf$ is $C^1$ and   $$d(P_tf)=(\delta P_t)(df)$$
for  any $C^1$ function $f$ with both $f$ and $df$ bounded.
\end{theorem}

\noindent {\bf Proof:}   
Let $(x,v)\in TM$. Take a geodesic curve $\sigma\colon  [0,\ell]\to M$ starting from $x$ with velocity $v$. By the strong 1-completeness, $F_t(\sigma(s))$ is a.s. differentiable with respect to $s$. So for almostly all $\omega$: 

$${f\left(F_t\left(\sigma(s), \omega\right)\right)
-f\left(F_t(x, \omega)\right) \over s} 
={1\over s} \int_0^s df
\left(T_{\sigma(r)}F_t\left(\dot \sigma(r),\omega\right) \right)dr.$$

\noindent
Let 

$$I_s
={1\over s} \int_0^s df\left(T_{\sigma(r)}F_s(\dot \sigma(r), \omega)\right)dr.$$

\noindent
We want to show: $\lim_{s\to 0}E I_s=E\lim_{s\to 0} I_s$.
By the strong 1-completeness, we know  $T_{\sigma(r)}F_t(\dot \sigma(r))$ is
 continuous in $r$ for almost all $\omega$. Thus:
\begin{eqnarray*}
E\lim_{s\to 0} I_s &=&
E\lim_{s\to 0} {1\over s} \int_0^s df(T_{\sigma(r)}F_t(\dot \sigma(r))(\omega)dr\\
&=& Edf(TF_t(v)).
\end{eqnarray*}
On the other hand, $Edf(T_{\sigma(r)}F_t(\dot \sigma(r)))$  is continuous in $r$ if
$|T_{\sigma(r)}F_t|$ is uniformly integrable in $r$ with respect to the probability measure $P$. This is so if

$$\sup_rE|T_{\sigma(r)}F_t|^{1+\delta}<\infty.$$

\noindent   Thus 

$$\lim_{s\to 0} {1\over s} \int_0^s Edf
\left (T_{\sigma(r)}F_t(\dot\sigma(r),\omega)\right)dr
=Edf(TF_t(v))$$

\noindent
and the proof is finished.  \hfill \rule{3mm}{3mm}

\bigskip

\begin{corollary}
Let $M$ be a complete Riemannian manifold.  Suppose our s.d.s. is complete. If
$$E\sup_{s\le t}|T_xF_s|<\infty$$
and
$$\sup_{x\in K} E|T_xF_t|^{1+\delta}<\infty$$
for each $t>0$, and $K$ compact,  
then $dP_tf=(\delta P_t)(df)$   for all functions $f$ with both $f$ and $df$ bounded.
 \label{co: Derivative 1}
 \end{corollary}

In terms of the coefficients of s.d.e., we have:
\begin{corollary} \label{co: Derivative 30}
Let $M$ be a complete Riemannain manifold. Suppose our s.d.s. is complete and
satisfies:
$$H_{1+\delta}(v,v)\le k|v|^2.$$
Then $d P_tf=\delta P_t(df)$ if both $f$ and $df$ are bounded. Here
\[ \begin{array}{ll}
H_p(v,v)=& 2<\nabla A(x)(v), v> +\sum_{i=1}^m <\nabla^2 X^i(X^i,v), v>\\
&+\sum_1^m <\nabla X^i(\nabla X^i(v)),v> +\sum_1^m <\nabla X^i(v),\nabla X^i(v)>\\
&+(p-2)\sum_1^m {1\over |v|^2} <\nabla X^i(v), v>^2.
\end{array}\]
\end{corollary}
\noindent{\bf Proof:}
From corollary ~\ref{co: applications 3} on page ~\pageref{co: applications 3},
we have strong 1-completeness. Furthermore
$$\sup_{x\in M} E|T_xF_t|^{1+\delta \over 2}<\infty$$
if $H_{1+\delta}(v,v)\le k|v|^2$.  \hfill \rule{3mm}{3mm}

\bigskip

\noindent{\bf Remark:}
For a stochastic differential equation on $R^n$ (in It\^o form),
\[ \begin{array}{ll}
H_p(v,v)=& 2<\nabla A(x)(v), v>  +\sum_1^m <\nabla X^i(v),\nabla X^i(v)>\\
&+(p-2)\sum_1^m {1\over |v|^2} <\nabla X^i(v), v>^2,
\end{array}\]
as on page ~\pageref{co: strong completeness Euclidean}. Also notice if 

$$|X(x)|\le k(1+|x|)$$
and
$$<A(x),x>\le k(1+|x|)|x|,$$
then the system is complete. Thus our result on $dP_tf=\delta P_t(df)$ improves
a theorem in \cite{ELbook}, where global Lipschitz continuity is assumed
of the coefficients.
Here is a brief proof for the nonexplosion claim we made:

 Let $x_0\in M$, $x_s=F_s(x)$. On $\{\omega: t<\xi(x,\omega)\}$,
\begin{eqnarray*}
|x_t|^2&=&|x_0|^2 +2\int_0^t<x_s, X(x_s)dB_s> +2\int_0^t <x_s, A(x_s)>ds\\
& & +\int_0^t |X(x_s)|^2 ds.
\end{eqnarray*}
Let $T_n(x)$ be the first exit time of $F_t(x)$ from the ball $B_n(0)$. Write $T_n=T_n(x_0)$ for simplicity. Then
\begin{eqnarray*}
E|x_{t\wedge T_n}|^2&=&|x_0|^2  +2E\int_0^{t\wedge T_n} <x_s, A(x_s)>ds 
+ E\int_0^{t \wedge T_n} |X(x_s)|^2 ds\\
&\le& |x_0|^2 +6k\int_0^{t\wedge T_n} (1+|x_s|^2) ds.
\end{eqnarray*}
By Grownall's inequality:
$$E|x_{t\wedge T_n}|^2 \le (|x_0|^2+6kt)e^{6kt}.$$

Thus $P\{T_n<t\}=0$, so there is no explosion.

 \bigskip
\begin{lemma} \label{le: Differentiate 10}
 Let $M$ be a complete Riemannian manifold. Consider a s.d.s. with
 generator $\A=\half \triangle^h$. 
Let $f\in C_K^\infty$. Then the following are equivalent:
\begin{enumerate}
\item
 $\heatsemif (df)=\delta P_t(df).$ 
\item
$\delta P_t(df)=d(P_tf).$
\end{enumerate}
\end{lemma}
This is because  $P_tf=\heatsemi f$ from section 6.2 and 
$\heatsemif (df)=d(\heatsemi f)$ from proposition $~\ref{pr: exchangeability}$.

\begin{corollary}\label{co: Stflour 1} Let $M$ be a complete Riemannian manifold.
Assume the conditions of the above theorem and $\A={1\over 2} \triangle^h$. Then
for $f\in C_K^\infty$,
 $$\delta P_t(df)=\heatsemif (df).$$
\end{corollary}

\bigskip

\noindent {\bf Remark:} We will show later that  if $d(P_tf)=\delta P_t(df)$ 
for  $f$ with $f$ and $df$ bounded, then there is no  explosion.   See proposition 
$~\ref{pr: nonexplosion}$ and its corollary. 

\bigskip

Finally as known in the compact case \cite{ELflour}, we also have: 

\begin{proposition}   \label{pr: Derivative 2}
Assume $E\left(|T_xF_t|\chi_{t<\xi(x)}\right)$ is finite and continuous in 
$t$ for $t\in[0,a]$. Here $a$ is positive constant.  Let $\phi$ be a 1-form in $C_K^\infty$. Then

$${\partial (\delta P_t \phi)\over \partial t}|_{t=0} ={\cal L}\phi$$
with  $\lim_{t\to 0} \delta P_t\phi(v)=\phi(v)$ for each $v\in T_xM$. 
Here $\cal L$ is as defined on page $~\pageref{sy: ell}$. 

In particular if $\A= {1\over 2}\triangle +L_Z$ and $\phi$ is closed:

$${\partial (\delta P_t\phi) \over \partial t} |_{t=0}
=\left({1\over 2} \triangle^{1}+L_Z\right)\phi.$$

If the s.d.s. considered is a gradient system, we do not require $\phi$ 
 to be closed in the above.
\end{proposition}

\noindent{\bf Proof:} Take $v_0\in T_{x_0}M$. 
Applying   It\^o formula to $\phi$ we have:
\begin{eqnarray*}
\lim_{t\to 0} {\delta P_t\phi -\phi \over t}(v_0)
=\lim_{t\to 0} {1\over t} \int_0^t E\left( {\cal L}\phi_{x_s}(v_s)
\chi_{t<\xi}\right)\, ds
={\cal L}\phi_{x_0}(v_0)
\end{eqnarray*}
since ${\cal L}$ is a local operator(so ${\cal L}\phi$ remain bounded and 
 continuous) and  $E|T_{x_0}F_s(v_0)|$ is continuous in $s$.

\section{Analysis of $\delta P_t$ for Brownian systems}

This section will be devoted to discussions of $\delta P_t$ in the special
 case of  $\A={1\over 2} \triangle^h$.  The situation here  is particular nice,
since  the generator are both self-adjoint and elliptic. However 
to start with we would like to mention the following relevant theorem for a 
general elliptic generator, part 2 of which   improves a theorem  of Elworthy 
 \cite{ELflow}(part 1,3 and 4 here are  new):

As usual let $\{U_n\}$ be a sequence of nested relatively compact open sets in $M$ 
with  $\bar U_n\subset U_{n+1}$ and $\cup U_n =M$. Denote by $T_n$  the first exit 
time of $F_t(x)$  from $U_n$.

\begin{theorem} 
$\label{th: Stflour}$
  Let $\A= {1\over 2}\triangle +L_Z$.  Let $\{\phi_t: t>0\}$ be a regular solution of the heat equation:
$${\partial \over \partial t}\phi_t ={1\over 2} \triangle^{1}\phi_t+ L_Z\phi_t$$
with $d\phi_t=0$ for all $t$.  Then 
$\phi_t(v_0)=(\delta P_t)\phi_0(v_0)$ if one of the following conditions hold:

\begin{enumerate}
\item  Suppose   $|\phi_s|$ and $|\nabla \phi_s|$ are uniformly bounded in $s$ in $[0, t]$ for any $t$, and assume the s.d.s. is complete with  $|\nabla X|$ bounded and
$$\int_0^t E|T_xF_s|^2 ds<\infty, \hskip 12pt x\in M.$$
\item  
Assume completeness, $\phi_0$ bounded,
 and $E\left(\sup_{s\le t}|T_xF_s|\right)<\infty$ for each $t$ and $x$. 
 \item 
Suppose  $|\phi_s|$ is uniformly bounded in $s$ in $[0, t]$ for all $t$, and assume the s.d.s. is complete with  
$\sup_{n}E\left(|TF_{ T_n(x)}|^{1+\delta}\chi_{T_n<t}\right)<\infty$  
for some constant  $\delta>0$ and each $t$.

\item 
Assume for each $t$, $|\phi_s|$ converges uniformly in $s\in [0,t]$ to zero as $x\to \infty$  and 
$$\limsup_n E|TF_{ T_n(x)}|\chi_{T_n<t}<\infty$$

\end{enumerate}
\end{theorem}

\noindent
{\bf Proof:}  (a).	 Let $x_0\in M$, and $v_0\in T_{x_0}M$. 
First assume completeness. Apply It\^o formula to $\phi_{T-t}(x)$ for fixed
 number $T>0$, to get:
$$\phi_{T-t}(v_t)
= \phi_T(v_0) +\int_0^t \nabla \phi_{T-s}(X(x_s)dB_s)(v_s)
	+\int_0^t \phi_{T-s} (\nabla X(v_s) dB_s).$$

\noindent 	Let $t=T$, we get:
\begin{eqnarray*}
\phi_0(v_T)&=&\phi_T(v_0) 
	+\int_0^T \nabla \phi_{T-s}(X(x_s)dB_s)(v_s)  \\
  &  &+\int_0^T \phi_{T-s} (\nabla X(v_s) dB_s).
\end{eqnarray*}

\noindent   Take expectations on both sides to get:
$$\phi_T(v_0)=E\phi_0(v_T)$$ under either of the first two  conditions.

\noindent
(b).	For the rest we apply the It\^o formula with stopping time:
\begin{eqnarray*}
\phi_{T-t\wedge T_n}  (v_{t\wedge T_n})
&=& \phi_T(v_0) +\int_0^{t\wedge T_n} \nabla \phi_{T-s}(X(x_s)dB_s)(v_s) \\
 &  &+\int_0^{t\wedge T_n} \phi_{T-s} (\nabla X(v_s) dB_s).
\end{eqnarray*}

\noindent
Setting  $t=T$, we get:
\begin{eqnarray*}
&&\phi_0(v_T)\chi_{T<T_n} +\phi_{T-T_n}(v_{T_n})\chi_{T>T_n}\\
&  &= \phi_T(v_0)
 +\int_0^{T\wedge T_n} \nabla \phi_{T-s}(X(x_s)dB_s)(v_s)   \\
  &  &+\int_0^{T\wedge T_n} \phi_{T-s} (\nabla X(v_s) dB_s).
\end{eqnarray*}
Taking expectations of both sides above to get: 
$$E\left(\phi_0(v_T)\chi_{T<T_n}\right) 
+E\left(\phi_{T-T_n}(v_{T_n})\chi_{T>T_n}\right)
= \phi_T(v_0).$$

\noindent
(c).	Assuming part 3, $\phi_{T-T_n}(v_{T_n})\chi_{T>T_n}$ is uniformly integrable. So

$$\lim_{n\to \infty}E\left(\phi_{T-T_n}(v_{T_n})\chi_{T>T_n}\right) 
=E\left(\lim_{n\to \infty}\phi_{T-T_n}(v_{T_n})|\chi_{T>T_n}\right)=0.$$
since $T_n\to \infty $ from completeness while $$E\phi_0(v_T)\chi_{T<T_n}\to \delta P_t(\phi_0)(v_0).$$

\noindent
(d).	For part 4: Let $\epsilon>0$. There is a number $N>0$ s.t. if $n>N$,
$$|\phi_t(x)|<\epsilon  \hskip 6pt \hbox{for all $t$ and $x\not \in U_n$ }.$$
Thus 
 $$ E \left(|\phi_{T-T_n}(v_{T_n})|\chi_{T>T_n}\right)
\le \epsilon E\left(|v_{T_n}|\chi_{T>T_n}\right).$$ 
 Letting $\epsilon\to 0$, we get
 $\limsup_n E |\phi_{T-T_n}(v_{T_n})|\chi_{T>T_n}|=0$.
Thus the result.  \hfill \rule{3mm}{3mm}

\bigskip

\noindent{\bf Remarks:} 
(1).  For gradient system, we do not need the assumption $d\phi_s=0$ in the theorem. 

 (2).Let $M$ be a complete Riemannian manifold and  $\A=\half \triangle^h$.  
This theorem  gives us the equivalence of $\delta P_t$ and $\heatsemif \phi$ for 
closed forms, while corollary ~\ref{co: Stflour 1} only gives us the equivalence
 for exact forms.

 (3). Let $\phi\in C_K^\infty$.  Then $\heatsemif\phi$ is a regular solution  to the heat 
equation.  But we know $d\heatsemif \phi={\rm e}^{\half t\triangle^h}d\phi$ from 
section 2.3.  So in the above theorem we only need  $d\phi_0=0$ instead of
 $d\phi_t=0$ for all $t$.

\begin{corollary}\label{co: Derivative Ricci}
Let $M$ be a manifold with Ric-2Hess(h) bounded from below. Consider a Brownian 
system with gradient drift. If for some constant $\delta>0$:
$$\sup_nE|TF_{T_n}|^{1+\delta}\chi_{T_n<t}<\infty,$$
 then $\heatsemi \phi=E\phi(TF_t)$ for bounded closed 1-forms $\phi$. Note the 
required inequality holds if $H_{1+\delta}(v,v)\le c|v|^2$ for some constsnt $c$.
\end{corollary}

\noindent{\bf Proof:} Note there is no explosion by \cite{BA86}.
Furthermore if Ricci-2Hess(h) is bounded from below by  a constant $c$, then
$$|W_t^h|\le {\rm e}^{ct}$$
as in \cite{ELflour}. Here $W_t^h$ is the Hessian  flow satisfying:

$${DW_t^h\over \partial t}
=-\half{\rm Ric}(W_t^h,-)^\# + <{\rm Hess}(h)(W_t^h),->^\#.$$

\noindent
 Therefore  $\heatsemi \phi=E\phi(W_t^h)$ for bounded 1-form $\phi$. 
See \cite{ELflow} for a  proof. Hence $\heatsemi\phi$ is uniformly bounded in $t$
 in finite intervals. 
Apply part 3 of theorem ~\ref{th: Stflour} and lemma ~\ref{le: applications 2}, 
we get the conclusion. End of the proof.  \hfill \rule{3mm}{3mm}

\bigskip

Note also if  $h=0$(use $W_t$ for $W_t^h$),  we may get more information on
 the heat semigroup:

\begin{eqnarray*}
|\heatsemi \phi|&\le& E|\phi|^2 E|W_t|^2\\
&\le& {\rm e}^{2ct} P_t(|\phi|^2).
\end{eqnarray*}
But $P_t$ has the $C_0$ property if Ricci is bounded from below. Thus 
$\heatsemi \phi(x)$ converges to zero uniformly in $t$ in finite intervals as
$x$ goes to infinity. 

For more discussions on ralations between the Ricci flow and the derivative flow,
see \cite{EL-YOR}.

\bigskip

The following result  does not assume completeness. It is interesting
 since it does not assume any  condition on the heat semigroup either. In fact this
is used in the next chapter to get a nonexplosion result. See proposition
~\ref{pr: nonexplosion}.

\begin{proposition} \label{pr: Derivative 3}
 Let $p>1$, take $q$  to be the conjugate number to $p$: 
${1\over p}+{1\over q}=1$.  Then  $\delta P_t$ is a  $L^p$
 semigroup,  if  for each $t>0$:

\begin{equation} \sup_M E\left( |T_xF_t|^q\chi_{t<\xi}\right)<\infty.
\end{equation} 

 Also if there is a number $a>0$ such that:
\begin{equation} \sup_{s\le a}\sup_M E\left(|T_xF_s|^q \chi_{s<\xi}\right)
=c<\infty.
\end{equation}
and  the map $t\mapsto |T_xF_t|\chi_{t<\xi}$ is continuous into 
$L^p(\Omega,\F,P)$ at $t=0$ for each $x$,
 then $\delta P_t$ is a strongly continuous $L^p$ semigroup.

Furthermore if also
 $\lim_{t\to 0} \sup_{x\in M} E\left(|T_xF_t|\chi_{t<\xi}\right)=1$, 
then  $\delta P_t$ is a strongly continuous semigroup in $L^r$  for
 $p\le r\le \infty$.
\end{proposition}

\noindent
{\bf Proof:}  Take a 1-form $\phi \in L^p$. Then
\begin{eqnarray*}
\int |\delta P_t\phi|^{p} e^{2h}dx 
 &\le& \int E\left(|\phi|^p_{F_t(x)}\chi_{t<\xi}\right) 
\left(E|T_xF_t|^q\chi_{t<\xi}\right)^{p \over q} e^{2h}dx   \\
 &\le& \left(\sup_x E|T_xF_t|^q\chi_{t<\xi}\right)^{p-1} 
\int E\left(|\phi|^p_{F_t(x)}\chi_{t<\xi}\right) e^{2h} dx \\
&=& \left(\sup_x E|T_xF_t|^q\chi_{t<\xi}\right)^{p-1} \int|\phi|_x^{p} e^{2h}
 dx.
\end{eqnarray*}

\noindent
The last equality comes from the fact that $e^{2h}dx$ is the invariant measure
 for $F_t$. So we have showed $\delta P_t\phi$ is in $L^p$. 

To show the strong continuity of $\delta P_t$ in $t$,  we  only need to prove 
 $\delta P_t\phi$ is continuous   for $\phi$ in $C_K^\infty$  by the uniform 
boundedness principle(see $P_{60}$ in \cite{Dun-Sch}). 

Take  $\phi\in C_K^\infty$. First we have pointwise continuity from the $L^p$ 
continuity  of $|TF_t|\chi_{t<\xi}$:
 $$\lim_{t\to 0} |E\phi(T_xF_t(v_0))\chi_{t<\xi}-\phi(x)(v_0))|^p=0$$
 for  all $x$ in $M$ and $v_0\in T_xM$. 

Next we show $\delta P_t\phi$ converges to $\phi$ in $L^p(M, e^{2h} dx)$.
 Let $t<a$,  then

\begin{eqnarray*}
|E\left(\phi(T_xF_t)\chi_{t<\xi}\right)|^p  
&\le& E\left(|\phi|_{F_t(x)}^p\chi_{t<\xi}\right) \sup_{t\le a}\sup_{x\in M}
\left(E\left(|T_xF_t|^q\chi_{t<\xi}\right)\right)^{p\over q}\\
&=& c^{p\over q}E\left(|\phi|_{F_t(x)}^p\chi_{t<\xi}\right).
\end{eqnarray*}

\noindent
By the invariance property, we also have:

$$\int_M E\left(|\phi|_{F_t(x)}^p\chi_{t<\xi}\right)\,e^{2h}dx
=\int_M |\phi|_x^p e^{2h} dx.$$
But
\begin{eqnarray*}
|\delta P_t\phi -\phi|^p &\le & c_p|\delta P_t\phi|^p + c_p|\phi|^p\\
&\le& c^{p\over q}c_pE\left(|\phi|_{F_t(x)}^p\chi_{t<\xi}\right)+ c_p|\phi|^p.
\end{eqnarray*}
Here $c_p$ is a constant. For the last term we may change order of taking 
expectation and taking limit in $t$. Thus by a standard comparison theorem,
 we have:

$$\lim_{t\to 0} \int_M |\delta P_t\phi -\phi|^p e^{2h}\, dx =
\int_M\lim_{t \to 0} |E\phi(v_t)\chi_{t<\xi}-\phi|^p e^{2h} dx.$$
 So $\delta P_t \phi$ converges to $\phi$ in $L^p$.

 Finally notice with the last assumption, we can prove that $\delta P_t$ is a 
strongly continuous semigroup  both on $L^p$ and $L^\infty$. It is a
 $L^\infty$ semigroup since for $\phi\in L^\infty$,
$$|\delta P_t\phi|_{L^\infty}\le 
\sup_{x\in M} E|T_xF_t\chi_{t<\xi}||\phi|_{L^\infty}.$$
Its strong continuity comes from the last assumption.

 By the Riesz-Thorin interpolation theorem $\delta P_t$ is a semigroup on
 $L^r$ for  $p\le r \le \infty$ . Furthermore:
$$\lim_{t\to 0} \int_M |\delta P_t \phi -\phi|^r e^{2h} dx
\le \lim_{t\to 0} \int |\delta P_t\phi -\phi|^p\,
|\delta P_t\phi -\phi|_{L^\infty}^{r-p} e^{2h} dx =0.$$
End of the Proof. \hfill \rule{3mm}{3mm}

\bigskip

\noindent{\bf Remark:}
This proof works whenever there is an invariant measure for $F_t$.

\begin{corollary}\label{co: Derivative 3}
For a gradient Brownian system with drift $\nabla h$, we have
$$\delta P_t\phi=\heatsemif \phi$$
for $\phi\in C_K^\infty$  if for all $t$ and  some constant $a>0$
$$\sup_M E|T_xF_t|^2\chi_{t<\xi}<\infty,$$
$$\sup_{s\le a}\sup_M E|T_xF_t|^2\chi_{t<\xi}<\infty,$$ 

\noindent
 and  the map $t\mapsto |T_xF_t|\chi_{t<\xi}$ is continuous into $L^2(\Omega,\F,P)$.

 In particular the conditions hold if there is a number $\delta>0$ such that
for all $t$:
\begin{equation}
\sup_{s\le t}\sup_{M} E|T_xF_t|^{2+\delta} \chi_{t<\xi(x)}<\infty.
\end{equation}
\end{corollary}

\noindent
{\bf Proof:} From proposition \ref{pr: Derivative 2}, the semigroup $\delta P_t$ has
generator $\half \triangle^h$ for gradient system on $C_K^\infty$. The result
follows from theorem ~\ref{th: properties of semigroups} on page 
\pageref{th: properties of semigroups}. \hfill \rule{3mm}{3mm}

\bigskip

Recall we defined $H_p$ for gradient h-Brownian system  (section 5.3) as follows:
\[\begin{array}{ll}
H_p(v,v)=&-{\rm Ric}(v,v)+2<{\rm Hess}(h)(v),v>\\
&+|\alpha(v,\cdot)|^2_{H,S}
+{p-2 \over |v|^2} |\alpha(v,v)|^2. \end{array}\]

Following Strichartz, we discuss the $L^p$ boundedness of heat semigroups for forms:

\begin{corollary} Let $M$ be a complete Riemannian manifold with Ricci-2Hess(h)
 bounded from below. Suppose there are  constants $k$ and  $\delta>0$ such that 
$$H_{1+\delta}(v,v)\le k_1|v|^2.$$
 Then $\heatsemi$ is $L^p$ bounded uniformly in $t$ in finite intervals $[0,T]$,
 for all $p$ between ${1+\delta\over \delta}$ and  infinity.
\end{corollary}

\noindent{\bf Proof:} Consider a gradient Brownian system with generator 
$\half \triangle^h$. It has no explosion if Ricci-2Hess(h) is bounded from below
\cite{BA86}.
Let $\phi\in C_K^\infty$, corollary ~\ref{co: Derivative Ricci} on page
 ~\pageref{co: Derivative Ricci}  gives us:
$$\heatsemi \phi=E\phi(v_t).$$

Let $\alpha={1+\delta \over \delta}$ be the conjugate number to $1+\delta$. 
 We have:

\begin{eqnarray*}
\int_M |E\phi(TF_t)|^\alpha e^{2h} dx&\le& \int_M E|\phi|^\alpha 
(E|TF_t|^{1+\delta})^{1\over \delta} e^{2h} dx\\
&\le& \sup_{t\le T} \sup_x \left(E|T_xF_t|^{1+\delta}\right)^{1\over \delta} 
\left(|\phi|_{L^\alpha}\right)^\alpha.
\end{eqnarray*}
However $\sup_{t\le T}\sup_M E|T_xF_t|^{1+\delta} <\infty$ if
 $H_{1+\delta}(v,v)\le k|v|^2$
as in section 5.3.
Thus the $L^p$ boundedness. \hfill \rule{3mm}{3mm}

Next we consider the $L^p$ contractivity of heat  semigroups 
(see  page \pageref{de: contraction} for definition). 
 Define $\delta P_t\phi$ for $k$ forms as follows:

\begin{equation}
\label{semigroups for higher}
\delta P_t\phi(v_0^1,\dots, v_0^k)
=E\phi(TF_t(v_0^1), \dots, TF_t(v_0^k))\chi_{t<\xi}.\end{equation}

\begin{lemma} Consider a  Brownian system with drift $\nabla h$ on a complete 
Riemannian manifold.  Suppose   $\delta P_t\phi=\heatsemi \phi$  for $k$ forms 
$\phi\in L^\infty$.
If there is a constant $\delta>0$ such that
$$\sup_xE|T_xF_t\chi_{t<\xi(x)}|^{k(1+\delta)} \le 1,$$
Then the heatsemigroup $\heatsemi$ on $k$ forms is contractive on $L^p$ 
(on $L^2\cap L^p$), for $p$ between ${1+\delta \over \delta }$ and  $\infty$.
\end{lemma}

\noindent{\bf Proof:} Let $\phi\in C_K^\infty$. 
By the argument in the proof of the last corollary, we have
$$|\heatsemi \phi|_{L^{{1+\delta \over \delta }}}\le |\phi|_{L^{{1+\delta \over \delta }}}.$$
On the other hand we have, for $\phi\in L^\infty$:
$$|\heatsemi \phi|_{L^\infty}\le |\phi|_{L^\infty}.$$
Thus $\heatsemi$ is both $L^{{1+\delta \over \delta }}$ and $L^{\infty}$
contractive by the uniform boundedness principle. Finally we apply the
  Reisz-Thorin
 interpolation theorem to get the required result.  \hfill \rule{3mm}{3mm}

\bigskip

Let $M$ be a compact manifold. Then $\delta P_t\phi=\heatsemif \phi$   for closed $C^2$ 1-form $\phi$,  and  $\delta P_t=\heatsemi $ for $C^2$ forms of all order 
if we are considering  gradient systems.
See \cite{ELflow} for detail. From this we have:

\begin{corollary} \label{stricharts:coho}
Let $M$ be a compact manifold. Consider a gradient Brownian system with
 generator $\half \triangle$. Then

(1). The first real cohomology\index{cohomology: theorem} group vanishes if
 the Ricci curvature is positive definite 
at one point and there is a number $\delta>0$ such that 
$$\sup_M E|T_xF_t|^{1+\delta}\le 1$$
when $t$ large.

(2). The cohomology in dimension $k$ vanishes  for a manifold whose curvature 
operator is positive definite at one point or for a 2-dimensional manifold 
which is not flat if our gradient Brownian system satisfies
$$\sup_ME|T_xF_t|^{k(1+\delta)}\le1.$$
when $t$ large.
\end{corollary}

\noindent{\bf Proof:}
We apply the following theorem from \cite{STRI86}: 
Suppose  the real cohomology in dimension $k$ is not trivial, then
the heat semigroup for $k$ forms is not $L^p$ contractive for $p\not = 2$ in the 
following cases:

(1). $n=2$ and the manifold is not flat.

(2). $k=1$ and the Ricci curvature is strictly positive at one point.

(3). The curvature operator is positive definite at one point.

Note we  have $L^p$ contractivity for some $p\not = 2$ from the assumption,
 thus finishing the proof.  \hfill \rule{3mm}{3mm}

\noindent
{\bf Example:}

Let $M=S^n(r)$ be the n-sphere  of radius $r$ in $R^{n+1}$. Then the second 
fundamental form $\alpha$ is given by:
$\alpha_x(u,v)=-{1\over r^2}<u,v>x$. The derivative flow for the  gradient Brownian
 motion on it satisfy:
$$E|v_t|^p=|v_0|^p+{p\over 2}\int_0^t E|v_s|^{p-2} H_p(v_s, v_s)ds,$$
and $H_p(v,v)=(p-n){|v|^2\over r^2}$.
So
$$E|v_t|^p=|v_0|^p\exp^{{p(p-n)\over 2|r|^2}t}.$$
See Elworthy \cite{ELflow}. Thus
$$\sup_M E|T_xF_t|^p\le \exp^{{p(p-n)\over 2|r|^2}t},$$
which is less or equal to $1$, when $t$ big if $p\le n$.
The above corollary confirms that the kth cohomology vanishes for the n-sphere
if $k<n$. However note that  the $n^{th}$ cohomoloy  of the sphere does not vanish. 

\medskip
    We'll come back to this topic in the next chapter.

\chapter{Consequences of moment stability}

\section{Introduction}

In the first section, we will show the first homotopy group vanishes if the Brownian motion on $M$ is strongly moment stable and satisfies certain regularity conditions, as for  compact manifolds  in  \cite{EL-survey}.  Also interesting here is the result on finite h-volume of the manifold.  In section 2, we look at 
the existence of harmonic functions and at cohomology given strong moment stability.

\section{Geometric consequences and vanishing of $\pi_1(M)$}

We shall first show if  the heat semigroup for 1-forms is 
"continuous" on $L^\infty$, then the h-Brownian motion on $M$ has no explosion. In the  following lemma  we only used  pointwise conditions, improving a result of
Bakry \cite{BA86} proved by essentially the same method.

\begin{lemma}\label{le: homotopy 1}
Let $M$ be a complete Riemannian manifold. Assume there is  a point $x_0\in M$ 
with a neighbourhood $U$ such that for each $x\in U$ there is a constant
 $C_t(x)$ with the following property: 
 $$|\heatsemif df|_x\le C_{t}(x)|df|_\infty, $$
 for  $f\in C_K^\infty$ and $t\le t_0$. Here $t_0$ is a positive constant.
Then the h-Brownian motion does not explode.
\end{lemma}

\noindent
{\bf Proof:}
Let $h_n$ be an increasing  sequence in $C_K^\infty$ such that
$\lim_{n\to \infty} h_n\to 1$, $0\le h_n\le 1$, and $|\nabla h_n|<{1\over n}$.
Such a sequence exists as shown in the appendix. Then 
$\heatsemi h_n(x)\to \heatsemi 1(x)$ for each $x$, since $\heatsemi h_n=P_th_n$ and
by the bounded convergence theorem. Here $\heatsemi 1$ is defined as in theorem 
$~\ref{th: STRana}$. 
By Schauder type  estimate as in the appendix, we have:

$$d\heatsemi h_n(x)\to d(\heatsemi 1)(x)$$
for all $x$ in $M$.  However for $t\le t_0$, and $x\in U$

\begin{eqnarray*}
 |d \heatsemi h_n(x)|
&=& |\heatsemif (dh_n)|_x \\
&\le& C_t(x)|dh_n|_{\infty} \le {C_t(x)\over n}\to 0.
\end{eqnarray*}
Thus $d(\heatsemi 1)(x)\equiv 0$ around $x_0$.
So
 $${\partial \over \partial t} (\heatsemi 1)(x_0)
=\half\triangle^h (\heatsemi 1)(x_0)=0.$$

This gives: $\heatsemi 1(x_0)=1$, for $t<t_0$.  But  $P_{t}1(x_0)=\heatsemi 1(x_0)$ 
from proposition ~\ref{pr: semigroups for functions 2}.  Thus 
$P_{t_0}1(x_0)=P\{t_0<\xi(x_0)\}=1$. Consequently $P_t1\equiv 1$ for all $t$.  Next
 we notice Brownian motion does not  explode if it does not explode at one point. The proof is finished. \hfill \rule{3mm}{3mm}

\bigskip

\begin{proposition} $\label{pr: nonexplosion}$
Consider a h-Brownian motion on a complete Riemannian manifold.
Assume there is a point $x_0\in M$ and a neighbourhood $U$ of $x_0$ such that:
$ E\left(|T_xF_t|\chi_{t<\xi}\right)<\infty$ for all $x\in U$ and $t<t_0$.  Here
$t_0>0$ is a constant. Then
 $P_t1=1$ if $dP_tf=(\delta P_t)(df)$ for $f\in C_K^\infty$.
\end{proposition}
\noindent{\bf Proof:} 
This is a direct consequence of lemma ~\ref{le: Differentiate 10} and the
lemma above. \hfill \rule{3mm}{3mm}

\bigskip

Note: If $dP_tf=(\delta P_t)(df)$ for  $f$ with both $f$ and $df$ bounded,
then set $f\equiv 1$,  the same argument in lemma
~\ref{le: homotopy 1} shows $P_t1=1$.

The proposition  above indicates there is not much point to discuss the possibility
 of $dP_t=(\delta P_t)d$ when there is explosion.  From this there  also arises an 
interesting question, which we have not answered yet:
If $E\left(|TF_t|\chi_{t<\xi}\right)<\infty$ for  $t<t_0$, does it hold for
 all $t$?

\begin{corollary}
\label{co: homotopy 1}
Let $(X, A)$ be a gradient s.d.s. with generator $\half \triangle^h$. Then it has no explosion if its derivative flow $TF_t$ satisfies the following conditions:
\begin{enumerate}
\item For each $t>0$, 
$$\sup_{ M}E|T_xF_t|^2\chi_{t<\xi(x)}<\infty$$
\item  There is a number $a>0$ such that:
$$\sup_{s\le a}\sup_{x\in M} \left(E|T_xF_s|^2\chi_{s<\xi(x)}\right)<\infty,$$
\item
The map $t\to |T_xF_t|\chi_{t<\xi(x)}$ is continuous into $L^2(\Omega, \F ,P)$. 
\end{enumerate}
\end{corollary}
\noindent
{\bf Proof:} This comes from  proposition ~\ref{pr: nonexplosion} and  corollary 
 ~\ref{co: Derivative 3}. \hfill \rule{3mm}{3mm}

\bigskip

 Note that  the conditions in the corollary can be checked in terms of the
extrinsic  curvatures of the manifold (see section ~\ref{applications}). Thus
 this gives us a result
 on nonexplosion of the Brownian motion(not necessarily gradient) on $M$. However it
is not clear that this improves Bakry's result: there is no explosion if
Ricci-Hess(h) is bounded from below.

Following Bakry \cite{BA86}, we get a finite volume result:

\begin{theorem} $\label{th: finite volume}$
Let $M$ be a complete Riemannian manifold.
 Assume the s.d.s. has generator $\half \triangle^{h}$ and  for each $t\ge 0$ and each compact set $K$ we have:

 (1). $dP_tf=(\delta P_t)(df)$,  for $f\in C_K^\infty$.

(2).  $$\int_0^\infty\sup_{x\in K}E|T_xF_s\chi_{s<\xi}|ds <\infty.$$
 Then the  h-volume of the manifold is finite.
\end{theorem}

{\bf Proof:}
Let $f\in C_K^\infty$ with nonempty support $K$.  Then $\lim_{t\to \infty}P_tf$
 exists in $L^2$ and is h-harmonic: $\triangle^h(\lim_{t\to\infty}P_tf)=0$.  This 
comes from the self-adjointness of $\triangle^h$ and an application of  the 
spectral theorem. 
Let $P_\infty f$ be the limit, then $\nabla (P_\infty f)=0$. So $P_\infty f$ must be  a constant.

Assume h-vol$(M)=\infty$, then $P_\infty f$ must be zero. We will prove this is impossible.  Take $g\in C_K^\infty$, then:
$$\int_M(P_\infty f-f)ge^{2h}\, dx
=\lim_{t\to \infty} \int_M(P_tf-f)ge^{2h}\, dx.$$
But
\begin{eqnarray*}
&&\int_M (P_tf-f)ge^{2h}dx 
	 =\int_M (\int_0^t{\partial\over \partial s} P_sf\,ds)g \,e^{2h}dx\\
&&=\half \int_0^t\left(\int_M (\triangle^hP_sf)ge^{2h}dx\right)ds
	=\half \int_0^t\int_M< d(P_sf), dg> e^{2h}dxds\\
&&=\half \int_0^t\int_M <P_s^{h,1}(df),dg> e^{2h}dxds
	= \half \int_0^t\int_M <df, P_s^{h,1}(dg)>e^{2h} dxds\\
&&=\half \int_0^t\int_K <df, P_s^{h,1}(dg)>e^{2h} dx\,ds
\end{eqnarray*}

\noindent
since $\triangle^{h,1}$ is self-adjoint and so the dual semigroup $\left(P_t^{h,1}\right)*$ equals $P_t^{h,1}$. 

Next note the stochastic dynamical system is complete under the assumption from
 proposition ~\ref{pr: nonexplosion}. 
Let
$$C_s=\sup_{x\in K} E\left(|T_xF_s|\right).$$
Thus

\begin{eqnarray*}
|\int_M (P_tf-f)ge^{2h}dx|
&\le& \half \int_0^t  \int_K |\nabla f|\, E\left(|T_xF_s|\right)\,
|\nabla g|_\infty e^{2h} dx\,ds\\
&\le& \half |\nabla g|_\infty\, \left(\int_0^t C_s ds\right) 
 \int_K|\nabla f| e^{2h} dx\\
&\le& \half |\nabla g|_\infty \,  |\nabla f|_{L^1}\, \int_0^t C_s ds.
\end{eqnarray*}
Let $g=h_n$, we get:

$$|\int_M (P_tf-f)h_n \,e^{2h}dx|
\le {1\over 2n} |\nabla f|_{L^1}\int_0^t C_s ds \to 0.$$

\noindent
But $\lim_{n\to \infty} \int_M (-fh_n)e^{2h} dx = - \int_M fe^{2h} dx$, 
from $|fh_n|\le |f| \in L^1$.  And  we can  choose a function 
$f\in C_K^\infty$ with $\int_M fe^{2h} dx\not = 0$. This gives  a
 contradiction. Thus the h-volume of $M$ must be finite. \hfill \rule{3mm}{3mm}

\bigskip

\begin{corollary} Let $M$ be a complete Riemannian manifold. Assume the generator is
$\half \triangle^h$ and there is no explosion. Then if for each $x$ and $t$:  
$$E\sup_{s\le t}|T_xF_s|<\infty,$$
and for each compact set $K$ 
$$\int_0^\infty \sup_{x\in K} E|T_xF_s| ds<\infty,$$
the manifold  has finite h-volume.
\end{corollary}

\noindent{\bf Proof:}
By theorem $~\ref{th: Stflour}$,  $dP_tf=\delta P_tf$ under the assumption. Applying
the above theorem, we get the conclusion.  \hfill \rule{3mm}{3mm}

\begin{corollary}
Let $M$ be a non-compact  complete Riemannian manifold  with
$${\rm Ric_x} > -{n \over n-1} {1\over \rho(x)^2}$$
for eaxh $x\in M$. Then if $d(P_tf)=(\delta P_t)(df)$, (for $F_t$ a Brownian
 motion on $M$),  it cannot be strongly moment stable(defined on page
 ~\pageref{moment stable}). Here
 $\rho$ denotes the distance function on $M$ from a fixed point $\cal O$.
\end{corollary}

\noindent{\bf Proof:}  This is a direct application of the following theorem from 
\cite{CH-GR-TA}: The volume of $M$ is infinite 
for  noncompact manifolds with the above condition on the Ricci curvature.  
But by theorem $~\ref{th: finite volume}$ strong moment stability implies finite 
volume. Note  the Brownian motion here has no explosion by the example on page
~\pageref{ex: infinity Ricci}. \hfill \rule{3mm}{3mm}

\bigskip

\section{Vanishing of $\pi_1(M)$}

\begin{theorem}\label{th: homotopy 2}
Let $M$ be a  Riemannian manifold with its injectivity radius bigger than a
positive number $c$. Assume we have an s.d.s. $(X,A)$ on $M$ which is strongly
1-complete( with $C^2$ coefficients), then the first homotopy group
$\pi_1(M)$ vanishes if for each compact set $K$,

$$\lim_{t\to\infty}\sup_{x\in K} E|T_xF_t|=0.$$

\end{theorem}

\noindent{\bf Proof:}
Take $\sigma$ to be a $C^1$ loop parametrized by arc length. Then 
$F_t\circ\sigma$ is a $C^1$ loop homotopic to $\sigma$ by the strong 
1-completeness. Let $\ell(\sigma_t)$ denote the length of $F_t(\sigma)$.

If we can  show $F_t\circ \sigma$ is contractible to a point in $M$ with
probability bigger than zero for some $t>0$, then the theorem is proved from
the definition:  $\pi_1(M)=0$ if every continuous loop is contractible to one
point.

\noindent
We have:

\begin{eqnarray*}
E\ell(\sigma_t)&=& E\int_0^{\ell_0} |T_{\sigma(s)}F_t(\dot \sigma(s))|ds\\
&\le& \ell_0 \sup_s E|T_{\sigma(s)}F_t|.
\end{eqnarray*}

\noindent
Thus
$$\lim_{t\to \infty} E\ell(\sigma_t)=0.$$

\noindent 
Take $t_0$  such that $E\ell(\sigma_{t_0})<{c\over 2}$. Then
 $\ell(\sigma_{t_0})<{c\over 2}$ with probability bigger than zero.
For such $\omega$ wth $\ell(\sigma_{t_0})(\omega)<{c\over 2}$, 
 $F_{t_0}\circ\sigma(\omega)$ is 
contained in a geodesic ball with radius smaller than $\half c$. Since the
 geodesic ball is  diffeomorphic to a ball in $R^n$, it is  contractible to 
 one point.  Thus $F_{t_0}\cdot \sigma(\omega)$  is contractible to one point
for  a set of  $\omega$ with probability bigger than $0$.  This finished the proof.
\hfill \rule{3mm}{3mm}

\bigskip

\begin{theorem}\label{th: homotopy}
Let $M$ be a complete Riemannian manifold. Consider a s.d.s. $(X,A)$ with generator
 $\half \triangle^h$. Suppose the s.d.s. is  strongly 1-complete and satisfies 
$dP_tf=(\delta P_t)(df)$ for $f\in C_K^\infty$.
Then the first homotopy group $\pi_1(M)$\index{$\pi_1(M)$} vanishes  if 

	$$\int_0^\infty \sup_{x\in K} E|T_xF_t| dt<\infty.$$

%\item
% $\int_0^\infty \limsup_{i\to\infty}\sup_{x\in K} E|TF_{t\wedge S_i^H}(x)|dt<\infty,%$  and 

% $\sup_i\sup_{x\in K} E|TF_{t\wedge S_i^H}(x)|^{1+\delta}<\infty$
%\end{enumerate}
%for some constant $\delta>0$. Here $S_i^H=\inf_{x\in H} S_i(x)$ is as defined before %%theorem
%$~\ref{th: dim1}$.

\end{theorem}

\noindent{\bf Proof:}
Let $\sigma: [0,\ell_0] \to M $ be a $C^1$ loop parametrized by arc length.  Then
 $F_t\circ\sigma$ is again a $C^1$ loop homotopic to $\sigma$ from the strong 
1-completeness. Denote by $\ell(\sigma_t)$ its length. 

 We only need to show $F_t\sigma(\omega)$ is contractible to one point for some
 $\omega$ and for some $t$.

First we claim there is a sequence of numbers $\{t_j\}$ converging to infinity 
such that:

\begin{equation}
E\ell(\sigma_{t_j}) \to 0. 
\label{eq:homo1}
\end{equation}
Since:
\begin{eqnarray*}
\int_0^\infty E\ell(\sigma_t)dt 
 &=&\int_0^\infty E\int_0^{\ell_0} |T_{\sigma(s)}F_t|dsdt\\
&=& \int_0^\infty\int_0^{\ell_0} E|T_{\sigma(s)}F_t|ds dt \\
&\le&\ell_0 \int_0^\infty \sup_s E|T_{\sigma(s)}F_t|dt<\infty.  
\end{eqnarray*}

\noindent
So $\liminf_{t\to \infty} E\ell(\sigma_{t}) =0$, giving $(\ref{eq:homo1})$.
Therefore $\ell(\sigma_{t_j})\to 0$ in probability.

Let $m=\exp^{2h}dx$ be the normalized invariant measure on $M$ for the
 process.    Let $K$ be a compact set in $M$ containing the image set of the
 loop $\sigma$ and  which has measure $m(K)>0$.  Let   $a>0$ be the infimum
 over $x\in K$ of the  injectivity radius at $x$.

By (\ref{eq:homo1}), there is a number $N$ such that for $j>N$, 

$$P\{\ell(\sigma_{t_j})>\half a\} < {m(K) \over 4}.$$
Note h-vol$(M)<\infty$ by theorem $~\ref{th: finite volume}$,  the ergodic 
theorem(see chapter 3) gives: 
$$\lim_{t\to \infty} P\{F_t(x)\in K\}=m(K)$$
for all $x\in M$.

Take a point $\tilde x$ in the image of the loop $\sigma$.  There exists a number $N_1$ such that if $j>N_1$, then: 

$$P\{F_{t_j}(\tilde x)\in K\}  >{m(K) \over 2}.$$
Thus
\begin{eqnarray*}
&&P\{\ell(\sigma_{t_j})<\half a, F_{t_j}(\tilde x)\in K\}\\
&  & =P\{F_{t_j}(\tilde x) \in K\} 
-P\{F_{t_j}(\tilde x)\in K, \ell(\sigma_{t_j})>\half a\}\\ 
&  &\ge P\{F_{t_j}(\tilde x) \in K\} -P\{\ell(\sigma_{t_j})>\half a\}\\
&  & > {m(K) \over 4}.
\end{eqnarray*}

But by the definition of the injectivity radius, there is a coordinate chart 
containing a geodesic ball of radius ${a\over 2}$ around $F_{t_j}(\tilde x) $. 
So the whole loop $F_{t_j}\circ\sigma$  is contained in the same chart with
probability $>{m(k)\over 4}$ (for t big), thus contractible to one point. 
\hfill \rule{3mm}{3mm}

\bigskip

\begin{corollary} Let $M$ be a complete Riemannian manifold. Assume nonexplosion
for the h-BM.
If $ E\left(\sup_{s\le t} |T_xF_s|\right)<\infty$ and
$$\int_0^\infty \sup_{x\in K} E|T_xF_t|dt<\infty$$
 for every compact set $K$, then we have 
$\pi_1(M)=0$. 

In particular if $E\left(\sup_{s\le t} |T_xF_s|\right)<\infty$ then $F_\cdot$ cannot be strongly moment stable given nonexplosion unless $M$ is simply connected.
\label{co: ho1}
\end{corollary}

\noindent
{\bf Proof:} Apply theorem ~\ref{th: dim1}, theorem ~\ref{th: Stflour}, and the 
theorem above.

\begin{corollary}
Let $M$ be a complete Riemannian manifold of finite h-volume. Assume the h-Brownian
 motion on $M$ is strongly 1-complete and satisfies:

$$\int_0^\infty \sup_{x\in K} E|T_xF_t|dt<\infty$$
for each compact set $K$. Then $\pi_1(M)=\{0\}$. 
\end{corollary}

\noindent{\bf Proof:}
This comes from the proof of theorem ~\ref{th: homotopy}. \hfill \rule{3mm}{3mm}

\bigskip

%Furthermore if we have strong completeness, then if 
%$E\left(\sup_{s\le t} |T_xF_s|\right)<\infty$, we cannot have moment stability, 
%due to the continuity of the moment exponents  $\mu_x(1)$.

\bigskip
Here is a corollary which generalizes a result of Elworthy and Rosenberg to 
noncompact manifolds. 
See \cite{EL-homotopy}. 

\begin{corollary}\label{co: homotopy 4} 
 Let $M$ be a complete Riemannian manifold.
Suppose   there is a h-Brownian  system on $M$  such that 
 $|\nabla X|$ is bounded and  $H_1(v,v)<-c^2|v|^2$ for  $c\not =0$. 
Then we have $\pi_1(M)=\{0\}$. In particular if $M$ is a submanifold of $R^n$,
we have $\pi_1(M)=\{0\}$ if its second fundamental form $\alpha$ is bounded, 
and $H_1(v,v)<-c^2|v|^2$.    Here 

\begin{eqnarray*}
H_p(v,v)&=& -\hbox{Ric}(v,v)+2<\hbox{Hess}(h)(v),v>\\
&+& \sum_1^m |\nabla X^i(x)|^2+(p-2)\sum_1^m {1\over |v|^2}<\nabla X^i(v),v>^2
\end{eqnarray*}
\end{corollary}

\noindent
{\bf Proof:} There is no explosion since Ricci-2Hess(h) is bounded from below
 from the assumptions. On the other hand $H_{1+\delta}$ is bounded above
since $|\nabla X|$ is bounded and so $d(P_tf)=(\delta P_t)(df)$ for $f$ in
$C_K^\infty$. See corollary \ref{co: Derivative 30}. The result follows from
theorem \ref{th: homotopy} since the s.d.s. is strongly moment stable from
$H_1(v,v)\le -c^2|v|^2$. The second part of the theorem follows from  the fact
that  the  sum of the last two terms in the formula for $H_p$ is 
$|\alpha(v,\cdot)|_{H,S}^2 +(p-2){1\over |v|^2} |\alpha(v,v)|^2$ for gradient
 Brownian systems.

\section{Vanishing of harmonic forms and cohomology}

We come back to the discussion on cohomology vanishing of page 
$~\pageref{stricharts:coho}$ and aim to extend some of the results in
 \cite{ELflow} on cohomology vanishing  given moment stability to noncompact 
manifolds. See  \cite{EL-RO88},\cite{ELflow}.

Let $C^\infty(\Omega^p)$ be the space of $C^\infty$ smooth $p$ forms on $M$.  A p-form $\phi$ is closed if $d\phi=0$,  exact if $\phi=d\psi$ for some p-1 form $\psi$. 
Here $d$ is the exterior differentiation defined in section $~\ref{differential forms}$. A $h$-{\it harmonic} form is a form with $\triangle^h\phi=0$. 
The $p^{th}$ {\it de Rham cohomology}\index{cohomology: definition} group $H^p(M,R)$ is defined to be the quotient group of the group of smooth closed $p$ forms by the group of $C^\infty$ exact forms:
$$H^p(M,R)
={{\rm Ker}\left(d: C^\infty(\Omega^p)\to C^\infty(\Omega^{p+1})\right)
\over {\rm Im}\left(d: C^\infty(\Omega^{p-1})\to C^\infty(\Omega^p)\right)}.$$
There is also the cohomology group $H_K^p(M,R)$ with compact supports:
$$H^p_K(M,R)
={{\rm Ker}\left(d: C^\infty_K(\Omega^p)\to C^\infty_K(\Omega^{p+1})\right)
\over
 {\rm Im}\left(d: C_K^\infty(\Omega^{p-1})\to C_K^\infty(\Omega^p)\right)}.$$

Let $h$ be a $C^\infty$ smooth function. There is the Hodge decomposition theorem(
see page $~\pageref{Hodge decomposition}$):
$$L^2\Omega^p=\overline{{\rm Im}(\delta^h)}\bigoplus 
\overline{{\rm Im}(d)}\bigoplus L^2(\cal H).$$
Let $\phi$ be a form  with $d\phi=0$ and decomposition:
$\phi=\alpha +\beta +H\phi$. 
Here $\alpha\in \overline{{\rm Im}(\delta)}$ and $\beta \in \overline{{\rm Im}(d)}$.
Then $\alpha=0$ since $d\alpha=0$ by $d\beta=0$ and $d(H\phi)=0$. 

Thus we have the Hodge's theorem:
every cohomology class has a unique harmonic representative if Im(d) is closed,
 in particular:
the dimension of $H^p(M,R)$, as a linear space, equals the dimension of the space of h-harmonic $L^2$ $p-$forms when $d$ has closed range.

The vanishing problem have been studied by conventional methods by e.g. Yau
 \cite{Yau76}.  The problem has been considered in a probabilistic context,
 see e.g. Vauthier \cite{Vauthier79}, Elworthy and
 Rosenberg \cite{EL-RO88},\cite{EL-RO91}.
The idea  we are using here is that the probabilistic semigroup $\delta P_t$ on
 forms often agrees with the heat semigroup (see chapter 6). Thus the
 existence of harmonic forms is  directly related to the  behaviour of
  diffusion processes and their derivatives.  In the following we follow
 Elworthy and Rosenberg's approach to get vanishing results for harmonic
 1-forms.  But we use $\delta P_t$ instead of 
the standard probabilistic formula obtained from the Weitzenbock formula.
However we do not intend to include
 all the possible results in this thesis, but only demonstrate the idea.
 A second theorem we give here follows from an approach of Elworthy
\cite{EL-survey}. This approach uses  integration of p-forms along singular
 p-simplices and fits very well with our definition of  
strong  p-completeness.

\begin{proposition}\label{pr: topo 1}
Let $M$ be a complete Riemannian manifold, consider a s.d.s. with generator
$\half \triangle^h$. 

(1).	Assume $\heatsemif\phi =\delta P_t\phi$ for closed 1-forms $\phi\in L^q$. Let
$p$ be the conjugate number to $q$. 
Then there are no nonzero $L^p$ h-harmonic 1-forms, if 
\begin{equation}
\lim_{t\to \infty} {1\over t} \log \sup_{x\in M} 
E\left(|T_xF_t|\chi_{t<\xi}\right)^{q}<0.
\label{eq: coh1}
\end{equation}

(2). Assume $\heatsemif \phi =\delta P_t\phi$ for closed bounded 1-forms
 $\phi$. 
 Then there  are  no bounded  h-harmonic 1-forms if for each $x\in M$
\begin{equation}
\lim_{t\to \infty}{1\over t} \log E\left(|T_xF_t|\chi_{t<\xi}\right)<0. 
 \label{eq: coh2}
\end{equation}
\end{proposition}

\noindent{\bf Proof:}
(1).	We have nonexplosion and finite volume in the first case according to 
the nonexplosion result on page $~\pageref{pr: nonexplosion}$. Let $\phi$ be
 a nonzero harmonic p-form in $L^p$. Then there is a point $x_0\in M$ with
 $\phi(x_0)\not =0$. Thus $\int |\phi|_x^p\,e^{2h} dx >0$ by continuity. So

\begin{eqnarray*}
0&=& \lim_{t\to \infty} {1\over t} \log \int |\phi|_x^p e^{2h} dx\\
&=& \lim_{t\to \infty} {1\over t} \log \int |\delta P_t\phi|_x^p\, e^{2h} dx\\
&=&\lim_{t\to \infty} {1\over t} \log \int |E\phi(T_xF_t)|^p\, e^{2h}dx\\
 &\le&\lim_{t\to \infty} {1\over t} \log \int
E|\phi|_{F_T(x)}^p\left(E|T_xF_t|^q\right)^{p\over q}e^{2h} dx\\
&\le& \lim_{t\to \infty} {1\over t} \log \sup_{x\in M}
 \left(E|T_xF_t|^q\right)^{p\over q}+\lim_{t\to \infty} {1\over t} \log 
\int E|\phi|_{F_t(x)}^p e^{2h} dx\\
&\le& \lim_{t\to\infty} {p\over q} {1\over t}\log \sup_{x\in M} 
E|T_xF_t|^q + \lim_{t\to \infty} {1\over t} \log \int |\phi|_x^p e^{2h} dx.
\end{eqnarray*}

\noindent
But $ \int|\phi|_x^p e^{2h}dx<\infty$,  giving a contradiction. 

\bigskip

(2).	The proof of  the second part is just as before.  First note we have
 nonexplosion. Let $\phi$ be a closed bounded harmonic 1-form. Let $x_0\in M$
 with $|\phi|_{x_0}\not = 0$. Then:

\begin{eqnarray*}
0&=&\lim_{t\to\infty} {1\over t}\log |\phi(x_0)|
=\lim_{t\to\infty} {1\over t}\log |\heatsemif \phi(x_0)|\\
&=&\lim_{t\to\infty}{1\over t} |E\phi(T_{x_0}F_t)|\\
&\le&\lim_{t\to \infty} {1\over t} \log |\phi|_{\infty} E\left(|T_{x_0}F_t|\right)\\
&\le&\lim_{t\to\infty} {1\over t} \log E|T_{x_0}F_t|.
\end{eqnarray*}

\noindent
But this is impossible from the assumption. End of the proof. \hfill \rule{3mm}{3mm}
\vspace{6mm}

\noindent{\bf Remarks} 
(1). Let $p<2$. Assume $(~\ref{eq: coh1})$.  Then   $\delta P_t=\heatsemif$ on
 $C_K^\infty$ implies $\delta P_t=\heatsemif$ on $ L^p$. 
Since in this case  $q={p\over p-1}>p$, so $\delta P_t$ is a strongly
continuous $L^p$ semigroup from proposition ~\ref{pr: Derivative 3}. We may then 
apply uniform boundedness principle.

(2). If we know that $M$ has finite h-volume to start with, then all bounded
 harmonic 
1-forms vanishes if  $E\left(sup_{s\le t} |T_xF_s|\right)<0$ and if $F_t$ is 
strongly moment stable from theorem ~\ref{th: dim1}.

Following Elworthy \cite{EL-survey} for the compact case, we have the following 
theorem:

\begin{theorem}\label{th: p-complete cohomology}
Let $M$ be a Riemannian manifold and assume there is a strongly p-complete s.d.s. 
with strong $p^{th}$-moment stability. Then all bounded closed p-forms are exact. 
In particular the natural map from $H_K^p(M,R)$ to $H^p(M,R)$ is trivial.
\end{theorem}

\noindent{\bf Proof:}
Let $\sigma$ be a singular p-simplex, and  $\phi$ a bounded closed p-form. We
 shall not distinguish  a singular simplex map from its image.
 Denote by $F_t^*\phi$ the pull back of the form  $\phi$ and 
$(F_t)_*\sigma=F_t\circ\sigma$. Then
$$\int_{(F_t)_*\sigma} \phi =\int_{\sigma} (F_t)^*\phi$$
by definition.
But $(F_t)_*\sigma$ is homotopic to $\sigma$ from the strong p-completeness. Thus:
$$\int_\sigma\phi=\int_{(F_t)_*\sigma} \phi.$$
This gives: 
$$\int_\sigma \phi =\int_\sigma (F_t)^*\phi.$$
Take expectations of  both sides to obtain:
\begin{eqnarray*}
E|\int_\sigma\phi| &=& \lim_{t\to \infty} E|\int_\sigma (F_t)^*\phi| \\
&\le&|\phi|_\infty \lim_{t\to\infty} \,\int_\sigma E|TF_t|^p  \\
&\le&|\phi|_\infty\lim_{t\to\infty} \sup_{x\in \sigma} E|TF_t|^p\\
&=&0
\end{eqnarray*}
from the  strong $p^{th}$ moment stability. 
Thus $\int_\sigma \phi=0$, and so  $\phi$ is exact by deRham's theorem.
 \hfill\rule{3mm}{3mm}

\begin{corollary}
Let $M$ be a complete Riemannian manifold. Assume there is a complete s.d.s.
on $M$  with strong $p^{th}$ moment stability and satisfying 
$\sup_{x\in K}E\left(\sup_{s\le t}|T_xF_s|^{p+1}\right)<\infty$ for each 
compact set $K$. Then all bounded closed p-forms  are exact. 
In particular Suppose  $M$ is
a closed submanifold of $R^m$ with its second fundamental form  bounded.
Then if  $H_p(v,v)\le -c^2|v|^2$ for some constant $c$, then the conclusion 
holds. Here  $H_p$ is as defined on  page ~\pageref{eq: H for BM}:

\begin{eqnarray*}
H_p(v,v)&=& -Ric(v,v)+2<\hbox{Hess}(h)(v),v>\\
 & & +|\alpha(v,-)|_{H,S}^2+(p-2){1\over |v|^2} |\alpha(v,v)|^2.
\end{eqnarray*}

\end{corollary}
\noindent{\bf Proof:}
Direct applications of the above theorem and theorem 
\ref{th: strong completeness}.

\section{Examples}

\noindent
{\bf Example 1} Let $M=R^n$, $h(x)=-|x|^2$. Then h-vol$(R^n)<\infty$. 
Furthermore we have:
$$H_1(v,v)=2<\hbox{Hess}(h)(v),v>=-2|v|^2.$$
Thus the s.d.s. on $R^n$:
$$dx_t=dB_t-\nabla h(x_t)dt=dB_t-2x_tdt$$
is strongly complete and strongly $p^{th}$ moment stable. 

\bigskip

\noindent{\bf Example 2} Let $B_t^1$, $B_t^2$ be independent Brownian motions
 on $R^1$. Then 
$$F_t(x,y)=(x+\int_0^t \exp^{B_s^2-{s\over 2}} dB_s^1, y\exp^{B_t^2-{t\over 2}})$$
is a Brownian flow on the hyperbolic space $H^2$. It is strongly complete and 
satisfies: 
$$E\sup_{s\le t}|T_xF_s|<\infty.$$
So $\delta P_t\phi=\heatsemi \phi$ for bounded 1-forms. 
Thus this Brownian system on $H^2$ is not strongly moment stable since 
$H^2$ 
has infinite volume(c.f. theorem ~\ref{th: homotopy}).
%(c.f. proposition ~\ref{pr: topo 1}).

\bigskip

\noindent  {\bf Example 3}\label{ex: Langevin}
Consider the Langevin equation\index{Langevin equation} on $R^2$:

$$dx_t=\gamma dB_t -cx_tdt.$$
\noindent
Here  $\gamma$ and $c$ are constants.  The solution can be written down explicitly:

$$x_t=x_0\exp^{-ct} +\gamma \int_0^t {\rm e}^{-c(t-s)}\, dB_s.$$

\noindent
It has Gaussian distribution and generator $\half\triangle^h$ for $h=-{cx^2 \over 2}$, and has no explosion. Its derivative flow is given by:
$$TF_t(v)=\exp^{-ct}v$$
and enjoys the following properties:

\begin{enumerate}
\item
 $$\sup_x E\left(\sup_{s\le t}|T_xF_s|^2\right)=1,$$
\item
Strong 1-completeness from theorem ~\ref{th: dim1},
\item
$$d(P_tf)=(\delta P_t)(df)$$
for $f\in C_K^\infty$ from theorem ~\ref{th: differentiate semigroups},
\item
$$\int_0^\infty \sup_x E|T_xF_s| ds =\int_0^\infty \exp^{-cs} ds=1,$$
\item
The solution process is recurrent and has $e^{2h} dx$ as  finite invariant
 measure, since h-vol$(R^n)<\infty$.
\end{enumerate}
If we consider the same equation on $M=R^2-\{0\}$ instead of on $R^2$, then
 all the
properties hold except for the strong 1-completeness as shall be shown below.

Clearly part 1 and part 4 hold. The solution is recurrent on $R^2-\{0\}$ since
it is recurrent on $R^2$. Furthermore $e^{2h}dx$ is still the invariant
 measure since  $R^2-\{0\}$ has negligible boundary and from the completeness
 of the s.d.e. on  $R^2-\{0\}$. With these the conclusion of
 proposition ~\ref{pr: ergodic property}  certainly holds. 

Suppose the process is strongly 1-complete. The ergodic property gives us
$\pi_1(M)=\{0\}$ from the proof of the theorem ~\ref{th: homotopy} and part 4 of
the properties. This is a contradiction. Thus we do not have the strong 
1-completeness.

Finally part 3 holds since on $R^2-\{0\}$, $(\delta P_t)(df)(v)=Edf(v) $ and 
$P_tf(x)=Ef(x+B_t)$  as on $R^2$.

\vspace{10mm}

In the following we look at some surfaces whose second fundamental forms
are bounded. We will show that the
 Hyperboloid satisfies our hypothesis which imply  $\pi_1(M)=0$, and Brownian
 motions on  both the torus and the cylinder cannot be strongly moment stable.
First we recall the basic theory:

Let $M$ be a surface in $R^3$ parametrized by $x=x(u,v)$. The unit normal
vector $\mu$ is given by:
$$\mu={x_u\times x_v\over |x_u\times x_v|}.$$
There is the shape operator $S\colon T_xM\to T_x M$ given by:
$$S(v)=\colon -D_v\mu.$$
Here $D$ denotes the differentiation on $R^3$.  The second
fundamental form 
$II(u,v)=-l(u,v)\mu$ is given in terms of the shape operator $S$:
\begin{equation}
l(u,v)=-<S(u),v>.
\end{equation}
Let 

\[ \begin{array}{lll}
e=<\mu,x_{uu}> \hskip 6pt &f=<\mu, x_{uv}> \hskip 6pt &g=<\mu, x_{vv}>\\
E=<x_u,x_u> \hskip 6pt &F=<x_u,x_v>  \hskip 6pt &G=<x_v,x_v>.
\end{array}\]

There is then the Weingarten equation:

\[\begin{array}{lll}
-S(x_u)=&{fF-eG\over EG-F^2}x_u +&{eF-fE\over EG-F^2} x_v\\
-S(x_v)=&{gF-fG\over EG-F^2}x_u+&{fF-gE\over EG-F^2}x_v.\\
\end{array}\]

%%Thus  we have:
%%\[\begin{array}{lll}
%%<S(x_u),x_u>=&{-fF+eG\over EG-F^2} E+ &{-eF+fF\over EG-F^2} F\\
%%<S(x_u),x_v>=&{-fF+eG\over EG-F^2}F+&{-eF+fE\over EG-F^2}G\\
%%<S(x_v), x_v>=&{fG-gF\over EG-F^2}F+&{gE-fF\over EG-F^2}G.\\
%%\end{array}\]

\noindent{\bf Example 4 }[Surface of revolution]

Consider the surface given by:
$$(c_1(s)\cos\theta, c_1(s)\sin\theta, c_2(s)).$$
For this surface:
$$E=\left[c_1(s)\right]^2, 
\hskip 6pt F=0,  \hskip 6pt G= \left\{\left[c_1\prime(s)\right]^2+\left[c_2\prime(s)\right]^2\right\}^2,$$
\[\begin{array}{ll}
e=-&{c_1(s)c_2\prime(s)\over \sqrt{[c_1\prime]^2+c_2\prime]^2}},\\
f=&0,\\
g=&{-c_1\prime(s)c_2\prime \prime(s) +c_1\prime\prime(s)c_2\prime(s)
\over \sqrt{[c_1\prime(s)]^2+[c_2\prime(s)]^2}}.\\
\end{array}\]
\noindent
So 
\[ S\left(\begin{array}{c}
{\partial f\over \partial \theta}\\{\partial f \over
\partial s}\end{array}\right)
=\left(\begin{array}{cc}
K_1 &0\\
0 &K_2\end{array}\right)
\left(\begin{array}{c}
{\partial f\over \partial \theta}\\ {\partial f \over \partial s}
\end{array}
\right).\]
Here 
\begin{eqnarray*}
&&K_1={e\over E} =  {-c_2\prime(s)\over c_1(s)\sqrt{[c_1\prime]^2+[c_2\prime]^2} },\\
&&K_2={g\over G}= {-c_1\prime(s)c_2\prime \prime(s) +c_1\prime\prime(s)c_2\prime(s)
\over \sqrt{[c_1\prime(s)]^2+[c_2\prime(s)]^2}^3}.\\ 
\end{eqnarray*}
The normal vector is:
$$({c_2\prime \cos\theta\over \sqrt{[c_1\prime]^2+[c_2\prime]^2}},
{c_2\prime \sin\theta \over \sqrt{[c_1\prime]^2+[c_2\prime]^2}},
-{c_1\prime \over \sqrt{[c_1\prime]^2+[c_2\prime]^2}}).$$

\bigskip

\noindent{\bf Example 4a} [The Hyperboloid] 

 We will show the surface
$$z^2-(x^2+y^2)=1$$
satisfies the conditions of theorem ~\ref{th: strong completeness Hessian}
 and theorem  ~\ref{th: homotopy}.
Consider the following parametrization: $(s\cos\theta, 
s\sin\theta, \sqrt{s^2+1})$.

\noindent Thus
$E=s^2$, $F=0$, $G={2s^2+1\over s^2+1}$,
$e=-{s^2\over \sqrt{1+2s^2}}$, $f=0$, $g=-{1\over (1+s^2)\sqrt{1+2s^2}}$ .

The unit normal vector is:
$$\mu=({s\cos\theta \over \sqrt{1+2s^2}}, {s\cos\theta \over \sqrt{1+2s^2}},
-{\sqrt{1+s^2}\over \sqrt{1+2s^2}}).$$
Also the Ricci curvature is given by: Ricci$(v)=K_1K_2|v|^2$, while
$K_1=-{1\over \sqrt{1+2s^2}}$, $K_2=-{1\over \sqrt{1+2s^2}^3}$.
Clearly the second fundamental form is bounded, thus
  the Brownian motion on the surface is strongly complete. 
Next we construct a Brownian motion with drift which is strongly moment
stable and thus verify $\pi_1(M)=0$.

\noindent
According to section ~\ref{applications}:

\begin{eqnarray*}
<\nabla X^i(v),v>&=&<A_x(v,<e_i,\mu>\mu) , v>\\
&=& -<e_i,\mu>\l(v,v).
\end{eqnarray*}

Let $h=-c|x|^2$.  Let grad$h$ and $\tilde{{\rm Hess}}(h)$ denote the gradient 
 and the  Hessian of $h$  in $R^3$ respectively.  Then for $x=(x_1, x_2,x_3)\in M$,

\begin{eqnarray*}
<\hbox{Hess}(h)(v),v>&=&<\nabla X(v)({\rm grad}h),v> +<\tilde{{\rm Hess}}h(v),v>\\
&=&-2c<\nabla X(v)(x),v>-2c|v|^2\\
&=&-2c \sum_1^3 x_i<e_i,\mu>l(v,v)-2c|v|^2\\
&=&{2c\over \sqrt{1+2s^2}}l(v,v)-2c|v|^2.\\
\end{eqnarray*}

So 
$$<\hbox{Hess}(h)(v),v>={2c\over \sqrt{1+2s^2}} l(v,v) -2c|v|^2.$$ 
Since $l(v,v)$ is negative and bounded, we may choose $c$ big such that
\begin{eqnarray*}
H_1(v,v)&=&-\hbox{Ric}(v,v)+|l(v,\cdot)|^2 
-{l(v,v)^2\over |v|^2} +2<\hbox{Hess}(h)(v),v>\\
&=&-K_1K_2|v|^2+|l(v,\cdot)|^2 
-{l(v,v)^2\over |v|^2} +{4c\over \sqrt{1+2s^2}}l(v,v)-4c|v|^2\\
&\le&-|v|^2.
\end{eqnarray*}
The system with the chosen drift is then strongly moment stable, so satisfies the
conditions of theorem ~\ref{th: homotopy}.

\noindent{\bf Example 4b}
The  torus given by 
$$((a+b\cos v)\cos u, (a+b\cos u)\sin u, b\sin u)$$
has:
$E=(a+b \cos u)^2$, $F=0$, $G=b^2$,  $e=\cos u (a+b\cos u)$, $f=0$,

\noindent  $g=b$,
$K_1=-{\cos u\over a+b\cos u}$, and $K_2=-{1\over b}$.
So the Brownian motion here is strongly complete, but it cannot be strongly moment 
stable since $\pi_1(M)$ $\not =$$0$. In fact for $a=b={1\over \sqrt{2}}$, the first
 moment is identically $1$ as calculated, e.g. in \cite{ELflow}.

\bigskip

\noindent{\bf Example 4c}
The cylinder $S^1\times(-\infty, \infty)$ parametrized by
$(\cos \theta, \sin \theta, s)$ has 
$K_1=-1$ and $K_2=0$.  The normal vector is: $$\mu=(\cos \theta, \sin\theta,0).$$
As for torus, the Brownian motion here is strongly
complete but not strongly moment stable. To convince ourself we will
try to add a drift as for the Hyperboloid. Let $h=-c|x|^2$.
Then
\begin{eqnarray*}
<\hbox{Hess}(h)(v),v>&=&-2cl(v,v)-2c|v|^2\\
&=&2c|v_1|^2-2c|v|^2=-2c|v_2|^2.
\end{eqnarray*}
Here $v=(v_1,v_2).$ Now
$$H_1(v,v)=v_2^2({v_1^2 \over |v|^2} -2c)
=|v_2|^2{(1-2c)v_1^2-2cv_2^2 \over v_1^2+v_2^2}.$$
Clearly the same argument does not work.

\chapter{Formulae for the derivatives of the solutions of the heat equations}

 \section{Introduction}

In chapter 6, we examined carefully $\delta P_t$ and the heat semigroup $\heatsemif$
 for one forms, and obtained some conditions to ensure $\delta P_t=\heatsemif$. In 
fact we may expect more once we know $\delta P_t$ does agree with $\heatsemif$. The
 semigroup $\heatsemif$ can be given in terms of the line integral 
~\ref{eq: Bismut 2} and a martingale, following from Elworthy \cite{ELflow} for
 compact manifolds. In particular we have a formula for the gradient of the 
logarithm of the heat kernel, extending Bismut's formula \cite{Bismut}. 
 See \cite{DA-EL-ZA} for  an infinite dimensional version of the formulae 
 by Da Prato, Elworthy, and Zabczyk (comment: their proof was inspired by the proof presented in this chapter,
 and follows from an earlier  private communication). 
 See also Norris \cite{Norris90} for another
approach.

The discussions for one forms also work well for higher order forms. We define
 $\delta P_t$ for $q$ forms as on page ~\pageref{semigroups for higher}, and
look briefly the relation between $\delta P_t$ and the heat semigroup $\heatsemi$
for $q$ forms. In the end we give a formula for the exterior derivative of
 $\heatsemi \phi$ in terms of the derivative flow of a h-Brownian motion and
$\phi$ itself.

In the following we write $P_t^h\phi=\heatsemi \phi$. If $\phi$ is a q-form,
we may use $P_t^{h,q}\phi$ instead of $P_t^h\phi$.

\section{For 1-forms}

Let $\phi$ be a 1-form, we define $\int_0^t \phi\circ dx_s$ to be the line 
integral of $\phi$ along Brownian paths as in \cite{ELflow}:
\begin{equation}\label{eq: Bismut 2}
\int_0^t\phi\circ dx_s
=\int_0^t\phi(X(x_s) dB_s) -\half \int_0^t\delta^h\phi(x_s)\, ds.
\end{equation}

\bigskip
Here is the formula for 1-form, which is a direct extension of the formula in 
\cite{ELflow} for compact manifolds:

\begin{proposition} \label{pr: Bismut 1}
Let $(X,A)$ be a complete $C^2$ stochastic dynamical system 
 on a complete
Riemannian manifold $M$ with generator $\half \triangle^h$.	
Suppose  for   closed 1-form  $\phi$ in $D(\bar\triangle^h)\cap L^\infty$,
$$(\delta P_t)\phi=\heatsemif \phi$$
 and  for each $x\in  M$, 

$$\int_0^t E|T_xF_s|^2 ds<\infty.$$ 

Then

\begin{equation}
\label{eq: Bismut for 1-forms}
P_t^{h,1}\phi(v_0)={1\over t} E\int_0^t \phi\circ dx_s\int_0^t
<X(x_s)dB_s, TF_s(v_0)>
\end{equation}

\noindent
 for all $v_0\in T_{x_0}M$.
\end{proposition}

\noindent{\bf Proof:}  
Following the proof for a compact manifold as in \cite{ELflow}.  Let
\begin{equation}
Q_t(\phi)=-\half \int_0^t P_s^{h}(\delta^h \phi)ds. 
\label{eq: Bismut function1}
\end{equation}

\noindent
Differentiate  equation $~\ref{eq: Bismut function1}$ to get:

$${\partial \over \partial t}Q_t\phi=-\half P_t^{h}(\delta^h\phi).$$
We also have:

\begin{eqnarray*}
d(Q_t\phi)&=&-\half \int_0^t d\delta^h (P_s^h\phi) ds\\
&=&\half \int_0^t \triangle^h(P_s^h \phi)ds\\
&=&P_t^{h}\phi -\phi
\end{eqnarray*}

\noindent
since $d\delta^h(P_s^h\phi)=P_s^{h}(d\delta^h \phi)$ is uniformly continuous 
in $s$ and  $$d(P_s^{h} \phi)=P_s^{h} d\phi=0$$
from proposition  $~\ref{pr: exchangeability}$.  
Consequently:

$$\triangle^h(Q_t(\phi))=-P_t^h(\delta^h\phi)+\delta^h\phi.$$

\noindent
 Apply It\^o  formula to $(t,x)\mapsto Q_{T-t}\phi(x)$, which is smooth 
because $P_s^h\phi$ is, to get:

\begin{eqnarray*}
Q_{T-t}\phi(x_t)&=&Q_T\phi(x_0)+\int_0^t d(Q_{T-s}\phi)(X(x_s)dB_s)\\
&  & +\half \int_0^t \triangle^{h} Q_{T-s}\phi(x_s)ds +\int_0^t{\partial
 \over \partial s}Q_{T-s}\phi(x_s)ds\\
&=& Q_T\phi(x_0)+\int_0^t P_{T-s}^{h} (\phi) (X(x_s)dB_s)-\int_0^t\phi\circ dx_s.
\end{eqnarray*}

\noindent
Let $t=T$. We obtain:

$$\int_0^T\phi\circ dx_s =Q_T(\phi)(x_0)+\int_0^TP_{T-s}^{h}(\phi)(X(x_s)dB_s),$$

\noindent
and thus 

$$E\int_0^T\phi\circ dx_s\int_0^T<X(x_s)dB_s, TF_s(v_0)>=E\int_0^T P_{T-s}^{h} \phi(TF_s(v_0))ds.$$

\noindent
But

\begin{equation}
E\int_0^T P_{T-s}^{h} \phi(TF_s(v_0))ds
=\int_0^T EP_{T-s}^{h} \phi(TF_s(v_0))ds.
\label{eq: 1-form 2}
\end{equation}

\noindent
Since by Fubini's theorem we only need to show $EP_{T-s}^{h}\phi(TF_s(v_0))$ is integrable with respect to the double integral:

\begin{eqnarray*}
\int_0^T E|P_{T-s}^{h} \phi (TF_s(v_0))| ds
&\le& |\phi|_\infty \int_0^T E|TF_T(v_0)| ds<\infty.
\end{eqnarray*}

\noindent
Next notice:

$$E\left( P_{T-s}^{h}\phi((TF_s(v_0))\right)=E\phi(TF_T(v_0))=P_T^{h}\phi(v_0).$$

\noindent
from the strong Markov property. We get: 
$$P_T^{h,1}\phi(v_0)={1\over T} E\left\{\int_0^T\phi\circ dx_s\int_0^T
<XdB_s, TF_s(v_0)>\right\}.$$
End of the  proof. \hfill \rule{3mm}{3mm}

\bigskip
\noindent{\bf Remark:}\label{remark}
If we assume  $\sup_x E|T_xF_t|^2<\infty$ for each $t$ in the proposition, we do 
not  need to assume $\phi\in L^\infty$. Since first  
 we have $\delta P_t\phi=\heatsemi \phi$ for such $\phi$ by the uniform boundedness
 principle and  also
equation (~\ref{eq: 1-form 2})  holds from the following argument:
 
\begin{eqnarray*}
&&\int_0^T E |P_{T-s}^{h}\phi (TF_s(v))| ds\\
&&= \int_0^T E|\{E\{\phi(TF_{s,T})|\F_s\}(TF_s(v))|  ds\\
&&\le\int_0^T E|\phi(TF_T(v))| \, ds\\
&&\le \sup_x E|T_xF_T(v)|^2\left(\int_0^T E|\phi|_{F_T(x)}^2 ds\right) \\
&&\le T\  E\left(|\phi|_{F_T(x)}^2\right) \sup_x E|T_xF_T(v)|^2<\infty.
\end{eqnarray*}

\noindent
But  $\int_ME|\phi|_{F_T(x)}^2 e^{2h} dx=\int |\phi|^2 e^{2h} dx<\infty$. So 
$E|\phi|_{F_T(x)}<\infty$ since $E|\phi|_{F_T(x)}^2$ $=P_T(|\phi|^2)(x)$ is
 continuous  in $x$. Thus we may still apply
Fubini's theorem to get: (~\ref{eq: 1-form 2}).

\bigskip

When $\phi=df$ for some function $f$, formula ~\ref{eq: Bismut for 1-forms}
 may be rewritten as:

\begin{equation}
d(P_t^{h}f)(v_0)=\frac 1 t E f(x_t)\int_0^t \langle TF_s(v_0), X(x_s)dB_s\rangle.
\label{eq: formulae in 1}
\end{equation}

In fact this works in a more general situations.  Here is a very intuitive proof by D. Elworthy and myself 
(let $BC^1$ be the space of bounded functions with  bounded continuous first derivative.)

\begin{theorem}\label{th:elementary}
Let $(X,A)$ be a complete nondegenerate stochastic dynamical system so
there is a right inverse map $Y(x)$ for   $X(x)$ each $x$ in $M$. 
Let $f$ be in $BC^1$ s.t.  $(\delta P_t)(df)=d(P_tf)$. 
Then for $v_0\in T_{x_0}M$:

\begin{equation}\label{formula}
d(P_tf)(v_0)={1\over t} Ef(x_t) \int_0^t \langle dB_s, Y(TF_s(v_0))\rangle
\end{equation}

\noindent
provided $\int_0^t <dB_s, Y(T_xF_s(v_0))>$ is a martingale for all $t$. 
Here $P_tf$ is a solution to ${\partial \over \partial t}=\half \sum X^iX^i+A$
with initial value $f$.
\end{theorem}

\noindent
{\bf Proof:}
Let $T>0$.  Apply  It\^o formula to the smooth map $(t,x) \mapsto P_{T-t}f(x)$:
 
$$P_{T-t}f(x_t)=P_Tf(x_0)+\int_0^t d P_{T-s}f(x_s)(XdB_s).$$
Letting t=T, we have:
$$f(x_T)=P_Tf(x_0)+\int_0^T d P_{T-s}f(x_s) (XdB_s).$$
So:
\begin{eqnarray*}
&&E f(x_T)\int_0^T \langle dB_s, Y(TF_s(v_0))\rangle \\
&=&E \int_0^Td P_{T-s}f\left(TF_s(v_0)\right)ds\\
&=& E\int_0^T(\delta P_{T-s}) df(TF_s(v_0))ds\\
&=& \int_0^T(\delta P_T) df(v_0)ds =T\ d(P_Tf)(v_0).
\end{eqnarray*}
End of the proof. \hfill \rule{3mm}{3mm}
\bigskip

Note $\int_0^t E \langle dB_s, Y(T_xF_s(v))\rangle $ is a martingale  if 
$$\int_0^t E|Y(T_xF_s(v))|^2ds<\infty.$$
In terms of the metric on $M$ determined by $Y$, this condition becomes:
$$\int_0^t E|T_xF_s|^2 ds<\infty.$$

\begin{corollary}
Let $p_t^{h}(x,y)$ be the heat kernel as defined on page ~\pageref{th: STRana}, 
then
$$\nabla \log p_t^h(\cdot, y_0)(x_0)
=E\{{1\over t} \int_0^t \left(TF_s\right)^*\left(X(x_s)dB_s\right) |x_t=y_0\}$$
under the assumptions of  proposition ~\ref{pr: Bismut 1}.
 \end{corollary}

\noindent
{\bf Proof:}  The proof is just as for compact case. See \cite{ELflow}.
Let $f\in C_K^\infty$.  Differentiate equation  $(~\ref{eq: integral kernel})$
on page ~\pageref{eq: integral kernel} to obtain:

$$d(P_t^hf)(v_0)=\int_M \langle \nabla p_t^h(-,y),v_0\rangle_{x_0}f(y)e^{2h}dy.$$

On the other hand, we may rewrite  equation $(~\ref{formula})$ as follows:

$$d(P_t^hf)(v_0)=\int_M  p_t^h(x_0,y)f(y)E\left\{{1\over t}
 \int_0^t \langle TF_s(v_0), X(x_s)dB_s\rangle | x_t=y\right\} e^{2h} dy$$
Comparing the last two equations, we get:
$$\nabla p_t^h(-,y)(x_0)=p_t^h(x_0,y)E\left\{{1\over t} \int_0^t TF_s^*(XdB_s)| x_t=y\right\}.$$ 
Thus finished the proof. \hfill \rule{3mm}{3mm}\bigskip

\noindent{\bf Remark:} Assume $\int_0^t \E|Y(T_xF_s)|^2ds<\infty$ for each
 $x\in M$. If
formula (~\ref{formula}) holds for $f\in C_K^\infty$, it holds for $f\in C_0^2(M)$,
the space of $C^2$ functions vanishing at infinity. 

\noindent Proof:
First assume $f$ positive. Take $f_n$ in $C_K^\infty$ converging to $f$. Then
$d(P_tf_n)\to d(P_tf)$ by Schauder type estimate as in the appendix. 
The convergence of the R.H.S. 
of the formula is also clear. Next take $f=f^+ -f^-$.

\section{For higher order forms and gradient Brownian systems}

Let $\alpha$ be a $p$ form, $\beta$ a $1$ form, the wedge of $\alpha$ and $\beta$ is a $p+1$ form defined as follows ( following the notations from \cite{AB-MA}):

$$(\alpha \wedge\beta)(v^1, \dots, v^{p+1})=\sum_{i=1}^{p+1} (-1)^{p+1-i}\beta(v^i)
\alpha(v^1,\dots, \widehat{v^i}, \dots, v^{p+1}).$$

\noindent 
 The symbol  $\hat { } \hskip 4pt  $    here means that the item below it is omitted.

The exterior differentiation of $\alpha$ is given  by:
\begin{equation}
d\alpha(v_1,\dots, v_{p+1})= \sum_{j=1}^{p+1} (-1)^{j+1} 
\nabla \alpha (v_j)(v_1, \dots, \widehat{v_j},\dots, v_{p+1}).
\end{equation}

\noindent
Let $\phi$ be a $q$ form. Let $v_0=(v_0^1,\dots, v_0^q)$, 
$v_t=(TF_t(v_0^1) ,\dots, TF_t(v_0^q))$.  Analogously to the case of 1-forms,
 we define:

\begin{equation}
 (\delta P_t)\phi(v_0)\index{$\delta P_t$}=E\phi(v_t)\chi_{t<\xi}.
\end{equation}

Then by the argument in proposition ~\ref{pr: Derivative 3}, $\delta P_t$
 is a $L^2$ semigroup if $\sup_x E\left(|T_xF_t|^{2q}\chi_{t<\xi}\right)<\infty$.

\begin{proposition} Let  $M$ be a complete Riemannian manifold. Let $(X,A)$ be 
 a complete gradient Brownian  system on $M$ with generator
 $\half \triangle^h$ and satisfying:
$$E\left(\sup_{s\le t} |TF_s|^q\right)<\infty.$$
Then
\begin{equation}
(\delta P_t)\psi=\heatsemi\psi
\label{eq: gradient 1}
\end{equation}
for bounded $q$-forms $\psi$.
\end{proposition}

The proof is as in $~\ref{th: Stflour}$. 
Note a similar result holds for a gradient Brownian system with a general drift.
In fact we could have a parallel discussion of the properties of $\delta P_t$ for
$q$ forms as in chapter 6. We should also point out that as for 1-forms, 
extra conditions on $TF_t$ will give $(~\ref{eq: gradient 1})$ for gradient
 systems, without assuming  nonexplosion(see corollary  ~\ref{co: homotopy 1} on
page ~\pageref {co: homotopy 1}). 

\bigskip

Given a $q$ form $\psi$, we define a $q$-$1$ form as follows:
 
\begin{equation}\begin{array}{ll}
\int_0^t\psi \circ d x_s(v_0)=& 	
	{1\over q}\int_0^t\psi\left(X(x_s)dB_s, TF_s(v_0^1), \dots ,TF_s(v_0^{q-1})\right)\\
	 &-{1\over 2}\int_0^t\delta^h \psi \left(TF_s(v_0^1), \dots ,TF_s(v_0^{q-1}\right)\, ds. 
\label{eq: h-forms 2}
\end{array}\end{equation}

\noindent
Here $v_0^q\in T_{x_0}M$, $i=1,2,\dots, q$, and $v_0=(v_0^1,\dots, v_0^q)$.

\begin{proposition} Let  $(X,A)$ be a complete gradient Brownian system on a 
complete Riemannian manifold $M$ with generator $\half \triangle^h$. Suppose
for  a closed  $q$ form $\psi\in D(\bar \triangle^{h,q})\cap L^\infty$,
 \begin{equation}
(\delta P_t)\psi=\heatsemi \psi,
\label{eq:chap8 21}
\end{equation}
and also for each $x$ in $M$,
$$E\int_0^t|T_xF_s|^{2q} ds <\infty.$$
Then:

\begin{equation}
P_t^{h,q}\psi(x_0)= (-1)^{q+1} \frac 1t
 E\int_0^t\psi\circ dx_s \wedge\int_0^t<X(x_s)dB_s, TF_s(\cdot)>.
\label{eq:chap8 20}
\end{equation}

\end{proposition}

\noindent{\bf Proof:}    Let 
\begin{equation}
Q_t(\psi)(v_0)= -\frac12\int_0^t(\delta^{h}P_s^{h,q-1}\psi)(v_0)ds.
\end{equation}

Notice   $P_t^{h,q}(\psi)$ is smooth on $[0,T]\times M$ by elliptic 
regularity, so

\begin{eqnarray*}
\frac \partial {\partial t} Q_t(\psi)&=&-\frac 12 \delta^h(P_t^h\psi),\\
	d(Q_t(\psi))&=&-\frac 12\int_0^td\delta^h(P_s^{h,q}\psi)ds, \\
	\delta^h Q_t(\psi)&=&-\half \int_0^t \delta^h\delta^h (P_s^{h,q}\psi)\, ds
=0.
\end{eqnarray*}

\noindent
Moreover, 
	$$d(Q_t(\psi))=  \frac 12\int_0^t\triangle^{h,q}(P_s^{h,q}\psi)ds       
	=P_t^{h,q}\psi - \psi.$$

\noindent
since $\triangle^{h,q} \psi= -d\delta^h \psi$.
Therefore:

  $$\triangle ^{h,q-1} (Q_t(\psi))
=-P_t^{h, q-1}(\delta^h\psi) +\delta^h\psi.$$
 
\noindent
Next we apply It\^o formula  (see page \pageref{eq: Ito formula for gradient}) 
to $(t,v)\vert \to Q_{T-t}(\psi)(v)$: 
\[\begin{array}{ll}
	Q_{T-t}\psi(v_t)=& Q_T\psi(v_0)
	+\int_0^t\nabla Q_{T-s}\psi(X(x_s)dB_s)(v_s)    \\
	&+      \int_0^tQ_{T-s}\psi((d\wedge)^{q-1}(\nabla X(\cdot)dB_s))v_s)   \\
	&+\frac 12\int_0^t\triangle^h Q_{T-s}\psi(v_s)ds
	+\int_0^t\frac \partial{\partial s}(Q_{T-s})\psi(v_s)ds.
\end{array}\]

\noindent
From  the  calculations above  we get:  

\begin{eqnarray*}
	Q_{T-t}\psi(v_t)&=&Q_T\psi(v_0)
	+\int_0^t\nabla Q_{T-s}\psi(X(x_s)dB_s)(v_s)     \\
	&+&  \int_0^tQ_{T-s}\psi((d\wedge)^{q-1}(\nabla X(\cdot)dB_s))v_s)   \\
	&+&   \frac 12 \int_0^t \delta ^h \psi(v_s) ds.
\end{eqnarray*}

\noindent
Let $t=T$,  then $Q_{T-t}(\psi)=0$. By definition and the equality above, we have:

\begin{equation}\begin{array}{ll}
\int_0^T\psi\circ dx_s(v_s)=&Q_T\psi(v_0)+ \frac 1q\int_0^T\psi(X(x_s)dB_s, v_s) \\
	 &+ \int_0^T\nabla Q_{T-s}\psi(X(x_s)dB_s)(v_s)    \\ 
	 &+\int_0^T Q_{T-s}\psi((d\wedge)^{q-1}(\nabla X(\cdot)dB_s))v_s).     
\label{eq: formulae 2}
\end{array}\end{equation}

We will calculate the expectation of each term of $\int_0^t \psi\circ dx_s$ after
 wedging with $\int_0^T<X(x_s)dB_s, TF_s(\cdot)>ds$.  It turns out that the first term and the last term vanishes. The latter is from equation ~\ref{eq: gradient}  for
a gradient system  on  page ~\pageref{eq: gradient}. 

Take $v_0=(v_0^1, \dots, v_0^q)$, write $ v_s^i=TF_s(v_0^i)$. 

 Denote by $w_s(\cdot)$ the linear map:

 $$w_s(\cdot)=\overbrace{ (TF_s(\cdot), \dots, TF_s(\cdot))}^{q-1}.$$ 
Then
\begin{eqnarray*}
	&E&  \int_0^T\psi(X(x_s)dB_s, w_s(\cdot)) 
		\wedge\int_0^T<X(x_s)dB_s, TF_s(\cdot)>(v_0)        \\
	&=& \sum_{i=1}^q (-1)^{q-i} E\int_0^T\psi(v_s^i, v_s^1, \dots, 
\widehat{v_s^i} , \dots , v_s^q) \, ds  \\
	&=& \sum_{i=1}^q (-1)^{q-i}(-1)^{i-1} E\int_0^T \psi(v_s^1,\dots , v_s^q)\,ds  \\
	&=& q(-1)^{q-1} E\int_0^T\psi(v_s^1, \dots, v_s^q)ds       \\
	&=&q(-1)^{q+1} \int_0^T P_s^{h} \psi (v)ds. 
\end{eqnarray*}

\noindent
The last step uses the assumption: $\int_0^t E|T_xF_s|^{2q}ds <\infty$. 
Similar calculation shows:

\begin{eqnarray*}
	 &&E\{\int_0^T\nabla Q_{T-s}\psi(X(x_s)dB_s)(w_s(\cdot))
		\wedge \int_0^T<X(x_s)dB_s, TF_s(\cdot)>\}(v)       \\
&=& \sum_{i=1}^q (-1)^{q-i} E\int_0^T \nabla (Q_{T-s}\psi)(v_s^i)
		 (v_s^1 , \dots  \hat{v_s^i} , \dots, v_s^q) \, ds     \\
	&=& (-1)^{q+1}E\int_0^T (d(Q_{T-s}\psi))(v_s^1, \dots, v_s^q)\, ds, \\
&=&(-1)^{q+1}  \int_0^T P_s^{h} \left(P_{T-s}^h(\psi)-\psi\right)\, (v)\, ds  \\
	&=& (-1)^{q+1} \left[T(P_T^{h}\psi)(v)  - \int_0^T P_s^h\psi(v)ds\right].
\end{eqnarray*}

\noindent
Comparing these with $~\ref{eq: formulae 2}$, we have: 

$$ P_T^{h,q}\psi =(-1)^{q+1}{1\over T} E\int_0^T \psi\circ dx_s
	\wedge \int_0^T<X(x_s)dB_s, TF_s(\cdot)>.$$

\noindent
End of the proof. \hfill \rule{3mm}{3mm}

Note: With an additional condition: $\sup_{x\in M}E|T_xF_s|^{2q}<\infty$,
the formula in the above proposition holds for forms which is not necessarily 
bounded. See the remark on page ~\pageref{remark}.
\bigskip

\begin{corollary}
 Let  $\psi=d\phi$  be a $C^2$ $q$ form, then

\begin{equation}\begin{array}{ll}
P_t^{h,q}(d\phi)=&(-1)^{q+1}\frac 1t
E\{\phi(\overbrace{TF_t(\cdot),\dots TF_t(\cdot)}^{q-1})\\
  &\wedge\int_0^t<X(x_s)dB_s, TF_s(\cdot)>\}.
\end{array}\end{equation}
\end{corollary}\

\noindent
{\bf Proof:}
  In this case,
\begin{eqnarray*}
\int_0^t\psi \circ dx_s
&=&\frac 1q \int_0^t d\phi(X(x_s)dB_s,
 \overbrace{TF_s(\cdot), \dots, TF_s(\cdot)}^{q-1})\\
& & +\frac 12 \int_0^t \triangle^h\phi(\overbrace{TF_s(\cdot), \dots, TF_s(\cdot}^{q-1})	 )\, ds.
\end{eqnarray*}

\noindent
There is  also the following equality:

\[\begin{array}{l}
	E\int_0^t d\phi(X(x_s)dB_s, \, TF_s(\cdot),\dots, TF_s(\cdot) )  
		\wedge \int_0^t <X(x_s)dB_s,\, TF_s(\cdot)>      \\
	=q E\int_0^t\nabla\phi\left (X(x_s)dB_s\right)TF_s(\cdot),\dots, 
TF_s(\cdot))\wedge \int_0^t <X(x_s)dB_s, TF_s(\cdot)>.\\
\end{array}\]

But by It\^o formula,

\[\begin{array}{ll}
\phi(v_t)=&\phi(v_0)+\int_0^t \nabla \phi(XdB_s)(v_s) \\
&+\half\int_0^t\triangle^h\phi(v_s)ds.\\
\end{array}\]
Here $v_t$ is the q-1 vector induced by $TF_t$ as is defined in the beginning of the section.
So

\begin{eqnarray*}
&&\int_0^t d\phi\circ dx_s\wedge\int_0^t <XdB_s, TF_s(\cdot)>\\
&=&E\{\phi(\overbrace{TF_t(\cdot),\dots TF_t(\cdot)}^{q-1})
  \wedge\int_0^t<X(x_s)dB_s, TF_s(\cdot)>\}. \end{eqnarray*}
End of the proof. \hfill \rule{3mm}{3mm}
\bigskip

This corollary can also be proved directly as in theorem ~\ref{th:elementary}. 
The factor ${1\over q}$ in the formula  may look odd, but it is due to that the
 tensors concerned is not symmetric.

\chapter{Appendix}

\begin{lemma}\cite{BA86}
Let $M$ be a complete Riemannian manifold. There exists an increasing sequence $\{h_n\}\subset C_K^\infty$ such that:

(1). $0\le h_n\le 1$

(2). $\lim_{n\to \infty} h_n(x)=1$, each $x$.

(3). $|\nabla h_n|\le {1\over n}$, for all $n$.
\end{lemma}

\noindent{\bf Proof:} This is standard result. Here's a proof from \cite{BA86}.
Let $M=R^1$. We may construct such a sequence $\{f_n\}$ as is well known. For a 
complete 
Riemannian manifold, there is a $C^\infty$  smooth function $f$ on $M$ such that
 $|\nabla f|\le 1$ and $\{|f|\le k\}$ is compact for all numbers $k$. 
Let $h_n=f_n\circ f$. Then $h_n$ is an increasing sequence and satisfies the
 requirements.

\subsection*{Schauder type estimates}

Let $M$ be a smooth differential manifold. Let $L$ be an elliptic differential 
operator on $M$. In local coordinates, we may write:
$$L=\sum_{i,j} a_{i,j}(x){\partial^2\over \partial x^i\partial x^j} +
\sum_i b_i(x){\partial \over \partial x^i} +c.$$
Here $(a_{i,j}(x))$ is symmetric positive definite $n\times n$ matrix for each $x$.
The coefficients are assumed to be $C^2$.

Consider the following differential equation of parabolic type:
\begin{equation}
Lu={\partial u\over \partial t}
\label{eq: parabolic}
\end{equation}
on a domain $\Omega\subset R^{n+1}$. A solution $u$ to equation
$(~\ref{eq: parabolic})$ is a function which is jointly continuous in $(t,x)$, 
$C^2$ in $x$ and $C^1$ in $t$ for $t>0$. A function which  satisfies the 
above regularity will be said to be in $C^{2,1}$. 

Let $D$ be a set of $R^{n+1}$, define:

\begin{equation}
|u|_{2,1}^D=\max_D(|u|+|D_xu|+|D_x^2u|+|D_tu|).
\end{equation}

\begin{theorem}
Let $B$ be a bounded domain in $R^n$, $D=B\times (0,1)$,  and $u$ a $C^{2,1}$
 solution to the parabolic equation. Let $D_1$ be a subdomain of $D$ with 
$d(\partial D_1,\partial D)$ denoting the distance from the boundary of 
$D_1$ to the boundary of $ D$.
Then
$$|u|_{2,1}^{D_1} \le k\max_D |u|.$$
Here $k$ is a constant depending only on the bounds of $a_{i,j}$, $da_{i,j}$,
 $db_i$, $dc$ in $D$ and   $d(\partial D_1,\partial D)$.
\end{theorem}

 See \cite{Friedman} $P_{64}$ for reference. 
By standard argument in analysis on uniform convergence, we conclude
(see \cite{AZ74} and $P_{89}$ of \cite{Friedman}): 

Let $u_n$ be an increasing sequence of (locally bounded $C^2$) solutions of the parabolic equation 
$(~\ref{eq: parabolic})$. Assume $\lim_{n\to \infty}u_n=u$ pointwise in $D$. 
Then  $u$ is also a solution and the convergence is  in fact uniformly in $C^1$
on compact subset  $D_1$ of $D$. Thus $D_x u_n\to D_xu$, and $D_tu_n\to D_tu$ in 
 $D_1$.

%\begin{theindex}
%\input thesis.index
%\end{theindex}
%

\printindex

\end{document}

%% file: picture.1
     %\documentstyle[12pt]{report}
      	%\newcommand{\A}{{\bf \cal A}}
	%\newcommand{\half}{{{1\over 2}}}
	%\begin{document}

\section*{Flow chart}

\begin{picture}(340,930)(-10,0)

\put(97,925){{\footnotesize various criteria}}
\put(137,920){\vector(0,-1){35}}
\put(139,905){{\tiny chapter 4}} 
\put(280,900){strong 1-complete}
\put(0,870){$E\sup_{s\le t}|TF_s|<\infty$ + complete}
\put(190,870){\vector(4,1){80}}
\put(190,870){\vector(4,-1){80}}
\put(210,850){{\tiny chap 6}}
\put(210,887){{\tiny chap 5}}
\put(280,840){$dP_tf=\delta P_t(df)$}

\put(220,825){{\tiny $\A=\half \triangle^h$}}
\put(260,808){\vector(-2,1){107}}
\put(200,790){$dP_t=\delta P_t$ and $E|T_xF_t|\chi_{t<\xi}<\infty$}
\put(80,770){\vector(3,4){68}}
\put(85,810){{\tiny chap 7}}
\put(0,750)
{Gradient BM+$\sup_{s\le t}\sup_{x\in M}E\left(|TF_s|^{2+\delta}\chi_{t<\xi}\right)
<\infty$}

\put(0,700){$\int_0^\infty \sup_{x\in K} E(|TF_s|\chi_{s<\xi})ds<\infty$}
\put(167,700){+$dP_t=\delta P_t$}
\put(232,700){\vector(1,0){50}}
\put(239,707){{\tiny chap 7}}
\put(293,700){finite h-volume}

\put(0,635){\shortstack{$\int_0^\infty \sup_{x\in K} E|TF_s|ds <\infty$\\
$\sup_{x\in K}E|TF_s|^{1+\delta}<\infty$}}
\put(135,645){+ strong 1-complete}
\put(246,645){\vector(1,0){50}}
\put(253,652){{\tiny Chap 7}}
\put(300,645){$\pi_1(M)=\{0\}$}

\put(0,580){$E(\sup_{s\le t} |TF_s|^{p+\delta})<\infty$ + complete}
\put(192,580){\vector(1,0){50}}
\put(202,587){{\tiny Chap 5}}
\put(260,580){strong p-completeness}
\put(295, 575){\vector(0,-1){65}}
\put(297,550){{\tiny chap 7}}
\put(170,530){{\footnotesize  + strong $p^{th}$ moment stable}}
\put(230, 500){bounded  closed    p-forms are exact}

\put(0,442){ \shortstack{$\int_0^t E|TF_s|^2ds<\infty$\\
 $dP_tf=\delta P_t(df)$ \\complete}}
\put(100,455){\vector(1,0){50}}
\put(112,457){{\tiny chapter 8}}
\put(155,455){\shortstack{$\nabla \log p_t^h(\cdot,y_0)(x_0)=E\left\{{1\over t}
 \int_0^t (TF_s)^* (XdB_s)| x_t=y_0\right\} $}}

\put(0,400){{\scriptsize The arrows here denote implications in the indicated 
chapter.}}
\end{picture}

   % \enddocument